\documentclass[12pt,twoside]{amsart}
\usepackage{amssymb,amsmath,amsfonts,latexsym,mathtools,tikz,tikz-cd}
\usepackage{amscd}
\usepackage{xcolor}
\usepackage{array,multirow,tabularx,longtable}
\usepackage[normalem]{ulem}
\usepackage{comment}
\usepackage{hyperref}
\usepackage[abbrev,alphabetic]{amsrefs}
\usepackage{url}
\usepackage[margin=1.25in]{geometry}

\usepackage{ascmac} 

\newcommand{\cred}{\color{black}}

\title[Fano threefolds in positive characteristic III]
{Fano threefolds in positive characteristic III} 
\author{Masaya Asai and Hiromu Tanaka} 
\address{Department of Mathematics, 
Graduate School of Science, 
Kyoto University, 
Kyoto 606-8502, JAPAN} 
\email{tanaka.hiromu.7z@kyoto-u.ac.jp}
\email{asaim.math@gmail.com}
\newcommand{\pr}[0]{{\operatorname{pr}}}
\newcommand{\et}[0]{{\operatorname{\acute{e}t}}}
\newcommand{\Bl}[0]{{\operatorname{Bl}}}

\newcommand{\red}[0]{{\operatorname{red}}}
\newcommand{\Ker}[0]{{\operatorname{Ker}}}

\newcommand{\Br}[0]{{\operatorname{Br}}}

\newcommand{\Tr}[0]{{\operatorname{Tr}}}

\newcommand{\Spec}[0]{{\operatorname{Spec}}}

\newcommand{\Supp}[0]{{\operatorname{Supp}}}
\newcommand{\Pic}[0]{{\operatorname{Pic}}}

\newcommand{\Ex}[0]{{\operatorname{Ex}}}
\newcommand{\Ext}[0]{{\operatorname{Ext}}}
\newcommand{\NE}[0]{{\operatorname{NE}}}
\newcommand{\Coker}[0]{{\operatorname{Coker}}}
\renewcommand{\Im}[0]{{\operatorname{Im}}}
\newtheorem{thm}{Theorem}[section]
\newtheorem{lem}[thm]{Lemma}
\newtheorem{cor}[thm]{Corollary}
\newtheorem{prop}[thm]{Proposition}

\newtheorem{step}{Step}

\theoremstyle{definition}

\newtheorem{dfn}[thm]{Definition}

\newtheorem{rem}[thm]{Remark}
      
\newtheorem{nota}[thm]{Notation}         

\newtheorem*{clm}{Claim}
\newtheorem{claim}[thm]{Claim}

\makeatletter
  
  \@addtoreset{equation}{thm}
\makeatother

\newenvironment{talign*}
 {\csname align*\endcsname}
 {\endalign}

\newcommand{\MO}{\mathcal{O}}
\newcommand{\R}{\mathbb{R}}
\newcommand{\Q}{\mathbb{Q}}
\newcommand{\Z}{\mathbb{Z}}
\newcommand{\G}{\mathbb{G}}

\newcommand{\F}{\mathbb{F}}
\renewcommand{\P}{\mathbb{P}}

\newcommand{\wt}{\widetilde}
\begin{document}

\maketitle

\begin{abstract}
We classify primitive Fano threefolds in positive characteristic whose Picard numbers are at least two.  
We also classify Fano theefolds of Picard rank two. 
\end{abstract}

\tableofcontents

\section{Introduction}

This article  is the third part of our series of papers. 
In the first and second parts \cite{TanI}, \cite{TanII}, we studied Fano threefolds $X$ in positive characteristic with $\rho(X)=1$, {\cred where $\rho(X)$ denotes the Picard number of $X$}. 
In this paper, we classify primitive Fano threefolds $X$ in positive characteristic with $\rho(X) \geq 2$ 
{\cred (for the definition of primitive Fano threefolds, see Subsection \ref{ss-notation}(6))}. 
More precisely, the main theorem is as follows. 

%

\begin{thm}[Theorem \ref{t-pic2-C}, Theorem \ref{t-pic2-D}, Theorem \ref{t-pic2-E}, 
Theorem \ref{t-pic3}]\label{t-main} 
{\cred Let $k$ be an algebraically closed field of characteristic $p>0$.} 
  Let $X$ be a primitive Fano threefold {\cred over $k$} with $\rho(X)\ge 2$.
  Then $X$ is isomorphic to one of the following varieties.
    \begin{center}
      \begin{longtable}{cccp{10.6cm}c}
No. &      $\rho(X)$ & $(-K_X)^3$ & Description  \\ \hline
2-2 &     $2$ & $6$ & a split double cover of $\mathbb{P}^2\times \mathbb{P}^1$ with $\mathcal L \simeq \MO_{\P^2 \times \P^1}(2, 1)$ \\ \hline
2-6 &      $2$ & $12$ & a smooth divisor on $\mathbb{P}^2\times \mathbb{P}^2$ of bidegree $(2,2)$ \\ \hline
2-6 &      $2$ & $12$ & a split double cover of $W$ with $\mathcal L^{\otimes 2} \simeq \omega_W^{-1}$ \\ \hline
2-8 &      $2$ & $14$ & a split double cover of $V_7=\mathbb{P}(\MO_{\mathbb{P}^2}\oplus \MO_{\mathbb{P}^2}(1))$ with $\mathcal L^{\otimes 2} \simeq \omega_{V_7}^{-1}$ \\ \hline
2-18 &      $2$ & $24$ & a split double cover of $\mathbb{P}^2\times \mathbb{P}^1$ 
      with $\mathcal L \simeq \MO_{\P^2 \times \P^1}(1, 1)$   \\ \hline
2-24 &      $2$ & $30$ & a smooth divisor on $\mathbb{P}^2\times \mathbb{P}^2$ of bidegree $(1,2)$ \\ \hline
2-32 &      $2$ & $48$ & $W$  \\ \hline
2-34 &      $2$ & $54$ & $\mathbb{P}^2\times \mathbb{P}^1$ \\ \hline
2-35 &      $2$ & $56$ & $V_7=\mathbb{P}(\MO_{\mathbb{P}^2}\oplus \MO_{\mathbb{P}^2}(1))$   \\ \hline
2-36 &      $2$ & $62$ & $\mathbb{P}(\MO_{\mathbb{P}^2}\oplus \MO_{\mathbb{P}^2}(2))$  \\ \hline
3-1 &      $3$ & $12$ & a split double cover of $\mathbb{P}^1\times \mathbb{P}^1\times \mathbb{P}^1$ 
      with $\mathcal L \simeq \MO_{\P^1 \times \P^1 \times \P^1}(1,1,1)$  \\ \hline
3-2 &      $3$ & $14$ & a smooth member of the complete linear system $|\MO_P(2)\otimes\pi^*\MO_{\mathbb{P}^1\times \mathbb{P}^1}(2,3)|$ on the $\mathbb{P}^2$-bundle $\pi\colon P=\mathbb{P}(\MO_{\mathbb{P}^1\times \mathbb{P}^1}\oplus \MO_{\mathbb{P}^1\times \mathbb{P}^1}(-1,-1)^{\oplus 2})\to \mathbb{P}^1\times \mathbb{P}^1$  \\ \hline
3-27 &      $3$ & $48$ & $\mathbb{P}^1\times \mathbb{P}^1\times \mathbb{P}^1$  \\ \hline
3-31 &      $3$ & $52$ & $\mathbb{P}(\MO_{\mathbb{P}^1\times \mathbb{P}^1}\oplus \MO_{\mathbb{P}^1\times \mathbb{P}^1}(1,1))$  
      \end{longtable}
  \end{center} 
  In the above table, we use the following notation and terminologies. 
  \begin{enumerate}
\item We say that $f : X \to Y$ is a split double cover if 
$f$ is a finite surjective morphism of projective normal varieites such that $\MO_Y \to f_*\MO_X$ splits as an $\MO_Y$-module homomorphism and the induced field extension $K(X) \supset K(Y)$ is of degree two. 
\item For a split double cover $f : X \to Y$, we set $\mathcal L := (f_*\MO_X/\MO_Y)^{-1}$, which is an invertible sheaf on $Y$ (Remark \ref{r-L-inv}). 
\item $W$ is a smooth prime divisor on $\P^2 \times \P^2$ of bidegree $(1, 1)$. 
Note that such a threefold is unique up to isomorphisms (Lemma \ref{l-bideg11}). 
  \end{enumerate}
\end{thm}

We also classify Fano threefolds with $\rho(X)=2$. 
Although this case has been treated already in \cite{Sai03}, 
their strategies are different.

\begin{thm}[Section \ref{s-pic2}]\label{t-main2} 
{\cred Let $k$ be an algebraically closed field of characteristic $p>0$.} 
  Let $X$ be a Fano threefold {\cred over $k$} with $\rho(X) =2$.
  Then $X$ is isomorphic to one of the following varieties. 
  \begin{center}
\begin{longtable}{ccp{10cm}c}
No. & $(-K_X)^3$ & Description  & Extremal rays\\ \hline
2-1 & $4$ & blowup of $V_1$ along an elliptic curve of degree $1$ & $D_1+E_1$\\ \hline
2-2 & $6$ & a split double cover of $\mathbb{P}^2\times \mathbb{P}^1$ with $\mathcal L \simeq \MO(2, 1)$ & $C_1+D_1$\\ \hline
2-3 & $8$ & blowup of $V_2$ along an elliptic curve of degree $2$ & $D_1+E_1$\\ \hline
2-4 & $10$ & blowup of $\P^3$ along a curve of genus $10$ {\cred and} degree $9$ & $D_1+E_1$\\ \hline
2-5 & $12$ & blowup of $V_3$ along an elliptic curve of degree $3$ & $D_1+E_1$\\ \hline
2-6 & $12$ & 
a smooth divisor on $\mathbb{P}^2\times \mathbb{P}^2$ of bidegree $(2,2)$, or 
a split double cover of $W$ with $\mathcal L^{\otimes 2} \simeq \omega_W^{-1}$ & $C_1+C_1$\\ \hline
2-7 & $14$ & blowup of $Q$ along a curve of genus $5$ {\cred and} degree $8$ & $D_1+E_1$\\ \hline
2-8 & $14$ &  a split double cover of $V_7=\mathbb{P}(\MO_{\mathbb{P}^2}\oplus \MO_{\mathbb{P}^2}(1))$ 
 with $\mathcal L^{\otimes 2} \simeq \omega_{V_7}^{-1}$ & $C_1+E_3\,{\rm or}\,E_4$\\ \hline
2-9 & $16$ & blowup of $\P^3$ along a curve of genus $5$ and degree $7$ & $C_1+E_1$\\ \hline
2-10 & $16$ & blowup of $V_4$ along an elliptic curve of degree $4$ & $D_1+E_1$\\ \hline
2-11 & $18$ & blowup of $V_3$ along a line & $C_1+E_1$\\ \hline
2-12 & $20$ & blowup of $\P^3$ along a curve of genus $3$ and degree $6$ & $E_1+E_1$\\ \hline
2-13 & $20$ & blowup of $Q$ along a curve of genus $2$ and degree $6$ & $C_1+E_1$\\ \hline
2-14 & $20$ & blowup of $V_5$ along an elliptic curve of degree $5$ & $D_1+E_1$\\ \hline
2-15 & $22$ & blowup of $\P^3$ along a curve of genus ${\cred 4}$ and degree $6$ & $E_1+E_3\,{\rm or}\,E_4$\\ \hline
2-16 & $22$ & blowup of $V_4$ along a conic & $C_1+E_1$\\ \hline
2-17 & $24$ & blowup of $\P^3$ along an elliptic curve curve of degree $5$ & $E_1+E_1$\\ \hline
2-18 & $24$ & a split double cover of $\mathbb{P}^2\times \mathbb{P}^1$ 
      with $\mathcal L \simeq \MO_{\P^2 \times \P^1}(1, 1)$ & $C_1+D_2$\\ \hline
2-19 & $26$ & blowup of $V_4$ along a line & $E_1+E_1$\\ \hline
2-20 & $26$ & blowup of $V_5$ along a cubic rational curve & $C_1+E_1$\\ \hline
2-21 & $28$ & blowup of $Q$ along a rational curve of degree $4$ & $E_1+E_1$\\ \hline
2-22 & $30$ & blowup of $\P^3$ along a rational curve of degree $4$ & $E_1+E_1$\\ \hline
2-23 & $30$ & blowup of $Q$ along an elliptic curve of degree $4$ & $E_1+E_3\,{\rm or}\,E_4$\\ \hline
2-24 & $30$ & a smooth divisor on $\mathbb{P}^2\times \mathbb{P}^2$ of bidegree $(1,2)$ & $C_1+C_2$\\ \hline
2-25 & $32$ & blowup of $\P^3$ along an elliptic curve of degree $4$ & $D_2+E_1$\\ \hline
2-26 & $34$ & blowup of $Q$ along a cubic rational curve  & $E_1+E_1$\\ \hline
2-27 & $38$ & blowup of $\P^3$ along a cubic rational curve & $C_2+E_1$\\ \hline
2-28 & $40$ & blowup of $\P^3$ along an elliptic curve of degree $3$ & $E_1+E_5$\\ \hline
2-29 & $40$ & blowup of $Q$ along a conic & $D_2+E_1$\\ \hline
2-30 & $46$ & blowup of $\P^3$ along a conic & $E_1+E_2$\\ \hline
2-31 & $46$ & blowup of $Q$ along a line & $C_2+E_1$\\ \hline
2-32 & $48$ & $W$ & $C_2+C_2$\\ \hline
2-33 & $54$ & blowup of $\P^3$ along a line & $D_3+E_1$\\ \hline
2-34 & $54$ & $\P^2 \times \P^1$ & $C_2+D_3$\\ \hline
2-35 & $56$ & $V_7=\mathbb{P}(\MO_{\mathbb{P}^2}\oplus \MO_{\mathbb{P}^2}(1))$ & $C_2+E_2$\\ \hline
2-36 & $62$ & $\mathbb{P}(\MO_{\mathbb{P}^2}\oplus \MO_{\mathbb{P}^2}(2))$ & $C_2+E_5$\\ \hline
      \end{longtable}
  \end{center} 
  In the above table, we use the following notation and terminologies in addition to (1)--(3) from Theorem \ref{t-main}. 
  \begin{enumerate}
  \setcounter{enumi}{3}
\item $Q$ is the smooth quadric hypersurace in $\P^4$. 
\item For $1 \leq d \leq 7$ with $d \neq 6$, 
let $V_d$ be a Fano threefold of index $2$ satisfying $(-K_{V_d}/2)^3 =d$. 
Note that $V_d$ is not determined uniquely by $d$. 
\item The centres of all the blowups are smooth curves. 
The degree of a curve $B$ on a Fano threefold $Y$ is defined as $-\frac{1}{r_Y}K_Y \cdot B$, 
where $r_Y$ denotes the index of $Y$. 
A line (resp. conic) is a smooth rational curve of degree one (resp. two). 
Note that if $|-\frac{1}{r_Y}K_Y|$ is very ample, then a curve of degree one (resp. two) is automatically a line (resp. conic). 
  \end{enumerate}

 \end{thm}

For more details, see Table \ref{table-pic2} in Section \ref{s table rho2}. 
The above theorems are positive-characteristic analogues of Mori--Mukai's results in characteristic zero \cite{MM81}, \cite{MM83}. 
The resulting tables are identical to those of characteristic zero. 
Our main strategy is the same as in characteristic zero, 
although we need to overcome some  obstructions. 
For the primitive case, we shall follow the thesis by Ott \cite{Ott}, which is extensively detailed. 
We here point out some differences. 

(I) {\em Conic bundles}.  
In characteristic zero, the discriminant divisor $\Delta_f$ of a threefold conic bundle $f: X \to Y$ is reduced and normal crossing. 
However, both the conclusions fail in characteristic two. 
Because of this, we shall need some minor modifications of the proofs. 
For example, if {\cred $\rho(X)=3$,} $Y = \P^1 \times \P^1$, and $\Delta_f \sim \MO_{{\cred \P^1 \times \P^1}}(2, 0)$, then we immediately get a contradiction in characteristic zero, because $\Delta_f$ does not contain any smooth rational curve as a connected component {\cred \cite[Proposition 6.2(3) and Proposition 6.3(1)]{MM83}}.

(II) {\em Double covers}. 
In characteristic two, double covers are more complicated than those in characteristic $p \neq 2$. 
Fortunately, all the double covers $f: X \to Y$ appearing in this paper can be checked to be split, i.e., $\MO_Y \to f_*\MO_X$ splits. 
Although it is hard to  control general double covers in characteristic two, 
split double covers are quite similar to the ones in characteristic $p \neq 2$, e.g., 
they are explicitly given, $\kappa(X) \geq \kappa(Y)$, etc (Section \ref{s-dc}).  

{\cred 
\begin{rem}\label{r intro existence}
The authors do not know whether 
all classes of Fano threefolds listed in Theorem \ref{t-main2} actually exist in arbitrary characteristic (e.g., 
the existence of Fano threefolds of No.\ 2-8 is not clear when $p=2$). 
\end{rem}
}

\vspace{5mm}

\textbf{Acknowledgements:} 
This paper is based on the master thesis \cite{Asa} of the first author, in which 
 primitive Fano threefolds of odd characteristic are classified. 
{\cred The authors 
thank the referee for reading the manuscript carefully and for suggesting several  improvements.} 
The second author was funded by JSPS KAKENHI Grant numbers JP22H01112 and JP23K03028.

\section{Preliminaries}

\subsection{Notation}\label{ss-notation}

In this subsection, we summarise notation used in this paper. 

\begin{enumerate}
\item We will freely use the notation and terminology in \cite{Har77} and \cite{KM98}. 
In particular, $D_1 \sim D_2$ means linear equivalence of Weil divisors. 
\item 
Throughout this paper, 
we work over an algebraically closed field $k$ 
of characteristic $p>0$ unless otherwise specified. 
\item For an integral scheme $X$, 
we define the {\em function field} $K(X)$ of $X$ 
as the local ring $\MO_{X, \xi}$ at the generic point $\xi$ of $X$. 
For an integral domain $A$, $K(A)$ denotes the function field of $\Spec\,A$. 

\item 
For a scheme $X$, its {\em reduced structure} $X_{\red}$ 
is the reduced closed subscheme of $X$ such that the induced closed immersion 
$X_{\red} \to X$ is surjective. 
\item We say that $X$ is a {\em variety} (over $k$) if 
$X$ is an integral scheme which is separated and of finite type over $k$. 
We say that $X$ is a {\em curve} (resp. a {\em surface}, resp. a {\em threefold})  
if $X$ is a variety over $k$ of dimension one (resp. {\em two}, resp. {\em three}). 
\item 
 We say that $X$ is a {\em Fano threefold}  
if $X$ is a three-dimensional smooth projective variety over $k$ such that $-K_X$ is ample. 
A Fano threefold $X$ is {\em imprimitive} if 
there exists a Fano threefold $Y$ and a smooth curve $B$ on $Y$ such that 
$X$ is isomorphic to the blowup $\Bl_B Y$ of $Y$ along $B$. 
We say that a Fano threefold $X$ is {\em primitive} if $X$ is not imprimitive. 
\end{enumerate}

\subsection{Split double covers}

\begin{dfn}\label{d-double-cover}
\begin{enumerate}
\item 
We say that a morphism $f: X \to Y$ is a {\em double cover}  if $f$ is a finite  surjective morphism of normal varieties such that the induced field extension $K(Y) \subset K(X)$ is of degree two. 
\item 
We say that a double cover $f: X \to Y$ is {\em split} 
if $\MO_Y \to f_*\MO_X$ splits as an $\MO_Y$-module homomorphism. 
\end{enumerate}
\end{dfn}

\begin{rem}\label{r-L-inv}
Let $f : X \to Y$ be a double cover, where $X$ and $Y$ are smooth varieties. 
\begin{enumerate}
    \item It is known that $\Coker(\MO_Y \to f_*\MO_X) = f_*\MO_X/\MO_Y$ is an invertible sheaf on $Y$ \cite[Lemma A.1]{Kaw21}. 
    We shall frequently use its inverse $\mathcal L := ( f_*\MO_X/\MO_Y)^{-1}$. 
    We have the following exact sequence: 
    \[
    0 \to \MO_Y \to f_*\MO_X \to \mathcal L^{-1} \to 0. 
    \]
    \item If $p \neq 2$, then the branch divisor $D$ satisfies 
    $\MO_Y(D) \simeq \mathcal L^{\otimes 2}$. 
\end{enumerate}
\end{rem}

\begin{lem}\label{l-dc-omega}
Let $f: X \to Y$ be a double cover of smooth projective varieties. 
For the invertible sheaf $\mathcal L := ( f_*\MO_X/\MO_Y)^{-1}$ (Remark \ref{r-L-inv}), 
it holds that $\omega_X \simeq f^*(\omega_Y \otimes \mathcal L)$. 
\end{lem}

\begin{proof}
See \cite[Proposition 0.1.3]{CD89}. 
\end{proof}

The following lemma gives criteria for the splitting of $f$. 

{\cred 
\begin{lem}\label{l-split-criterion1}
Let $f: X \to Y$ be a double cover of smooth projective varieties. 
Set  $\mathcal L := ( f_*\MO_X/\MO_Y)^{-1}$, which is an invertible sheaf (Remark \ref{r-L-inv}). 
Assume that $p=2$ and $f$ is separable. 
Then $H^0(Y, \mathcal L) \neq 0$. 
\end{lem}
}

\begin{proof}
Note that the trace map $\Tr : f_*\MO_X \to \MO_Y$ is nonzero, 
because the field-theoretic trace map $\Tr : K(X) \to K(Y)$ 
is given by $\Tr(b) = b + \sigma(b)$ for the Galois involution $K(X) \to K(X)$, and hence $\Tr(b) \neq 0$ for $b \in K(X) \setminus K(Y)$. 
As the composite $\MO_Y$-module homomorphism $\MO_Y \hookrightarrow f_*\MO_X \xrightarrow{\Tr} \MO_Y$ is zero, 
we get a nonzero $\MO_Y$-module homomorphism $\mathcal L^{-1} \simeq f_*\MO_X/\MO_Y \to \MO_Y$. 
Therefore, $H^0(Y, \mathcal L) \neq 0$. 
\end{proof}

\begin{lem}\label{l-split-criterion}
Let $f: X \to Y$ be a double cover of smooth projective varieties. 
Set  $\mathcal L := ( f_*\MO_X/\MO_Y)^{-1}$, which is an invertible sheaf (Remark \ref{r-L-inv}). 
Assume that one of the following (1)--(4) holds. 
\begin{enumerate}
\item $p \neq 2$. 
    \item $H^1(X, \mathcal L)=0$. 
       \item {\cred $f$ is separable and $H^1(Y, \MO_Y(E)) =0$ for {\cred every}  effective Cartier divisor $E$ on $Y$.}
    \item {\cred $f$ is inseparable and $Y$ is $F$-split}, i.e., $\MO_Y \to F_*\MO_Y$ splits as an $\MO_Y$-module homomorphism.
\end{enumerate}
Then $f$ is a split double cover, i.e., $\MO_Y \to f_*\MO_X$ splits. 
\end{lem}

\begin{proof}
It is well known that $f$ is split when (1) holds. 
{\cred In what follows, we assume $p=2$.} 
Note that $\MO_Y \to f_*\MO_X$ splits if and only if 
\[
0 \to \MO_Y \to f_*\MO_X \to \mathcal L^{-1} \to 0 
\]
splits. Since the  extension class corresponding to this exact sequence is contained in $H^1(Y, \mathcal L)$, (2) implies that $f$ is split.
Note that (3) implies (2) by Lemma \ref{l-split-criterion1}.



Assume (4). 
{\cred Then} 
the absolute Frobenius morphism $F : Y \to Y$ of $Y$ factors through 
$f: X \to Y$: 
\[
F : Y \to X \xrightarrow{f} Y, 
\]
and hence the splitting of $\MO_Y \to F_*\MO_Y$ implies the splitting of $\MO_Y \to f_*\MO_X$. 
Therefore, {\cred (4) implies that $f$ is split}.
\qedhere

\end{proof}


\subsection{Brauer groups}

\begin{dfn}[\cite{CTS21}*{Definition 3.2.1}]\label{d-Brauer}
For a scheme $X$, 
the {\em Brauer group} of $X$ is defined by 
\[
\Br(X) := H^2_{\et}(X, \G_m),
\]
where $H^2_{\et}(-)$ denotes the \'etale cohomology. 
For a ring $R$,  we set 
\[
\Br(R) := \Br(\Spec\,R) =H^2_{\et}(\Spec\,R, \G_m). 
\]
\end{dfn}

\begin{prop}\label{p-Br-vanish}
The following hold. 
\begin{enumerate}
\item If $Y$ is a regular noetherian integral scheme, 
then we have an injective group homorphism $\Br(Y) \hookrightarrow \Br(K(Y))$. 
\item If $C$ is a smooth curve over $k$, then $\Br(C)=\Br(K(C))=0$. 
\item If $Y$ is a smooth projective rational variety over $k$, then $\Br(Y) =0$. 
\end{enumerate}
\end{prop}

\begin{proof}
The assertion (1) follows from \cite[Corollaire 1.10]{Gro95}. 
The assertion  (2) holds by (1) and \cite[Theorem 1.2.14]{CTS21}. 
The assertion (3) holds by \cite[Corollary 6.2.11]{CTS21}. 
\end{proof}



\begin{prop}\label{p-Br-Pn}
Let $f: X \to Y$ be a projective morphism of smooth varieties. 
Fix $r \in \Z_{>0}$. 
Assume that $\Br(Y)=0$ and {\cred every}  fibre of $f$ is isomorphic to $\mathbb P^r$. 
Then there exists a vector bundle $E$ of rank $r+1$ such that $X$ is isomorphic to $\mathbb P(E)$ over $Y$. 
\end{prop}

For the reader's convenience, we include a sketch of a proof. 
For more details, we refer to websites \cite{Mur} and \cite{Mus}, which the following proof is based on. 

\begin{proof}
{\cred Note that $f$ is flat \cite[Theorem 23.1]{Mat86}.} 
Fix a closed point $y \in Y$ and set $A :=\MO_{Y, y}$. 
For the completion $\widehat{A} = \widehat{\MO}_{Y, y}$, 
we have an isomorphism $X \times_Y \Spec\,\widehat{A} \simeq \P^r_{\widehat A}$ over $\widehat A$ \cite[Corollary 1.2.15]{Ser06}. 
By Artin's approximation theorem, we obtain an 
isomorphism 
$X \times_Y \Spec\,A^h \simeq \P^r_{A^h}$ over the henselisation $A^h = \MO_{Y, y}^h$. 
Therefore, there exists an \'etale surjective morphism $Y' \to Y$ such that $X \times_Y Y' \simeq \P^r_{Y'}$.

By the exact sequence 
\[
1 \to \G_m \to {\rm GL}_{r+1} \to {\rm PGL}_{r+1} \to 1
\]
on \'etale topology, we obtain the following exact sequence: 
\[
H^1_{\et}(Y, {\rm GL}_{r+1}) \xrightarrow{\alpha} H^1_{\et}(Y, {\rm PGL}_{r+1}) \to H^2_{\et}(Y, \G_m) = \Br(Y)=0. 
\]
This implies that the ${\rm PGL}_{r+1}$-torsor $f : X \to Y$ comes from 
${\rm GL}_{r+1}$-torsor. 
Since {\cred a} ${\rm GL}_{r+1}$-torsor is nothing but a projective space bundle $\P_Y(E)$, the assertion holds by the following facts. 
\begin{itemize}
    \item $H^1_{\et}(X, {\rm GL}_{r+1})$ can be considered as a set of (Zariski) locally free sheaf of rank $r+1$ \cite[Expos\'e XI, Corollaire 5.3]{SGA1II}. 
    \item $H^1_{\et}(Y, {\rm PGL}_{r+1})$ consists of the isomorphism classes of 
    $\pi : W \to Y$ such that $W$ is a flat proper morphism such that 
    $W \times_Y Y' \simeq \P^{r}_{Y'}$ for some \'etale surjective morphism $Y' \to Y$. 
    \item $\alpha(E) = \P_Y(E)$.   
\end{itemize}
\end{proof}

\section{Extremal rays and Mori fibre spaces}

\subsection{Types of extremal rays}


\begin{dfn}\label{d-cont}
Let $X$ be a smooth projective threefold. 
Let $R$ be a $K_X$-negative extremal ray of $\overline{\NE}(X)$. 
By \cite[(1.1.1)]{Kol91}, there exists a unique morphism $f\colon  X \to Y$, called the {\em contraction (morphism)} of $R$, 
to a projective normal variety $Y$ such that the following hold. 
  \begin{enumerate}
    \item $f_*\MO_X = \MO_Y$. 
    \item For {\cred every}  curve $C$ on $X$, $[C]\in R$ if and only if $f(C)$ is a point.
  \end{enumerate}
\end{dfn}

\begin{dfn}\label{d-length}
Let $X$ be a smooth projective threefold. 
Let $R$ be an extremal ray of $\overline{\NE}(X)$ 
and let $f: X \to Y$ be the contraction of $R$. 
We set
\[ \mu_R := \min \left\{ (-K_X\cdot C)_X \mid C \mathrm{~is~a~rational~curve~on~} X \mathrm{~with~} [C]\in R \right\}, \]
which is called {\em  the length of an extremal ray} $R$. 
We say that $\ell$ is an {\em extremal rational curve} 
if $\ell$ is a rational curve on $X$ such that 
$[\ell] \in R$ and $-K_X\cdot \ell =\mu_R$. 
\end{dfn}

\begin{dfn}\label{d-type-ext-ray}
Let $X$ be a smooth projective threefold. 
Let $R$ be a $K_X$-negative extremal ray of $\overline{\NE}(X)$ 
and let $f: X \to Y$ be the contraction of $R$. 
\begin{enumerate}
\item $R$ and $f$ are called of {\em type $C$} if $\dim Y = 2$. 
\begin{itemize}
\item $R$ and $f$ are called of {\em type $C_1$} if 
$f$ is not smooth. 
\item $R$ and $f$ are called of {\em type $C_2$} if 
$f$ is smooth. 
\end{itemize}
\item $R$ and $f$ are called of {\em type $D$} if $\dim Y = 1$. 
Let $X_K$ be the generic fibre of $f$, which is a regular projective surface over $K:=K(Y)$. 
\begin{itemize}
\item $R$ and $f$ are called of {\em type $D_1$} if $1 \leq K_{X_K}^2 \leq 7$. 
\item $R$ and $f$ are called of {\em type $D_2$} if $K_{X_K}^2 =8$. 
\item $R$ and $f$ are called of {\em type $D_3$} if $K_{X_K}^2 =9$. 
\end{itemize}
\item $R$ and $f$ are called of {\em type $E$} if $\dim Y =3$. 
Set $D:=\Ex(f)$. 
\begin{itemize}
\item $R$ and $f$ are called of {\em type $E_1$} if 
$f(D)$ is a curve. 
\item  $R$ and $f$ are called of {\em type $E_2$} if 
$f(D)$ is a point, $D \simeq \P^2$, and $\MO_X(D)|_D \simeq \MO_D(-1)$. 
\item $R$ and $f$ are called of {\em type $E_3$} if 
$f(D)$ is a point, $D$ is isomorphic to a smooth quadric surface $Q$ on $\P^3$, and $\MO_X(D)|_D \simeq \MO_{\P^3}(-1)|_Q$. 
\item $R$ and $f$ are called of {\em type $E_4$} if 
$f(D)$ is a point, $D$ is isomorphic to a singular quadric surface $Q$ on $\P^3$, and $\MO_X(D)|_D \simeq \MO_{\P^3}(-1)|_Q$. 
\item  $R$ and $f$ are called of {\em type $E_5$} if 
$f(D)$ is a point, $D \simeq \P^2$, and $\MO_X(D)|_D \simeq \MO_D(-2)$. 
\end{itemize}
\item $R$ and $f$ are called of {\em type $F$} if $\dim Y = 0$. 
\end{enumerate}
\end{dfn}

\begin{rem}\label{r-length}
Let $X$ be a smooth projective threefold. 
Let $R$ be a $K_X$-negative extremal ray of $\overline{\NE}(X)$ 
and let $f: X \to Y$ be the contraction of $R$. 
\begin{enumerate}
\item Assume that $R$ is of type $C$. 
Then the following hold, 
{\cred because $f^{-1}(y)$ is a conic on $\P^2$ and we have $-K_X \cdot f^{-1}(y) = 2$ for every closed point $y \in Y$}. 
\begin{itemize}
    \item $R$ is of type $C_1$ if and only if $\mu_R =1$. 
    \item $R$ is of type $C_2$ if and only if $\mu_R =2$. 
\end{itemize}
\item Assume that $R$ is of type $D$. 
Then the following hold \cite[Proposition 3.17]{TanII}. 
\begin{itemize}
    \item $R$ is of type $D_1$ if and only if $\mu_R =1$. 
        \item $R$ is of type $D_2$ if and only if $\mu_R =2$.
        \item $R$ is of type $D_3$ if and only if $\mu_R =3$. 
\end{itemize}
\item Assume that $R$ is of type $E$. 
Then the following hold \cite[Proposition 3.22]{TanII}. 
\begin{itemize}
\item If $R$ of type $E_1$, then $\mu_R=1$. 
\item If $R$ of type $E_2$, then $\mu_R=2$. 
\item If $R$ of type $E_3$, then $\mu_R=1$. 
\item If $R$ of type $E_4$, then $\mu_R=1$. 
\item If $R$ of type $E_5$, then $\mu_R=1$. 
\end{itemize}
\end{enumerate}
\end{rem}

\subsection{Del Pezzo fibrations}


\begin{prop}\label{p-dP-fibres}
Let $X$ be a smooth projective 
threefold and let $f : X \to Y$ be 
a contraction of a $K_X$-negative extremal ray of type $D$.  
For $K := K(Y)$ and its algebraic closure $\overline K := \overline{K(Y)}$, 
let $X_K$ and $X_{\overline K}$ be the generic fibre and the geometric generic fibre of $f$, respectively. 
Then the following hold. 
\begin{enumerate}
\item For a closed point $y \in Y$, the effective Cartier divisor $f^*y$ is a prime divisor, i.e., {\cred every}  scheme-theoretic fibre of $f$ is geometrically integral. 
\item $X_{\overline K}$ is a projective canonical del Pezzo surface, i.e., 
$X_{\overline K}$ is a projective normal surface 
such that ${\cred X}_{\overline K}$ has at worst canonical singularities and $-K_{{\cred X}_{\overline K}}$ is ample. 
\item $1 \leq K_{X_K}^2 =K_{X_{\overline K}}^2  \leq 9$. 
\item $\Pic X_{\overline K}$ is a finitely generated free abelian group. 
Furthermore, $\Pic X_K \simeq \Z$. 
\item 
Assume that $K_{X_K}^2  = 9$. 
Then the following hold. 
\begin{enumerate}
\item $X_K \simeq \mathbb P^2_K$. In particular, general fibres are $\mathbb P^2$. 
\item There exist Cartier divisors $D$ on $X$ and $E$ on $Y$ such that $-K_{X} \sim 3D+f^*E$. 
\end{enumerate}
\item 
Assume that $K_{X_K}^2  = 8$. 
Then the following hold. 
\begin{enumerate}
\item If $X_K$ is smooth over $K$ and $y$ is a general closed point,  then 
both $X_{\overline K}$ and $X_y$ are isomorphic to $\mathbb P^1 \times \mathbb P^1$. 
\item If $X_K$ is not smooth over $K$ and $y$ is a general closed point, 
then $p =2$ and both $X_{\overline K}$ and $X_y$ are isomorphic to $\mathbb P(1, 1, 2)$. 
\item There exist Cartier divisors $D$ on $X$ and $E$ on $Y$ such that $-K_{X} \sim 2D+f^*E$. 
\end{enumerate}

\end{enumerate}
\end{prop}

\begin{proof}
The assertions (1)--(4) follow from 
\cite[Proposition 3.5]{TanII}. 
As for (5) and (6), see 
\cite[Lemma 3.14 and Lemma 3.15]{TanII}. 
\qedhere




\end{proof}

\begin{rem}
By \cite[Proposition 14.7]{FS20}, there actually exists a del Pezzo fibration $f: X \to Y$ such that $X$ is smooth and the generic fibre is not smooth. 
\end{rem}

\begin{prop}\label{p-dP-fibres2}
Let $X$ be a smooth projective 
threefold with an extremal ray $R$ of type $D$. 
Let $f : X \to Y$ be the contraction of $R$. For $K := K(Y)$, 
$X_K$ denotes the generic fibre of $f$. 
\begin{enumerate}
\item If $K_{X_K}^2 =9$, then $f$ is a $\P^2$-bundle, i.e., 
there exists a locally free $F$ on $Y$ of rank $3$ such that 
$X$ is isomorphic to $\P_Y(F)$ over $Y$. 
In particular, {\cred every}  fibre is $\P^2$. 
\item If $K_{X_K}^2 =8$, then there exists a $\P^3$-bundle $\P(E)$ over $Y$ and a closed immersion $j:X \hookrightarrow \P(E)$ such that 
$j_y(X_y) \subset \P^2$ is a quadric surface on $\P^2_k$ for every point $y \in Y$. 
\end{enumerate}
\end{prop}

\begin{proof}
Let us show (1). 
By Proposition \ref{p-dP-fibres}(5), we have $-K_X \sim 3D + f^*E$ for some Cartier divisors $D$ on $X$ and $E$ on $Y$, respectively. 
We use Fujita's $\Delta$-genus: $\Delta(V, L) = \dim V + L^{\dim V} -h^0(V, L)$. 
It holds that 
\[
\Delta(X_{\overline K}, \MO_X(D)|_{X_{\overline K}}) = 2 + 
{\cred (\MO_X(D)|_{X_{\overline K}})}^{2} -h^0(X_{\overline K}, \MO_X(D)|_{X_{\overline K}})=0. 
\]
Then the assertion (1) follows from the same argument as in \cite[Corollary 5.4]{Fuj75}. 
Similarly, the assertion (2) holds by using Proposition \ref{p-dP-fibres}(6) (cf. \cite[Corollary 5.5]{Fuj75}). 
\end{proof}

\subsection{Contraction morphisms}

\begin{lem}\label{firstlem}
  Let $X$ be a Fano threefold. 
  Then the following hold. 
  \begin{enumerate}
\item $\Pic X \simeq \Z^{\oplus \rho(X)}$. 
    \item  $H^i(X, \MO_X)=0$ for all $i>0$. In particular, $\chi(X, \MO_X)=1$. 
    \item  $-K_X\cdot c_2(X) =24$.
    \item For {\cred every}  effective divisor $D$ on $X$,
          \[ D \cdot c_2(X) = 6\chi(D, \MO_D)+6\chi(D, \MO_X(D)|_D)-2D^3-(-K_X)^2\cdot D\]
  \end{enumerate} 
\end{lem}

\begin{proof}
The assertions (1) and (2)  follow from  
\cite[Theorem 2.4]{TanI} (cf. \cite[Corollary 3.7]{Kaw21}). 
The assertion (3) holds by $\chi(X, \MO_X)=1$  and $\chi(X, \MO_X)=\frac{1}{24}(-K_X\cdot c_2(X))$, 
{\cred where the latter one follows from the Riemann--Roch theorem}. 

Let us show (4). 
By $\chi(X, \MO_X)=1$, the Riemann--Roch theorem implies  
  \begin{align*}
    \chi(\MO_X(D))=&1+\frac{1}{12}((-K_X)^2 +c_2(X)) \cdot D  +\frac{1}{4}(-K_X\cdot D^2)+\frac{1}{6}D^3,\\
    \chi(\MO_X(-D))=&1-\frac{1}{12}((-K_X)^2 +c_2(X)) \cdot D+\frac{1}{4}(-K_X\cdot D^2)-\frac{1}{6}D^3.  
  \end{align*}
  Hence 
  \[\chi(\MO_X(D))-\chi(\MO_X(-D))=\frac{1}{6}((-K_X)^2 +c_2(X)) \cdot D+\frac{1}{3}D^3.\]
  On the other hand, 
  we have 
\[
\chi(X, \MO_X(-D)) = \chi(X, \MO_X)-\chi(D, \MO_D), \quad \chi(X, \MO_X(D)) = \chi(X, \MO_X) + \chi(D, \MO_D(D)), 
\]
where $\MO_D(D) := \MO_X(D)|_D$. 
Therefore, we obtain 
\[
\chi(X, \MO_X(D))-\chi(X, \MO_X(-D))=\chi(D, \MO_D(D))+\chi(D, \MO_D),
\]
which implies 
\[
\chi(D, \MO_D(D))+\chi(D, \MO_D) = \frac{1}{6}((-K_X)^2 +c_2(X)) \cdot D+\frac{1}{3}D^3, 
\]
as required. 
\end{proof}

\begin{prop}\label{p-cont-ex}
Let $X$ be a smooth projective threefold and let $f : X \to Y$ be the contraction of a $K_X$-negative extremal ray $R$ of $X$. 
Let $\ell$ be a curve with $[\ell ] \in R$. 
Then the following sequence 
\[ 
0 \longrightarrow \Pic(Y) \xrightarrow{f^*} \Pic(X) \xrightarrow{(-\cdot \ell)} \mathbb{Z}  \]
is exact,     where $(- \cdot \ell)(D):=D\cdot \ell$ for $D\in \Pic(X)$.
\end{prop}

\begin{proof}
See \cite[Proposition 3.12]{TanII}. 
\end{proof}

\begin{cor}\label{picYrk1}
  Let $X$ be a Fano threefold with $\rho(X)=2$ and let $R$ be an extremal ray of $\NE(X)$ with the corresponding contraction $f\colon X\to Y$. Then $\Pic\,Y \simeq \Z$. 
\end{cor}

\begin{proof}
The assertion follows from 
$\Pic(X) \simeq \Z^2$ (Lemma \ref{firstlem}) and Proposition \ref{p-cont-ex}. 
\end{proof}

\begin{prop}\label{p-cont-ex2}
Let $X$ be a smooth projective threefold and let $f : X \to Y$ be the contraction of a $K_X$-negative extremal ray $R$ of $X$. 
Then the following hold. 
\begin{enumerate}
\item If $\dim Y \neq 1$, then $R^if_*\MO_X=0$ for {\cred every}  $i>0$. 
\item If $X$ is a Fano threefold, then $R^if_*\MO_X=0$ for {\cred every}  $i>0$. 
Furthermore, $H^j(X, \MO_X) \simeq H^j(Y, \MO_Y)$ for {\cred every}  $j \geq 0$. 
\end{enumerate}
\end{prop}

\begin{proof}
(1) 
If $\dim X_y \leq 1$ for {\cred every}  closed point $y \in Y$, 
then $R^if_*\MO_X=0$ holds by \cite[Theorem 0.5]{Tan15}. 
We may assume that $\dim X_y \geq 2$ for some closed point $y \in Y$. 
In particular, $\dim Y =0$ or $\dim Y =3$, 
as we assume $\dim Y \neq 1$. 
If $\dim Y=0$, then we obtain $R^if_*\MO_X=0$ by Lemma \ref{firstlem}. 

Assume that $\dim Y=3$. 
By $\dim X_y \geq {\cred 2}$, 
$f(E)$ is a point for $E:=\Ex(f)$. 
Then we have an exact sequence for {\cred every}  $r \in \Z_{\geq 0}$: 
\[
0 \to \MO_X(-(r+1)E) \to \MO_X(-rE) \to \MO_X(-rE)|_E \to 0. 
\]
Since $\MO_X(-E)|_E$ is ample (Definition \ref{d-type-ext-ray}) and $E$ is a normal projective toric surface, we have that $H^i(E, \MO_X(-rE)|_E)=0$ for every $i>0$. 
Therefore, we obtain surjections: 
\[
R^if_*\MO_X(-mE)) \to \cdots \to R^if_*\MO_X(-E)) \to R^if_*\MO_X.
\]
By the Serre vanishing theorem, we have $R^if_*\MO_X(-mE)) =0$ for some $m>0$, as required.

(2) By (1), we may assume that $\dim Y=1$. 
Pick general closed point{\cred s} $Q_1, ..., Q_n$ on $Y$. 
Set $S_1 := f^*Q_1, ..., S_n := f^*Q_n$, and $D_Y:=Q_1+\cdots + Q_n$. 
By the Serre vanishing theorem 
  \begin{align*}
    H^q(Y,R^if_*(f^*\MO_Y(D_Y)))=H^q(Y,R^if_*\MO_X\otimes \MO_Y(D_Y)) =0.
  \end{align*}
  for all $q>0$ and $n\gg 0$. 
  Hence, by the Leray spectral sequence, we have 
  \[H^0(Y,(R^if_*\MO_X)\otimes \MO_Y(D_Y))\simeq H^i(X,f^*\MO_Y(D_Y)).\]
  Hence it is enough to show that $H^i(X,f^*\MO_Y(D_Y))=H^i(X,\MO_X(S_1+\cdots + S_n))=0$ 
  for every $n > 0$. 
  We have an exact sequence 
  \[
0\longrightarrow \MO_X(-\sum_{\ell=1}^n S_{\ell}) \longrightarrow \MO_X \longrightarrow \bigoplus_{\ell=1}^n \MO_{S_{\ell}} \longrightarrow 0.\]
Since $Q_1+ \cdots + Q_n \sim D'_Y$ for some divisor $D'_Y$ such that $Q_1, ..., Q_n \not\in \Supp\,D'_Y$, 
it holds that 
\[
0\longrightarrow \MO_X \longrightarrow \MO_X(f^*D_Y) \longrightarrow \bigoplus_{\ell=1}^n \MO_{S_{\ell}}\longrightarrow 0
\]
{\cred is exact.} 
Since each $S_i$ is a canonical del Pezzo surface (Proposition \ref{p-dP-fibres}(2)), we have $H^i(S_{\ell},\MO_{S_{\ell}})=0$ for every $\ell$ and every $i>0$. 
It holds that $H^i(X, \MO_X)=0$ for {\cred every}  $i>0$, 
which implies $H^i(X, \MO_X(f^*D_Y))=0$ for {\cred every}  $i>0$. 
\end{proof}

\begin{cor}\label{euler_cont}
  Let $X$ be a Fano threefold with $\rho(X) \geq 2$. 
  For an extremal ray $R$ of $\NE(X)$, let   $f\colon X\to Y$ be the contraction of $R$. 
  Let $D_Y$ be an effective Cartier divisor on $Y$ and set $D:=f^*D_Y$. 
  Then we have
  \[
  \chi(D, \MO_D)=1-\chi(Y, \MO_Y(-D_Y)) \quad \text{and}\quad 
  \chi(D, \MO_X(D)|_D)=\chi(Y, \MO_Y(D_Y))-1.
  \]
\end{cor}

\begin{proof}
Set $\MO_D(D) := \MO_X(D)|_D$. 
  We have 
  \begin{align*}
    \chi(D, \MO_D) &= \chi(X, \MO_X)-\chi(X, \MO_X(-D)), \\
    \chi(D, \MO_D(D)) &= \chi(X, \MO_X(D))-\chi(X, \MO_X).
  \end{align*}
Recall that $\chi(X, \MO_X)=1$  (Lemma~\ref{firstlem}).
  By the projection formula, 
  we have 
  \[R^if_*\MO_X(D)\simeq R^if_*(f^*\MO_Y(D_Y)) \simeq \MO_Y(D_Y)\otimes R^if_*\MO_X\]
  for all $i\ge 0$.
  We obtain $f_*\MO_X(D)\simeq \MO_Y(D_Y)$, and $R^if_*\MO_X(D)= 0$ for all $i>0$ (Proposition \ref{p-cont-ex2}).
  Therefore, by the Leray spectral sequence, the following holds for every $i \geq 0$: 
  \[H^i(X,\MO_X(D))\simeq H^i(Y,f_*\MO_X(D))\simeq H^i(Y,\MO_Y(D_Y)).\]
  Hence $\chi(X, \MO_X(D))=\chi(Y, \MO_Y(D_Y))$, which implies 
  \[
  \chi(D, \MO_D(D))=\chi(Y, \MO_Y(D_Y))-1.
  \]
  Similarly, we get  $\chi(D, \MO_D)=1-\chi(Y, \MO_Y(-D_Y))$.
\end{proof}

\begin{lem}\label{l-(-K)^3-drops}
Let $Y$ be a smooth projective threefold and 
let $f: X \to Y$ be a blowup along a smooth curve $B$ on $Y$. 
Assume that $-K_X$ is ample. 
Then $(-K_X)^3 < (-K_Y)^3$. 
\end{lem}

\begin{proof}
Set $E :=\Ex(f)$. 
Since $-K_X$ is ample, we have $(-K_X)^2 \cdot E >0$. 
This, together with \cite[Lemma 3.21]{TanII}, implies the following: 
\begin{enumerate}
\item $(-K_X)^3 = (-K_Y)^3 +2K_X \cdot B +2g(B) -2$. 
\item $-K_X \cdot B -2g(B)+2  = (-K_X)^2 \cdot E>0$
\end{enumerate}
By (1), it is enough to show $-K_X \cdot B > g(B)-1$. 
Since $-K_X$ is ample, this holds when $g(B) \leq 1$. 
Hence we may assume that $g(B) > 1$. 
In this case, the required inequality $-K_X \cdot B > g(B)-1$ holds by (2): 
$-K_X \cdot B > 2g(B) -2 > g(B)-1$.   
\end{proof}

\subsection{Conic bundles}\label{ss-CB}


In this subsection, we recall some terminologies and results on conic bundles. 
For more details, we refer to \cite{Tan-conic} and \cite[Subsection 3.3]{TanII}.

\begin{dfn}
We say that $f: X \to S$ is a {\em conic bundle} if 
$f : X \to S$ is a flat projective morphism of noetherian schemes 
such that $X_s$ is isomorphic to a conic on $\mathbb P^2_{\kappa(s)}$ 
{\cred for every $s \in S$}. 
\end{dfn}

If $X$ is a smooth projective threefold and $R$ is an extremal ray of type $C$, 
then its contraction $f: X \to Y$ is a conic bundle to a smooth projective surface $Y$. 
For the definition of $\Delta_f$, we refer to \cite[3.10]{TanII}. 

\begin{prop}\label{conic-embedding}
Let $f :X \to S$ be a conic bundle, where $X$ and $S$ are smooth varieties. 
Then the following hold. 
\begin{enumerate}
\item $f_*\MO_X = \MO_S$. 
\item $f_*\omega_X^{-1}$ is a locally free sheaf of rank $3$. 
\item $\omega_X^{-1}$ is very ample over $S$, and hence it defines a closed immersion 
$\iota: X \hookrightarrow  \mathbb P( f_*\omega_X^{-1})$ over $S$. 
\end{enumerate}
\end{prop}

\begin{proof}
See \cite[Lemma 2.5 and Proposition 2.7]{Tan-conic}. 
\qedhere 
\end{proof}








\begin{prop}\label{conic_disc-num}
Let $f : X \to S$ be a conic bundle, where $X$ is a smooth projective threefold and $S$ is a smooth projective surface. 
Then the following hold. 
\begin{enumerate}
\item $f_*K_X \sim -2S$. 
\item 
$-f_*(K_{X/S}^2) \equiv \Delta_f$. 
\item 
For every divisor $D$ on $S$, it holds that 
\[
K_X^2 \cdot f^*D = -4 K_S \cdot D -\Delta_f \cdot D. 
\]
\end{enumerate}
\end{prop}

\begin{proof}
See \cite[Proposition 3.11]{TanII}. 
\qedhere
\end{proof}

\section{Basic properties of primitive Fano threefolds}

\subsection{Extremal rays on primitive Fano threefolds}

\begin{thm}[Cone theorem]\label{cone}
  Let $X$ be a Fano threefold.
  Then 
  there exist finitely many extremal rays $R_1, ..., R_n$ of $\NE(X)$ such that 
  \[
  \NE(X) = R_1 + \cdots + R_n. 
  \]
\end{thm}

For each $1 \leq i \leq n$, there exists an extremal rational curve $\ell_i$ such that 
$R_i = \R_{\geq 0}[\ell_i]$ (Definition \ref{d-length}). 

\begin{proof}
  See \cite[Theorem 1.2]{Mor82}.
\end{proof}



\subsubsection{Type C}


\begin{thm}\label{t-wild-cb}
Let $X$ be a Fano threefold. 
Assume that there exists a conic bundle $f: X \to Y$ such that any  fibre of $f$ is not smooth. 
Then the following hold. 
\begin{enumerate}
\item $p=2$. 
\item One of the following holds. 
\begin{enumerate}
    \item 
    $Y \simeq \P^2$ and $X$ is isomorphic to 
    a prime divisor on $\P^2 \times \P^2$ of bidegree $(1, 2)$. 
    Furthermore, $\rho(X) =2$ and $X$ is primitive. 
    \item  
    $Y \simeq \P^1 \times \P^1$ and $X$ is isomorphic to 
    a prime divisor on $P:=\P_Y( \MO(0, 1) \oplus \MO(1, 0) \oplus \MO)$ 
    which is linearly equivalent to $\MO_P(1)^{\otimes 2}$, 
    where $\MO_P(1)$ denotes the tautological bundle with respect to the $\P^2$-bundle structure $P \to Y$. 
    Furthermore, $\rho(X)=3$ and $X$ is imprimitive. 
\end{enumerate}
\end{enumerate}
In particular, if $X$ is a primitive Fano threefold and there exists an extremal ray of type $C$ whose contraction $f: X \to Y$ has no smooth fibres, 
then {\rm (a)} of {\rm (2)} holds. 
\end{thm}

\begin{proof}
See \cite[Corollary 8 and Remark 10]{MS03}. 
\end{proof}

\begin{lem}\label{Y-rat}
Let $X$ be a Fano threefold with an extremal ray $R$ of type $C$. 
Let $f: X \to Y$ be the contraction of $R$. 
Then  $Y$ is a smooth rational surface.
\end{lem}

\begin{proof}
By Theorem \ref{t-wild-cb}, we may assume that $f$ is generically smooth. 
Note that $Y$ is a smooth projective surface \cite[Main Theorem 1.1, (1.1.3.1)]{Kol91}. 
Since the geometric generic  fibre $X_{\overline{\eta}}$ is smooth, $(X_{\overline{\eta}}, 0)$ is $F$-pure. 
Then it follows from \cite[Corollary 4.10(2)]{Eji19} that $-K_Y$ is big. 
In particular, the Kodaira dimension of $Y$ is negative, i.e., $Y$ is a smooth ruled surface. 
On the other hand, $Y$ is rationally chain connected, because so is $X$ \cite[Ch. V, Theorem 2.13]{Kol96}. 
Therefore, $Y$ is a smooth rational surface. 
\end{proof}

\begin{thm}\label{dimy=2}
Let $X$ be a Fano threefold 
and let $R$ be an extremal ray of type $C$.  
Let $f \colon X\to Y$ be the contraction of $R$. 
Then $f$ is a conic bundle and $Y$ is a smooth rational surface. 
Furthermore, if $f$ is generically smooth, then 
the following hold for 
the discriminant divisor $\Delta_f$, 
the length of extremal ray $\mu_R$, and an extremal rational curve $\ell$ of $R$.  
  \begin{center}
    \begin{longtable}{cp{6cm}cp{5cm}}
      type of $R$ & $f$& $\mu_R$ & $\ell$ \\ \hline
      \multirow{2}{*}{$C_1$} & 
$\Delta_f$ is a nonzero effective divisor, 
if $f$ is generically smooth 
      & \multirow{2}{*}{$1$}& 
 \multirow{2}{*}{a curve in a singular fibre}\\ \hline
      \multirow{1}{*}{$C_2$}  & 
       \multirow{1}{*}{$\Delta_f=0$ and f is a $\P^1$-bundle} & \multirow{1}{*}{$2$} & a fibre  of $f$%
    \end{longtable}
  \end{center}

\end{thm}

\begin{proof}
We may assume that $f$ is generically smooth. 
Then  general fibres of $f$ are $\mathbb P^1$. 

Assume that $f$ is of type $C_1$, i.e., $f$ has a singular fibre. 
Then $\Delta_f \neq 0$. 
By $-K_X \cdot f^{-1}(y)=2$, we have that $-K_X \cdot \ell =1$ and $\mu_R =1$. 

Assume that $f$ is of type $C_2$, i.e., every fibre of $f$ is smooth. 
Then $\Delta_f =0$ and {\cred every}  fibre is isomorphic to $\mathbb P^1$. 
Since $Y$ is a smooth projective rational surface, we have $\Br (Y)=0$ (Proposition \ref{p-Br-vanish}(2)). 
It follows from Proposition \ref{p-Br-Pn} that $X$ is isomorphic to $\mathbb P(E)$ over $Y$ 
for some vector bundle $E$ of rank $2$. 
In this case, we have $\mu_R=2$ and $-K_X \cdot \ell =2$ for {\cred every}  fibre $\ell$. 
\end{proof}

\begin{cor}\label{fib_of_C}
  Let $X$ be a Fano threefold. 
  Let $R = \R_{\geq 0}[\ell]$ be an extremal ray of type $C$ with an extremal rational curve $\ell$. 
  Let $f\colon X\to Y$ be the contraction of $C$. 
  Then all the fibres of $f$ are numerically equivalent to $(2/\mu_R)\ell$, 
  where $\mu_R$ denotes the length of the extremal ray $R$. 
\end{cor}

\begin{proof}
Assume that $R$ is of type $C_1$. 
By Theorem \ref{dimy=2}, 
we have $\mu_R=1$ and $f^{-1}(y) \equiv 2 \ell$ for {\cred every}  closed point $y$. 
Note that this holds even if no fibre of $f$ is smooth. 

Assume that 
$R$ is of type $C_2$. 
Then we have $\mu_R =2$ and $f^{-1}(y) \equiv \ell$ for {\cred every}  closed point $y$. 
\end{proof}

\subsubsection{Type D}



\begin{thm}\label{dimy=1}
Let $X$ be a Fano threefold 
and let $R$ be an extremal ray of type $D$.  
Let $f \colon X\to Y$ be the contraction of $R$. Fix a closed point $y$ and $X_y := f^{-1}(y)$ denotes the scheme-theoretic fibre over $y$. 
Then $Y \simeq \P^1$ and $X_y$ is a projective {\cred Gorenstein} 
surface on $X$ such that $\omega_{X_y}^{-1}$ is ample. 
Furthermore, 
the following hold for 
the length of extremal ray $\mu_R$ and an extremal rational curve $\ell$ of $R$.  
  \begin{center}
    \begin{tabular}{clcl}
      type of $R$ & $f$& $\mu_R$ & $\ell$ \\ \hline
      $D_1$  & $1 \leq K_{X_y}^2 \leq 7$ & $1$ & \\ \hline
      $D_2$  & $K_{X_y}^2 =8$ & $2$ & a line on a fibre\\ \hline
      $D_3$  & $K_{X_y}^2 =9$ & $3$ & a line on a fibre\\
    \end{tabular}
  \end{center}
\end{thm}

\begin{proof}
By Proposition \ref{p-cont-ex2}, we obtain $Y \simeq \P^1$. 
The remaining assertions follow from Remark \ref{r-length} and 
Proposition \ref{p-dP-fibres}. 
\end{proof}

\subsubsection{Type E}


\begin{lem}\label{l-E1-prim}
Let $X$ be a Fano threefold and let $f : X \to Y$ be the contraction of an extremal ray $R$ of type $E_1$. 
Note that $C:=f(D)$ is a smooth curve on $Y$, $Y$ is a smooth projective threefold, and $f$ coincides with the blowup along $C$. 
Assume that $-K_Y$ is not ample. 
Then the following hold. 
  \begin{enumerate}
    \item $C\simeq \mathbb{P}^1$.
    \item $\mathcal{N}^*_{C/Y}\simeq \MO_{\P^1}(1)\oplus \MO_{\P^1}(1)$.
    \item $D\simeq \P^1 \times \P^1$.
    \item $\MO_X(K_X)|_D \simeq  \MO_{\P^1\times\P^1}(-1,-1)$ and $\MO_X(D)|_D \simeq \MO_{\P^1\times\P^1}(-1,-1)$.
  \end{enumerate}  
\end{lem}

\begin{proof}
It follows from Proposition~\ref{standard_3}(4)  that 
  \begin{equation}\label{e1-E1-prim} 
    -K_Y\cdot C +2-2g(C) =(-K_X)^2 \cdot D >0, 
  \end{equation}
where the inequality holds because $-K_X$ is ample. 
For $\mathcal{E} :=\mathcal{N}^*_{C/Y}\otimes \MO_C(-K_Y)$ and $D' :=\mathbb{P}(\mathcal{E})$, 
it follows from \cite[Ch. II, Lemma 7.9]{Har77} that 
\begin{itemize}
\item $\pi^{\prime}\colon D^\prime  \xrightarrow{\simeq, \varphi} D \xrightarrow{\pi} C$, and 
\item $\MO_{D'}(1) \simeq \varphi^*\MO_D(1) \otimes (\pi')^* \MO_C(-K_Y)$. 
\end{itemize}
It holds that  
\[
\mathcal{E} = \mathcal{N}^*_{C/Y}\otimes \MO_C(-K_Y) 
\simeq \mathcal{N}^*_{C/Y}\otimes \MO_C(-K_C)\otimes \bigwedge^2 \mathcal{N}_{C/Y} 
\simeq \MO_C(-K_C)\otimes \mathcal{N}_{C/Y}.
\]
 where the first isomorphism holds by \cite[Ch. II, Proposition 8.20]{Har77} and the second one follows from \cite[Ch. II, Exercise 5.16(b)]{Har77}. 
  Moreover, 
  \begin{align*}
    \MO_{D^\prime}(1) &\simeq \varphi^*\MO_D(1)\otimes (\pi')^*\MO_C(-K_Y) \\
    &\simeq \varphi^*\MO_D(1)\otimes (\varphi^*\pi^*\MO_C(-K_Y)) \\
    &\simeq \varphi^*(\MO_D(1)\otimes \pi^*\MO_C(-K_Y)) \\
    &\simeq \varphi^*( (\MO_X(-D) \otimes f^*\MO_Y(-K_Y))|_D)  \\
    &\simeq \varphi^*(\MO_X(-K_X)|_D), 
  \end{align*}
where the fourth and fifth isomorphisms follow from (2) and (1) of Proposition~\ref{standard_3}, respectively. 
  Since $-K_X$ is ample, also $\MO_{D^\prime}(1)$ is ample. 
  Hence $\mathcal{E}$ is an ample vector bundle on $C$ by \cite[Ch. III, Theorem 1.1]{Har70}.
  
\begin{clm} \label{c-E1-prim} 
It holds that  $-K_Y\cdot C \le 0$.  
\end{clm}

\begin{proof}[Proof of Claim] 

Suppose the contrary, i.e., $-K_Y\cdot C>0$. 
In order to derive a contraction, it suffices to prove that $-K_Y$ is ample. 
To this end, it is enough to show (i) and (ii) below \cite[Ch. I, Proposition 4.6]{Har70}.  
\begin{enumerate}
\item[(i)] $-K_Y \cdot Z>0$ for {\cred every}  curve $Z$ on $Y$. 
\item[(ii)] $-K_Y$ is semi-ample, i.e., $|-mK_Y|$ is base point free for some $m \in \Z_{>0}$. 
\end{enumerate}
  
Let us show (i). 
Fix a curve $Z$ on $Y$. 
If $Z=C$, then $-K_Y \cdot Z >0$ holds by our assumption. 
  Assume that $Z\neq C$. 
  Let $\widetilde{Z}$ be the proper transform of $Z$ on $X$.
Since $-K_X$ is ample and $D$ is an effective divisor with $\widetilde{Z} \not\subset \Supp\,D$, 
we have 
  \begin{align*}
    -K_Y\cdot Z =  -f^*K_Y \cdot \widetilde Z = (-K_X+D) \cdot \widetilde Z>0.
  \end{align*}
  Hence (i) holds. 

Let us show (ii). 
By (i), $-K_Y$ is nef.  
By $-K_Y\cdot C>0$, $\MO_Y(-K_Y)|_C$ is ample.
  Hence $\pi^*(\MO_Y(-K_Y)|_C)\simeq f^*(\MO_Y(-K_Y))|_D$ is semi-ample.
  By $f^*(-K_Y)=-K_X+D$, it follows from \cite[Theorem 3.2]{CMM14} that $f^*(-K_Y)$ is semi-ample.
  Since $X$ and $Y$ are normal, $-K_Y$ is semi-ample 
  {\cred \cite[Lemma 2.11(4)]{CT20}}.  
Thus (ii) holds. 
This completes the proof of Claim.
\end{proof}

By (\ref{e1-E1-prim}) and Claim, we have $g(C)=0$. 
Hence (1) holds. 
Since $\mathcal{E}=\MO_C(-K_C)\otimes \mathcal{N}_{C/Y}$ is a vector bundle of rank 2 on $C \simeq \mathbb P^1$, 
we can write
\[
\mathcal{E}\simeq \MO_{\mathbb{P}^1}(a)\oplus \MO_{\mathbb{P}^1}(b).
\]
for some $a,b\in \mathbb{Z}$. 
Since $\mathcal{E}$ is ample, we have $a>0$ and $b>0$ \cite[Ch. III, Corollary 1.8]{Har70}. 
By $\MO_{D^\prime}(1)\simeq \varphi^*(\MO_X(-K_X)|_D)$, we obtain 
\begin{align*}
a+b &=  \deg_C(\mathcal{E}) \\
&= c_1(\MO_{D^\prime}(1))^2 \\
  &= c_1(\varphi^*(\MO_X(-K_X)|_D))^2 \\
  &= (-K_X)^2 \cdot D \\
  &= -K_Y\cdot C +2 -g(C) \\
  &\le 2, 
\end{align*}
where the second equality holds by Proposition~\ref{standard_2}(4) 
and the fifth one follows from Proposition~\ref{standard_3}(4). 
Hence we get $a=b=1$ and 
$\mathcal{E}\simeq \MO_{\mathbb{P}^1}(1)\oplus \MO_{\mathbb{P}^1}(1)$. 

By $\mathcal{E}\simeq \MO_C(-K_C)\otimes \mathcal{N}_{C/Y}$, we have 
\begin{align*}
  \mathcal{N}_{C/Y}^* \simeq \mathcal{E}^* \otimes \omega_C^{-1} \simeq (\MO_{\mathbb{P}^1}(-1)\oplus \MO_{\mathbb{P}^1}(-1))\otimes \MO_{\mathbb{P}^1}(2) \simeq \MO_{\mathbb{P}^1}(1)\oplus \MO_{\mathbb{P}^1}(1).
\end{align*}
Thus (2) holds. 
By Proposition \ref{standard_3}(2), we obtain 
\[
D \simeq \mathbb P(\mathcal{N}_{C/Y}^*) \simeq \mathbb P(\MO_{\mathbb{P}^1}(1)\oplus \MO_{\mathbb{P}^1}(1)) \simeq \mathbb P^1 \times \mathbb P^1.
\]
Thus (3) holds.

Let us show (4). 
By $D\simeq \mathbb{P}^1\times \mathbb{P}^1$, we can write 
\[\MO_X(D)|_D\simeq \MO_{\mathbb{P}^1\times \mathbb{P}^1}(m,n)\]
for some $m,n\in \mathbb{Z}$.
It holds that 
\[
  2mn = c_1(\MO_{\mathbb{P}^1\times \mathbb{P}^1}(m,n))^2 
= c_1(\MO_X(D)|_D)^2 
\]
\[
= D^3 
= \deg_C(\mathcal{N}^*_{C/Y}) 
= \deg_{\mathbb{P}^1}(\MO_{\mathbb{P}^1}(1)\oplus \MO_{\mathbb{P}^1}(1)) 
= 2,
\]
where the fourth equality follows from Proposition~\ref{standard_3}(3). 
Hence we have $(m, n)=(1, 1)$ or $(m, n)=(-1, -1)$.  
Since $D \cdot \zeta <0$ for an $f$-exceptional curve $\zeta \subset D$, 
we obtain $(m, n)=(-1, -1)$. By the adjunction formula $\MO_X(K_X+D)|_D \simeq \MO_D(K_D) \simeq \MO_{\P^1 \times \P^1}(-2, -2)$, we obtain $\MO_X(K_X)|_D \simeq \MO_{\P^1 \times \P^1}(-1, -1)$. 
Thus (4) holds. 
\end{proof}

\begin{thm}\label{dimy=3}
Let $X$ be a primitive Fano threefold 
and let $R$ be an extremal ray of type $E$.  
Let $f \colon X\to Y$ be the contraction of $R$. 
Then $f$ is a birational morphism to a projective normal threefold $Y$. 
Furthermore, the following hold for 
the length $\mu_R$ of {\cred the} extremal ray {\cred $R$} and an extremal rational curve $\ell$ of $R$.  
  \begin{center}
      \begin{longtable}{cp{6cm}cp{4cm}}
      type of $R$ & $f$ and $D$ & $\mu_R$ & $\ell$ \\ \hline
      \multirow{5}{*}{$E_1$} & $Y$ is smooth, & \multirow{5}{*}{$1$}& $P\times\P^1$  \\
       & $C = f(D)\simeq \P^1$, & &\\
       & $\mathcal{N}^*_{C/Y}\simeq \MO_{\P^1}(1)\oplus \MO_{\P^1}(1),$ &&\\
       & $D\simeq \P^1\times\P^1,$ &&\\
       & $\MO_D(D)\simeq \MO_{\P^1\times \P^1}(-1,-1)$ &&\\ \hline
      \multirow{4}{*}{$E_2$}  & $Y$ is smooth, & \multirow{4}{*}{$2$} & a line on $D$ \\
       & $f(D)$ is a point, && \\
       & $D\simeq \P^2$, &&\\
       & $\MO_D(D)\simeq \MO_{\P^2}(-1)$ && \\ \hline
       \multirow{5}{*}{$E_3$} & $f(D)$ is a point, & \multirow{5}{*}{$1$} & $P\times\P^1$ or $\P^1\times Q$ \\
       & $D\simeq \P^1\times\P^1$, & & on $D$ $(P,Q\in \P^1)$ \\
       & $\MO_D(D)\simeq \MO_{\P^1\times\P^1}(-1,-1)$, &&\\
       & $P\times \P^1$ and $\P^1\times Q$ are numerically && \\
       & equivalent on $X$ for all $P, Q\in \P^1$ && \\ \hline
       \multirow{5}{*}{$E_4$} & $f(D)$ is a point, & \multirow{5}{*}{$1$} & a line on $D$ \\
       & $D$ is isomorphic to  && \\
       & the singular quadric surface && \\
       & in $\P^3$, && \\
       & $\MO_D(D)\simeq \MO_D \otimes \MO_{\P^3}(-1)$ && \\ \hline
       \multirow{3}{*}{$E_5$} & $f(D)$ is a point, & \multirow{3}{*}{$1$} & a line on $D$  \\ 
       & $D\simeq \P^2$, && \\ 
       & $\MO_D(D)\simeq \MO_{\P^2}(-2)$ && \\
      \end{longtable}
  \end{center}
    In particular, $\MO_X(-D)|_D$ is ample. 
\end{thm}

\begin{proof}
{\cred 
If $R$ is not of type $E_1$, then the assertion follows from 
\cite[Proposition 3.22]{TanII}. 
If $R$ is of type $E_1$, then the assertion holds by Lemma \ref{l-E1-prim}. 
Note that the columns $\mu_R$ and $\ell$ can be confirmed by Remark \ref{r-length}.} 
\end{proof}


\begin{cor}\label{run}
  Let $X$ be a primitive Fano threefold with an extremal ray $R$ of type $E$. 
Let $f: X \to Y$ be the contraction of $R$ and set $D:=\Ex(f)$. 
Then {\cred every}  effective divisor $Z$ on $D$ is semi-ample. 
\end{cor}

\begin{proof}
If $D \simeq \P^2$ or $D \simeq \P^1 \times \P^1$, then the assertion is clear. 
Hence we may assume that $D$ is a singular quadric surface in $\P^3$ (Theorem \ref{dimy=3}). 
In this case, the assertion follows from the fact that $D$ is $\Q$-factorial and $\rho(D) =1$ 
(it is well known that $D$ is a projective toric surface which is obtained by contracting the $(-2)$-curve on $\P_{\P^1}( \MO_{\P^1} \oplus \MO_{\P^1}(2))$). 
\end{proof}

\begin{cor}\label{E-to-pt}
  Let $X$ be a Fano threefold with an extremal ray $R$ of type $E_2, E_3, E_4$, or $E_5$. 
Let $f: X \to Y$ be the contraction of $R$ and set $D:=\Ex(f)$. 
Let $g\colon X\to \mathbb{P}^1$ be a morphism. 
Then $g(D)$ is a point.
\end{cor}

\begin{proof}
Assume that $g(D)$ is not a point, i.e., $g(D) = \mathbb P^1$. 
Then there exist curves $\Gamma_1$ and $\Gamma_2$ on $D$ 
such that $g(\Gamma_1)$ is a point and $g(\Gamma_2) =\mathbb P^1$. 
We have $g^*\MO_{\mathbb P^1}(1) \cdot \Gamma_1 =0$ and $g^*\MO_{\mathbb P^1}(1) \cdot {\cred \Gamma_2} >0$. 
Since the numerical equivalence classes $[\Gamma_1]$ and $[\Gamma_2]$ lie on $R\setminus \{0\}$, 
we can find $\lambda\in \mathbb{R}_{>0}$ satisfying $\Gamma_1 \equiv \lambda\Gamma_2$, 
which is  a contradiction.
\end{proof}

\subsection{Picard groups}

\begin{thm}\label{t-ex-surje}
  Let $X$ be a Fano threefold with $\rho(X) \geq 2$. 
  Let $R$ an extremal ray of $\NE(X)$ and let $f\colon X\to Y$ be the contraction of $R$. 
  Pick an extremal rational curve $\ell$ with $R = \R_{\geq 0}[\ell]$. 
  Then the sequence
  \[ 0 \longrightarrow \Pic\,Y \stackrel{f^*}{\longrightarrow} \Pic\,X \stackrel{(-\cdot \ell)}{\longrightarrow} \mathbb{Z} \longrightarrow 0 \]
is exact.   In particular, $\rho(X)= \rho(Y)+1$. 
\end{thm}

\begin{proof}
  By Proposition \ref{p-cont-ex}, it is enough to find a divisor $F$ on $X$ 
  such that $F \cdot \ell =1$. 
  If $\mu_R=1$, then we get  $-K_X\cdot\ell =\mu_R=1$. 
  Therefore, we are done for the case when the type of $R$ is one of 
  $C_1, D_1, E_1, E_3, E_4, E_5$ (Remark \ref{r-length}). 
  The remaining cases are $C_2, D_2, D_3$, and $E_2$.

  Assume that $R$ is of type $C_2$. 
  By Theorem~\ref{dimy=2}, $X$ is isomorphic to a $\mathbb{P}^1$-bundle over $Y$, associated to some locally free sheaf of rank $2$.
  Moreover, $\ell$ is a fibre of this bundle (Theorem~\ref{dimy=2}).
  Hence $\MO_X(1)\cdot \ell=1$ by Lemma~\ref{standard_1}(1).

  Assume that $R$ is of type $D_3$. 
  Then there exist Cartier divisors $D$ on $X$ and $E$ on $Y$ such that $K_X \sim 3D +f^*E$ (Proposition \ref{p-dP-fibres}(5)). 
  Then $3 = -K_X \cdot \ell = 3D \cdot \ell$, i.e., $D \cdot \ell =1$. 
  By using Proposition \ref{p-dP-fibres}(6), 
  the same argument works for the case when $R$ is of type $D_2$.

Assume that $R$ is of type $E_2$. 
Then the exceptional divisor $D$ on $X$ satisfies $D \simeq \mathbb{P}^2, \MO_D(D)\simeq \MO_D(-1)$, and $\ell$ is a line on $D$ \cite[Proposition 3.22]{TanII}. 
Hence $(-D\cdot \ell)=-(D|_D\cdot\ell)_D=-(c_1(\MO_{\mathbb{P}^2}(-1))\cdot c_1(\MO_{\mathbb{P}^2}(1)))_{\mathbb{P}^2}=1$.
\end{proof}



\begin{cor}\label{no_E1}
  Let $X$ be a Fano threefold with $\rho(X)=2$. Then 
  $X$ is primitive if and only if $X$ has no extremal ray of type $E_1$.
\end{cor}

\begin{proof}
If $X$ has no extremal ray of type $E_1$, 
then it is clear that $X$ is primitive. 
Suppose that $X$ is primitive and there exists an extremal ray $R$ of type $E_1$. 
It suffices to derive a contradiction. 
Let $f\colon X\to Y$ be the contraction of $R$. 
Then 
$Y$ is a smooth projective threefold and $f$ is the blowup  along a smooth curve $C$. 
It suffices to prove that $Y$ is Fano, because $X$ is a primitive Fano threefold. 


  By Corollary~\ref{picYrk1}, there is an ample divisor $L$ on $Y$ which generates $\Pic\,Y\simeq \mathbb{Z}$. Hence we can write $-K_Y\sim \alpha L$ for some $\alpha\in \mathbb{Z}$.
  By the Bertini theorem, we can find a curve $Z$ on $Y$ with $C \cap Z = \emptyset$. In particular, its proper transform $\tilde{Z}$ on $X$ is disjoint from $D := \Ex(f)$. Then 
  \begin{align*}
 (-K_X)\cdot \tilde{Z} = (f^*(-K_Y)-D)\cdot \tilde{Z} 
=-K_Y\cdot Z =\alpha L\cdot Z.
  \end{align*} 
  Since $-K_X$ and $L$ are ample, we get $\alpha>0$. Hence $-K_Y\sim \alpha L$ is ample. 
\end{proof}

\subsection{Existence of conic bundle structures}

\begin{lem}\label{has-CorD}
Let $X$ be a primitive Fano threefold with $\rho(X) \geq 2$. 
Then  $X$ has an extremal ray of type $C$ or $D$. 
\end{lem}

\begin{proof}
Suppose the contrary.
Let $R_1, \dots , R_n$ be all the extremal rays. 
We can find curve{\cred s} $\ell_1, ..., \ell_n$ such that 
\[
\NE(X) = \sum_{i=1}^n R_i = \sum_{i=1}^n \R_{\geq 0}[\ell_i]. 
\]
For each $i \in \{1, ..., n\}$, 
let $f_i : X \to Y_i$ be the contraction of $R_i$, 
which is of type $E$. 
Set $D_i := \Ex(f_i)$. 
Note that each $D_i$ is a projective normal toric surface (Theorem \ref{dimy=3}). 

  \begin{clm}
    $D_1 \cap D_i = \emptyset$ for all $2 \le i \le n$.
  \end{clm}
  \begin{proof}[Proof of Claim]
    Suppose that  $D_1\cap D_i \neq \emptyset$ for some $2 \leq i \leq n$.
    Since $D_i|_{D_1}$ is a nonzero effective Cartier divisor on $D_1$, 
there exists a curve $Z$ such that $Z \subset D_1 \cap D_i$. 
 
Since $\MO_X(-D_1)|_{D_1}$ is ample (Theorem \ref{dimy=3}), we have 
\[
D_1 \cdot Z = \MO_X(D_1) \cdot Z  = (\MO_X(D_1)|_{D_1}) \cdot Z <0. 
\] 
On the other hand, it follows from Corollary~\ref{run} that 
there exist $m \in \Z_{>0}$ and an effective Cartier divisor $Z'$ on $D_i$ such that 
$mZ \sim Z'$ and $\Supp\,Z'$ does not contain any irreducible component of $D_1 \cap D_i$. 
Then it holds that 
\[
D_1 \cdot Z = (D_1|_{D_i}) \cdot Z = \frac{1}{m}(D_1|_{D_i}) \cdot Z' \geq 0, 
\] 
    which is a contradiction. 
    This completes the proof of Claim. 
  \end{proof}

Fix $s \in \Z_{>0}$ such that $|-s K_X|$ is very ample. 
Pick two general members $H_1, H_2 \in |-s K_X|$ such that $H_1 \cap H_2$ is a smooth curve. 
By $\mathrm{NE}(X) = \R_{\geq 0} [\ell_1] + \cdots + \R_{\geq 0}[\ell_n]$, 
we have $H_1 \cdot H_2 = H_1 \cap H_2\equiv \sum_{i=1}^n a_i \ell_i$ for some $a_1, ..., a_n \in \R_{\geq 0}$. 
We then obtain the following contradiction: 
  \[
0< D_1 \cdot (-sK_X)^2=D_1\cdot H_1\cdot H_2  = a_1 D_1\cdot\ell_1  + \sum_{i=2}^n a_i D_1\cdot \ell_i =a_1 D_1 \cdot \ell_1 \leq 0, 
\]
where Claim implies $D_1\cdot \ell_2= \cdots =D_1 \cdot \ell_n =0$ and the inequality 
$D_1 \cdot \ell_1<0$ follows from the ampleness of $\MO_X(-D_1)|_{D_1}$ (Theorem \ref{dimy=3}). 
\end{proof}

\begin{lem}\label{l-finite-morph}
Let $X$ be a primitive Fano threefold. 
Let $R_1$ and $R_2$ be two distinct extremal rays. 
For each $i \in \{1, 2\}$, 
let $f_i: X \to Y_i$ be the contraction of $R_i$. 
Then the following hold. 
\begin{enumerate}
\item 
The induced morphism $f_1 \times f_2 : X \to Y_1 \times Y_2$ is a finite morphism. 
\item If $f_1$ is of type $C$ and $f_2$ is of type $E$, then the composite morphism 
\[
f_1|_D : D \overset{j}{\hookrightarrow} X \xrightarrow{f_1} Y_1
\]
is a finite surjective morphism, where $D :=\Ex(f_2)$ is the $f_2$-exceptional prime divisor and $j:D \hookrightarrow X$ denotes the induced closed immersion. 
\end{enumerate}

\end{lem}

\begin{proof}
Let us show (1). 
Suppose that $f_1 \times f_2 : X \to Y_1 \times Y_2$ is not a finite morphism. 
Then there exists a curve $\Gamma$ on $X$ such that $(f_1 \times f_2)(\Gamma)$ is a point. 
Then both $f_1(\Gamma)$ and $f_2(\Gamma)$ are points. 
However, this implies $[\Gamma] \in R_1 \cap R_2 =\{0\}$, which is a contradiction. 
Thus (1) holds. 

Let us show (2). 
By $\dim D  = \dim Y_1$, it suffices to show that $f_1|_D$ is a finite morphism. 
Suppose that $f_1|_D : D \to Y_1$ is not a finite morphism. 
Then there exists a curve $\Gamma$ on $X$ such that $\Gamma \subset D$ and $f_1(\Gamma)$ is a point. 
If $R_2$ is of type $E_2, E_3, E_4$, or $E_5$, then also $f_2(\Gamma)$ is a point. 
This is a contradiction: $[\Gamma] \in R_1 \cap R_2 =\{0\}$. 
Hence $R_2$ is of type $E_1$. 
Since $\MO_X(-D)|_D$ is ample, we have that $D \cdot \Gamma <0$. 
Let $\Gamma'$ be a general fibre of $f_1 : X \to Y_1$, so that $\Gamma' \not\subset D$. 
Hence we get $D \cdot \Gamma'\geq 0$. 
On the other hand, we have $\Gamma' \equiv \lambda \Gamma$ for some $\lambda \in \R_{>0}$, 
which implies $D \cdot \Gamma \geq 0$. 
This is a contradiction. 
\end{proof}

\begin{prop}\label{hasC}
Let $X$ be a primitive Fano threefold with $\rho(X) \geq 2$. 
Then  $X$ has an extremal ray of type $C$.
\end{prop}

\begin{proof}
By Lemma \ref{has-CorD}, we may assume that there exists an extremal ray $R_1$ of type $D$. 
Then we have $\rho(X)=2$ (Theorem \ref{t-ex-surje}). 
Let $R_2$ be the other extremal ray. 
For each $i \in \{1, 2\}$, let $f_i : X \to Y_i$ be the contraction of $R_i$, where $Y_1 =\mathbb P^1$. 
Suppose that $R_2$ is of type $E$ or $D$. 
It suffices to derive a contradiction. 
If $R_2$ is of type $D$, then Lemma \ref{l-finite-morph}(1) leads to a contradiction. 

The problem is reduced to the case when $R_2$ is of type $E$. 
By Corollary~\ref{no_E1}, $R_2$ is not of type $E_1$, and hence $R_2$ is of type $E_2, E_3, E_4,$ or $E_5$. 
Set $D_2 := \Ex(f_2)$. 
It follows from Corollary~\ref{E-to-pt} that $f_1(D_2)$ is a point. 
Pick a curve $Z$ with $Z \subset D_2$. 
Then both $f_1(Z)$ and $f_2(Z)$ are points. 
Hence we get $[Z] \in R_1 \cap R_2 = \{0\}$, which contradicts $-K_X \cdot Z>0$. 
\end{proof}

\subsection{$\rho(X) \leq 3$}


\begin{thm}\label{t-cb-base}
Let $X$ be a primitive Fano threefold. 
Let $R$ be an extremal ray of type $C$ and 
let $f: X \to S$ be the contraction of $R$. 
Then $S \simeq \mathbb P^2$ or $S \simeq \mathbb P^1 \times \mathbb P^1$. 
\end{thm}

\begin{proof}
If $f$ is a wild conic bundle, 
then the assertion follows from Theorem \ref{t-wild-cb}. 
Hence we may assume that $f$ is generically smooth. 
Note that $S$ is a smooth projective rational surface (Theorem \ref{dimy=2}). 
By the classification of smooth projective rational surfaces, 
it suffices to show that  there exists no curve $E$ on $S$ such that $E^2 < 0$. 
Suppose that there is a curve $E$ on $S$ such that $E^2 < 0$. 
Let us derive a contradiction. 

Fix a curve $C$ on $X$ such that $C \subset f^{-1}(E)$ and $f(C)=E$. 
We have 
\[
\NE(X) = \sum_{i=1}^n R_i =\sum_{i=1}^n \R_{\geq 0}[\ell_i], \qquad R_i = \R_{\geq 0}[\ell_i]
\]
for the extremal rays $R_1=R, R_2, ..., R_n$ and curves $\ell_1, ..., \ell_n$ on $X$. 
We then have 
  \[C\equiv\sum_{i=1}^n a_i \ell_i\]
  for some $a_1, ..., a_n \in \R_{\geq 0}$. 
  
  Now, we can write $f_*C = bE$ for some $b \in \mathbb{Q}_{>0}$.
  By $E^2 <0$, we get 
  \[
  \sum_{i=1}^na_i E\cdot f_*\ell_i =
  E \cdot \sum_{i=1}^na_i f_*\ell_i = E\cdot f_*C =b E^2 <0.\]
  Hence, possibly after permuting the indices,  we have $E\cdot f_*\ell_2<0$ (note that $f_*\ell_1 =0$). 
  Let $f_2 : X \to Y_2$ be the contraction of $R_2$. 
By $E\cdot f_*\ell_2<0$ and $E^2<0$, we get   $f(\ell_2)=E$. 
By the projection formula, we obtain 
  \begin{equation}\label{ne_cu}
   f^*E\cdot \ell_2=E \cdot (f_*\ell_2)  <0.
  \end{equation}

  Assume that $R_2$ is of type $C$ or $D$. 
  Then the fibres of $f_2 : X \to Y_2$ are curves or surfaces. 
 By $f^{-1}(E) \neq X$, there is a curve $Z$ on $X$ such that 
\begin{enumerate}
\item[(i)] $Z \not\subset f^{-1}(E)$ and 
\item[(ii)] $f_2(Z)$ is a point. 
\end{enumerate}
By (ii), we can write $Z\equiv c \ell_2$ for some $c\in \mathbb{R}_{>0}$. 
Therefore, we obtain 
\[
f^*E \cdot \ell_2 =c^{-1} f^*E \cdot Z \geq 0, 
\]
where the inequality holds by (i). 
This contradicts (\ref{ne_cu}).

Assume that $R_2$ is of type $E$. 
Set $D_2 := \Ex(f_2)$. 
Since $D_2$ is covered by curves $Z$ with $[Z] \in R_2 = \R_{\geq 0}[\ell_2]$, 
it follows from (\ref{ne_cu}) that $D_2 \subset f^{-1}(E)$. 
By $f(D_2) \subset f(f^{-1}(E)) =E$, 
there exists a curve $C'$ on $D_2$ such that $f(C')$ is a point. 
We get $[C'] \in R$. 
Since $\MO_X(-D_2)|_{D_2}$ is ample (Theorem \ref{dimy=3}), 
we have 
\[
D_2 \cdot C'<0.
\]
On the other hand, 
{\cred it follows from $f(D_2) \subset E$ that} 
$D_2$ is disjoint from $f^{-1}(s)$ for a general closed point $s \in S$, and hence we get 
\[
D_2 \cdot f^{-1}(s) =0. 
\]
This is a contradiction, because we have $[C'] \in R$ and $[f^{-1}(s)] \in R$. 
\end{proof}

\begin{thm}\label{pic-over2}
  Let $X$ be a primitive Fano threefold. 
Then the following hold. 
{\cred 
\begin{enumerate}
\item $\rho(X) \leq 3$.   
\item If $\rho(X) =2$, then the following hold. 
\begin{enumerate}
\item[(2a)] There exists an extremal ray of type $C$. 
\item[(2b)] If $f: X \to S$ is a contraction of an extremal ray of type $C$, then 
$S \simeq \P^2$. 
\end{enumerate}
\item If $\rho(X) =3$, then the following hold. 
\begin{enumerate}
\item[(3a)] There exists an extremal ray of type $C$. 
\item[(3b)] If $f: X \to S$ is a contraction of an extremal ray of type $C$, then 
$S \simeq \P^1 \times \P^1$ and $f$ is generically smooth. 
\item[(3c)] Every extremal ray of $X$ is of type $C$ or $E_1$. 
\end{enumerate}
\end{enumerate}
}

\end{thm}



\begin{proof}
{\cred 
By Proposition \ref{hasC}, 
(2a) and (3a) hold. 
Then (1), (2b), and (3b) follow from $\rho (X)  = \rho(S)+1$, 
Theorem \ref{t-wild-cb}, and Theorem \ref{t-cb-base}.}


Let us show {\cred (3c)}. 
Suppose that {\cred $\rho(X)=3$ and} there exists an extremal ray $R$ which is not of type $C$ nor $E_1$. 
Let us derive a contradiction. 
{\cred 
Note that $R$ is not of type $D$, 
as otherwise we would get $\rho(X) = 2$ (Theorem \ref{t-ex-surje}), 
which contradicts $\rho(X)=3$. 
Hence} 
the type of $R$ is $E_2, E_3, E_4$, or $E_5$. 
Let $f: X \to Y$ be the contraction of $R$. 
Set $D := \Ex(f)$. 
There exists an extremal ray $R'$ of type $C$ (Proposition \ref{hasC}). 
Let $f': X \to \mathbb P^1 \times \mathbb P^1$ be the contraction of $R'$. 
For each $i \in \{1, 2\}$, consider 
\[
g_i : X \xrightarrow{f'} \mathbb P^1 \times \mathbb P^1 \xrightarrow{{\rm pr}_i} \mathbb P^1. 
\]
By Corollary \ref{E-to-pt}, we have that $g_i(D)$ is a point for each $i \in \{1, 2\}$. 
Hence also $f'(D)$ is a point. 
However, this is a contradiction, because {\cred every}  fibre of $f'$ is one-dimensional. 
\qedhere

\end{proof}

\section{Fano threefolds with $\rho(X)=2$}\label{s-pic2}

In this section, we classify Fano threefolds with $\rho(X)=2$. 
If $X$ has an extremal ray of type C and its contraction is not generically smooth 
(i.e., a wild conic {\cred bundle}), 
then $X$ is a prime divisor on $\P^2 \times \P^2$ of bidegree $(1, 2)$ (Theorem \ref{t-wild-cb}). 
In what follows, we always assume that the contraction of type C is  generically smooth. 

{\cred 
Recall that a Fano threefold $X$ with $\rho(X)=2$ 
has an extremal ray of type $C$ or $E_1$ (Theorem \ref{pic-over2}). 
We shall treat the former (resp. latter) case 
in Subsections \ref{ss-rho2-CC}, \ref{ss-rho2-CD}, \ref{ss-rho2-CE} (resp. Subsections \ref{ss-rho2-E1C}, \ref{ss-rho2-E1D}, \ref{ss-rho2-E1E}).} 

\begin{nota}\label{n-pic-2}
Let $X$ be a Fano threefold with $\rho(X)=2$. 
By Theorem \ref{cone}, $\mathrm{NE}(X)$ has exactly two extremal rays $R_1$ and $R_2$. 
For each $i \in \{1, 2\}$, 
let 
\[
f_i : X \to Y_i
\] 
be the contraction of $R_i$ and 
let $\ell_i$ be an extremal rational curve with $R_i = \R_{\geq 0}[\ell_i]$. 
Set $\mu_i := -K_X \cdot \ell_i$, which is the length of $R_i$.  
We have 
\[
\NE(X) = R_1 + R_2  =\R_{\geq 0}[\ell_1] + \R_{\geq 0}[\ell_2].
\]
It follows from  Corollary~\ref{picYrk1} that $\Pic (Y_i) \simeq \Z$ for each $i \in \{1, 2\}$. 
Let $L_i$ be an ample Cartier divisor on $Y_i$ such that $\MO_{Y_i}(L_i)$ generates $\Pic\,Y_i$. 
Set $H_i := f_i^*L_i$. 
\begin{itemize}
\item If $R_i$ is of type $C$, 
then we assume that $f_i$ is generically smooth and 
$\Delta_{f_i}$ denotes its discriminant divisor (Subsection \ref{ss-CB}). 
\item If $R_i$ is of type $D$, 
then $d_i := (-K_X)^2 \cdot H_i$. 
Note that $1 \leq d_i \leq 9$ (Theorem \ref{dimy=1}).  
\item 
Assume that $R_i$ is of type $E$. 
Set $D_i := \Ex(f_i)$. 
If $R_i$ is of type $E_1$, then 
$B_i := f_i(D_i)$. 
If $R_i$ is of type $E_1, E_3,$ or $E_4$ 
(resp. $E_2$, resp. $E_5$), 
then let $r_i$ be the largest positive integer which divides $-K_X+D_i$ (resp. $-K_X+2D_i$, resp. $-2K_X+D_i$). 
Note that $r_i$ coincides with the index of the Fano threefold $Y_i$ when $R_i$ is of type $E_1$ or $E_2$. 
\end{itemize}
\end{nota}

\subsection{The Picard group and the canonical divisor}\label{ss-rho2-generalities}


\begin{lem}\label{l-Hi-prim}
We use Notation \ref{n-pic-2}. 
Then, for each $i \in \{1, 2\}$, 
$H_i$ is a primitive element in $\Pic\,X$, i.e., 
there exists no pair $(s, H)$ such that 
$s \in \Z_{\geq 2}$, $H$ is a Cartier divisor on $X$, and $H_i \sim sH$. 
\end{lem}

\begin{proof}
  Suppose that there exist a Cartier divisor $H$ on $X$ and 
  an integer $s \geq 2$ such that $H_i \sim sH$. 
  By $H_i \cdot \ell_i= f_i^* L_i \cdot \ell_i = 0$, 
  we obtain $H\cdot \ell_i =0$. 
By Theorem \ref{t-ex-surje}, there exists a Cartier divisor $M$ on $Y_i$ 
such that $H \sim f_i^*M$. 
Then we get $0 \sim H_i - sH \sim f_i^*L_i -s f_i^*M$, which implies $L_i \sim s M$. 
This contradicts the fact that $\MO_{Y_i}(L_i)$ is a primitive element in $\Pic\,Y_i$. 
\end{proof}

\begin{lem}\label{pic2-numequiv}
We use Notation \ref{n-pic-2}. 
Fix $i \in \{1, 2\}$ and assume that $R_i$ is of type $C$. 
Then 
  \[
H_i^2\equiv \frac{2}{\mu_i}\ell_i. 
\]
\end{lem}

\begin{proof}
By $Y_i\simeq\mathbb{P}^2$ (Theorem~\ref{pic-over2}(1)), we have $\MO_X(L_i) \simeq \MO_{\mathbb{P}^2}(1)$. 
Hence $L_i^2$ is a point on $\mathbb{P}^2$. 
Then $H_i^2=f_i^*L_i^2$ is a fibre of $f_i$, which is numerically equivalent to $\frac{2}{\mu_i}\ell_i$ (Corollary~\ref{fib_of_C}).
\end{proof}

We now compute $c_2(X)\cdot H_i$ for each $i \in \{1, 2\}$. 

\begin{lem}\label{c_times_H}
We use Notation \ref{n-pic-2}. Fix $i \in \{1, 2\}$. 
Then the following hold. 
\begin{enumerate}
\item 
Assume that $R_i$ is of type $C$. 
Recall that $\Delta_{f_i}$ denotes the discriminant divisor of the 
generically smooth conic bundle $f_i : X \to Y_i =\P^2$. 
Then  $c_2(X) \cdot H_i$ satisfies the following. 
  \begin{center}
    \begin{tabular}{c||c|c}
      type of $R_i$ & $C_1$ & $C_2$ \\ \hline
      $c_2(X)\cdot H_i$ & $6+\deg\Delta_{f_i}$ & $6$ \\
    \end{tabular}
  \end{center}
Moreover, $\deg\Delta_{f_i} \geq 1$ if $R_i$ is of type $C_1$.
\item 
Assume that $R_i$ is of type $D$. Then  $c_2(X) \cdot H_i$ satisfies the following. 
  \begin{center}
    \begin{tabular}{c||c|c|c}
      type of $R_i$ & $D_1$ & $D_2$ & $D_3$ \\ \hline
      $c_2(X)\cdot H_i$ & $12-(-K_X)^2 \cdot H_i$ & $4$& $3$\\
    \end{tabular}
  \end{center}
\item 
Assume that $R_i$ is of type $E$. 
Set $D_i := \Ex(f_i)$. 
Then $c_2(X) \cdot H_i$ satisfies the following: 
  \begin{center}
    \begin{tabular}{c||c|c|c|c}
      type of $R_i$ & $E_1$ & $E_2$ & $E_3$ or $E_4$ & $E_5$ \\ \hline
      $c_2(X)\cdot H_i$ & $\frac{24}{r_i} + \deg B_i$ & $\frac{24}{r_i}$ & $\frac{24}{r_i}$& $\frac{45}{r_i}$\\
    \end{tabular}
  \end{center}
where we set 
$B_i := f_i(D_i)$ and $\deg B_i := \frac{-K_{Y_i} \cdot B}{r_i}$ for the case when $f_i$ is of type $E_1$. 
\end{enumerate}
\end{lem}

\begin{proof}
Let us show (1) {\cred and (2)}. 
Assume that $R_i$ is of type $C$ {\cred (resp. $D$)}. 
Then we have $f_i\colon X\to Y_i =\mathbb{P}^2$ 
{\cred (resp. 
$f_i\colon X\to Y_i =\mathbb{P}^1$)}. 
Since $L_i$ is an ample generator of $\Pic\,\mathbb{P}^2$ 
{\cred (resp. $\Pic\,\mathbb{P}^1$)}, 
we may assume that $L_i$ is a line on $Y_i =\mathbb P^2$ {\cred (resp. a point on $\P^1$)}.  
  Since $H_i =f_i^*L_i$ is an effective divisor, 
  it follows from Lemma~\ref{firstlem}(4) that 
  \[
c_2(X)\cdot H_i = 6\chi(H_i, \MO_{H_i}) + 6\chi(H_i, \MO_X(H_i)|_{H_i}) -2H_i^3 - (-K_X)^2\cdot H_i. \]
We have $H_i^3 = f_i^*L_i^3 =0$.

Assume that $R_i$ is of type $C$. 
It holds that 
\begin{itemize}
\item $\chi(H_i, \MO_{H_i})=1-\chi(\P^2, \MO_{\mathbb{P}^2}(-L_i))=1-\chi(\P^2, \MO_{\mathbb{P}^2}(-1))=
1$ (Corollary~\ref{euler_cont}), 
\item $\chi(H_i, \MO_X(H_i)|_{H_i})=\chi(\P^2, \MO_{\mathbb{P}^2}(L_i))-1=\chi(\P^2, \MO_{\mathbb{P}^2}(1))-1=2$ (Corollary~\ref{euler_cont}), and 
\item 
$(-K_X)^2\cdot H_i=-4K_{\mathbb{P}^2}\cdot L_i-\Delta_{f_i}\cdot L_i =12-\deg(\Delta_{f_i})$ (Proposition~\ref{conic_disc-num}). 
\end{itemize}
To summarise, we get  
\[
c_2(X)\cdot H_i = 6+12-0-(12-\deg(\Delta_{f_i})) 
= 6+\deg(\Delta_{f_i}).
\]
Thus (1) holds. 

\medskip

Assume that $R_i$ is of type $D$. 
It holds that 
\begin{itemize}
\item $\chi(H_i, \MO_{H_i})=1-\chi(\P^1, \MO_{\mathbb{P}^1}(-L_i))=1-\chi(\MO_{\mathbb{P}^1}(-1))=1$ (Corollary~\ref{euler_cont}), and 
\item $\chi(H_i, \MO_X(H_i)|_{H_i})=\chi(\P^1, \MO_{\mathbb{P}^1}(L_i))-1=\chi(\MO_{\mathbb{P}^1}(1))-1=1$ (Corollary~\ref{euler_cont}). 
\end{itemize}
To summarise, we obtain 
  \begin{align*}
    c_2(X)\cdot H_i =6+6-0- (-K_X)^2 \cdot H_i =12-(-K_X)^2 \cdot H_i.
  \end{align*}
Thus (2) holds.

  \medskip

Let us show (3). 
{\cred 
Assume that $R_i$ is of type $E$. Set
\[
{\small 
a:= 
\begin{cases}
1 \quad (\text{if $R_i$ is of type $E_1, E_3$, or $E_4$})\\ 
2 \quad (\text{if $R_i$ is of type $E_2$})\\
\frac{1}{2} \quad (\text{if $R_i$ is of type $E_5$})\\
\end{cases}
\epsilon := 
\begin{cases}
1 \quad (\text{if $R_i$ is of type $E_1, E_2, E_3$, or $E_4$})\\ 
2 \quad (\text{if $R_i$ is of type $E_5$}). 
\end{cases}
}
\]
Then it is easy to see that 
$K_X \equiv f_i^*K_{Y_i} + a D_i$ and 
\[
\epsilon(-K_X + aD_i) \sim -\epsilon f_i^*K_{Y_i} \sim r_iH_i. 
\]
  By Lemma~\ref{firstlem}(3), we obtain 
  \[
 c_2(X)\cdot H_i = \frac{\epsilon}{r_i}(c_2(X)\cdot(-K_X) +a c_2(X)\cdot D_i) 
 =\frac{\epsilon }{r_i}(24+a c_2(X)\cdot D_i).
\]
It follows from Lemma~\ref{firstlem}(4) that 
  \[
  c_2(X)\cdot D_i=6\chi(D_i, \MO_{D_i}) + 6\chi(D_i, \MO_X(D_i)|_{D_i})-2D_i^3-(-K_X)^2\cdot D_i. 
  \]
In what follows, we apply case study to compute $c_2(X) \cdot D_i$.} 

\medskip

Assume that $R_i$ is of type $E_1$. 
Let  $g(B_i)$ be the genus of $B_i$. 
We have $d_i := \deg B_i = (-\frac{1}{r_i}K_{{\cred Y_i}}) \cdot B_i$. 
By \cite[Lemma 3.21(2)]{TanII}, we get 
\begin{itemize}
\item $(-K_X)^2 \cdot D_i =r_id_i -2g(B_i)+2$, 
\item $(-K_X) \cdot D_i^2 = 2g(B_i)-2$, and 
\item $D_i^3 =  -r_id_i +2-2g(B_i)$.
\end{itemize}
It holds that $\chi(D_i, \MO_{D_i}) = 1 - g(B_i)$ and 
\[
\chi(D_i, \MO_X(D_i)|_{D_i}) = \chi(D_i, \MO_{D_i}) +\frac{1}{2} (D_i|_{D_i}) \cdot (D_i|_{D_i}-K_{D_i}) 
\]
\[
= 1-g(B_i) + \frac{1}{2} (-K_X) \cdot D_i^2 = 0. 
\]
We then get 
\begin{align*}
c_2(X)\cdot D_i 
&=6\chi(D_i, \MO_{D_i}) + 6\chi(D_i, \MO_X(D_i)|_{D_i})-2D_i^3-(-K_X)^2\cdot D_i\\
&=6(1-g(B_i)) + 0 -2(-r_id_i +2-2g(B_i)) - (r_id_i -2g(B_i)+2) =r_id_i. 
  \end{align*}
To summarise,  we obtain 
\[
{\cred 
c_2(X)\cdot H_i = 
\frac{\epsilon }{r_i}(24+a c_2(X)\cdot D_i)
= \frac{1}{r_i}(24 + 1 \cdot r_id_i) = \frac{24}{r_i} + \deg B_i,} 
 \]
which completes the proof for the case when $R_i$ is of type $E_1$.

\medskip

Assume that $R_i$ is of type $E_2$. 
By \cite[Proposition 3.22]{TanII}, 
we get  
$\mu_i = 2, D_i\simeq \mathbb{P}^2, D_i^3 =1, 
(-K_X)^2 \cdot D_i =4$, and  
$\MO_X(D_i)|_{D_i}\simeq \MO_{D_i}(-1)$. 
It holds that $\chi(D_i, \MO_{D_i}) =1$ 
and $\chi(D_i, \MO_{D_i}(D_i))=0$. 
We then get 
\begin{align*}
c_2(X)\cdot D_i 
&=6\chi(D_i, \MO_{D_i}) + 6\chi(D_i, \MO_X(D_i)|_{D_i})-2D_i^3-(-K_X)^2\cdot D_i\\
&=6 + 0 -2 - 4 =0. 
  \end{align*}
Thus we obtain 
\[
{\cred 
c_2(X)\cdot H_i = 
\frac{\epsilon }{r_i}(24+a c_2(X)\cdot D_i)
= \frac{1}{r_i}(24 + 0) = \frac{24}{r_i},} 
 \]
which completes the proof for the case when 
$R_i$ is of type $E_2$.

  
\medskip 

Assume that $R_i$ is of type $E_3$ or $E_4$. 
It follows from  \cite[Proposition 3.22]{TanII} 
that $\mu_i =1$, $(-K_X)^2 \cdot D_i=2$, $D_i^3=2$, 
$D_i$ is isomorphic to a possibly singular quadric surface in $\mathbb{P}^3$, and 
 $\MO_X(D_i)|_{D_i} \simeq \MO_{\P^3}(-1)|_{D_i}$.  
It holds that $\chi(D_i, \MO_{D_i})=h^0(D_i, \MO_{D_i})=1$ 
and 
$\chi(D_i, \MO_{D_i}(D_i))=\chi(D_i, \MO_{\mathbb P^3}(-1)|_{D_i})=0$. 
Therefore, we obtain 
\begin{align*}
c_2(X)\cdot D_i 
&=6\chi(D_i, \MO_{D_i}) + 6\chi(D_i, \MO_X(D_i)|_{D_i})-2D_i^3-(-K_X)^2\cdot D_i\\
&=6 + 0-4-2=0. 
  \end{align*}
We then get 
\[
{\cred 
c_2(X)\cdot H_i = 
\frac{\epsilon }{r_i}(24+a c_2(X)\cdot D_i)
= \frac{1}{r_i}(24 + 0) = \frac{24}{r_i},} 
\]
which 
completes the proof for the case when 
$R_i$ is of type $E_3$ or $E_4$.

\medskip

Assume that $R_i$ is of type $E_5$. 
It follows from  \cite[Proposition 3.22]{TanII} 
that 
$\mu_i=1$, $D_i\simeq \mathbb{P}^2$, 
$D_i^3 =4$, $(-K_X)^2 \cdot D_i = 1$, and 
$\MO_X(D_i)|_{{\cred D_i}} \simeq \MO_{\mathbb{P}^2}(-2)$. 
It holds that 
$\chi(D_i, \MO_{D_i})=1$ and 
$\chi(D_i, \MO_X(D_i)|_{D_i})=0$. 
We obtain
\begin{align*}
c_2(X)\cdot D_i 
&=6\chi(D_i, \MO_{D_i}) + 6\chi(D_i, \MO_X(D_i)|_{D_i})-2D_i^3-(-K_X)^2\cdot D_i\\
&= 6+0-8-1= -3. 
  \end{align*}
We then get 
\[
{\cred 
c_2(X)\cdot H_i = 
\frac{\epsilon }{r_i}(24+a c_2(X)\cdot D_i)
= \frac{2}{r_i}(24 +\frac{1}{2} \cdot (-3)) = \frac{45}{r_i},} 
\]
which 
completes the proof for the case when 
$R_i$ is of type $E_5$. 
Thus (3) holds. 
\end{proof}

\begin{lem}\label{c_times_H2}
We use Notation \ref{n-pic-2}. Fix $i \in \{1, 2\}$. 
\begin{enumerate}
\item If $R_i$ is of type $C_1$, then $7 \leq c_2(X)\cdot H_i \leq 17$. 
\item If $R_i$ is not of type $E_1$, then $1 \leq c_2(X) \cdot H_i \leq 24$ or $c_2(X) \cdot H_i =45$. 
\item If $R_i$ is neither of type $E_1$ nor $E_5$, 
then $1 \leq c_2(X) \cdot H_i \leq 24$. 
\end{enumerate}
\end{lem}

\begin{proof}
Let us show (1).  
Assume that $R_i$ is of type $C_1$. 
By the same argument as in the proof of Lemma~\ref{c_times_H}(1), 
\[
c_2(X)\cdot H_i = 18-(-K_X)^2\cdot H_i.
\]
  Since $-K_X$ is ample, we get $c_2(X)\cdot H_i \le 17$. 
  On the other hand, it follows from 
  Lemma~\ref{c_times_H}(1) that 
  \[ 
  c_2(X)\cdot H_i=6+\deg \Delta_{f_i}\]
  and $\deg \Delta_{f_i}>0$. 
  In particular, $c_2(X)\cdot H_i\ge 7$.
 Thus (1) holds. 

Let us show (2). 
If $R_i$ {\cred is} not of type $D_1$, then these assertions follow from (1) and 
Lemma~\ref{c_times_H}. 
Assume that $R_i$ is of type $D_1$. 
Then  Lemma~\ref{c_times_H}{\cred (2)} implies 
  \[c_2(X)\cdot H_i=12-(-K_X)^2 \cdot H_i.\]
  Since $1\le (-K_X)^2 \cdot H_i \le 7$, we have $5 \le c_2(X)\cdot H_i \le 11$. Thus (2) holds.  
  Similarly, (3) holds. 
\end{proof}


\begin{lem}\label{l-basis}
We use Notation \ref{n-pic-2}. 
Set 
\[
a :=  |\Pic(X) / (\mathbb{Z}H_1 
+
\mathbb{Z}H_2)|, 
\]
which is the cardinality of the set $\Pic(X) / (\mathbb{Z}H_1 +  \mathbb{Z}H_2)$. 
Then the following hold. 
\begin{enumerate}
\item $H_1$ and $H_2$ are linearly 
{\cred independent} over $\R$ in $(\Pic\,X)\otimes_{\Z} \R/\equiv$. 
Moreover, $a \in \Z_{>0}$ 
{\cred 
and $\mathbb{Z}H_1 \oplus \mathbb{Z}H_2 = \mathbb{Z}H_1 + \mathbb{Z}H_2$}.    
\item $\Pic(X) / (\mathbb{Z}H_1 \oplus \mathbb{Z}H_2)\simeq \mathbb{Z} / (H_2\cdot \ell_1)\mathbb{Z} \simeq \Z/(H_1 \cdot \ell_2)\Z$. 
\item $a = H_2 \cdot \ell_1 = H_1 \cdot \ell_2$. 
\item  $-aK_X \sim \mu_2 H_1 + \mu_1 H_2.$
\item $24a = \mu_2 H_1 \cdot c_2(X) + \mu_1 H_2 \cdot c_2(X)$. 
\end{enumerate}
\end{lem}

\begin{proof}
Let us show (1). 
Assume 
  \[
  \lambda_1H_1+\lambda_2H_2 \equiv 0
  \]
  for $\lambda_1, \lambda_2\in\mathbb{R}$.
  By $H_1=f^*_1L_1$, we have $H_1\cdot \ell_1= 0$. 
  Hence we obtain 
$\lambda_2 H_2\cdot \ell_1 =0$. 
Since $f_2(\ell_1)$ is a curve, we have $H_2 \cdot \ell_1 >0$. 
We then get $\lambda_2=0$, which implies $\lambda_1=0$. 
Therefore, $H_1$ and $H_2$ are linearly independent over $\mathbb{R}$. 
{\cred In particular, we get $\mathbb{Z}H_1 \oplus \mathbb{Z}H_2 = \mathbb{Z}H_1 + \mathbb{Z}H_2$.} 
Since $\mathbb{Z}H_1\oplus \mathbb{Z}H_2$ is a free subgroup of rank $2$ of $\Pic\,X \simeq \Z^{\oplus 2}$, the quotient group 
\[
\Pic(X) / (\mathbb{Z}H_1\oplus \mathbb{Z}H_2)
\]
is a finite abelian group.
Thus (1) holds. 


Let us show (2) and (3). 
By $H_2 \cdot \ell_1 > 0$ and $H_1 \cdot \ell_2 > 0$, 
(2) implies (3). 
We now prove (2). 
    The composite group homomorphism 
    \[ \varphi: \Pic\,X\stackrel{(-\cdot \ell_1)}{\longrightarrow}\mathbb{Z} \longrightarrow \mathbb{Z}/(H_2\cdot\ell_1)\mathbb{Z} \]
    is surjective, because so is each map (Theorem \ref{t-ex-surje}).  
    It suffices to show that $\Ker(\varphi) = \mathbb{Z}H_1\oplus \mathbb{Z}H_2$. 
The inclusion $\Ker(\varphi) \supset \mathbb{Z}H_1\oplus \mathbb{Z}H_2$ is clear. 
It is enough to show the opposite inclusion. 
Take $H \in \Ker(\varphi)$. 
Then $H\cdot \ell_1 \in (H_2\cdot \ell_1)\mathbb{Z}$, and hence there exists $b \in \Z$ such that 
    \[ 
    H\cdot \ell_1= b (H_2\cdot \ell_1).
    \]
    Thus we obtain $H-bH_2\in \mathrm{Ker}(-\cdot\ell_1) =\mathrm{Im}(f^*_1)=\mathbb{Z}H_1$ by Theorem \ref{t-ex-surje}. 
    Then we can write $H\sim c H_1+b H_2$ for some $c \in \mathbb{Z}$. 
    Thus (2) and (3) hold. 

Let us show (4). 
Since  $H_1$ and $H_2$ are linearly independent over $\mathbb{R}$ by (1), 
$\{H_1, H_2\}$ is a basis of $\mathrm{N}^1(X)_\mathbb{R}$. 
Then  we can write 
  \[
  -K_X\equiv \beta_1 H_1 + \beta_2 H_2
  \]
  for some $\beta_1, \beta_2 \in \mathbb{R}$. 
By $\mu_1=-K_X\cdot\ell_1$ and $H_2 \cdot \ell_1 =a$, we get 
  \begin{align*}
    \mu_1=-K_X\cdot\ell_1=\beta_1(H_1\cdot\ell_1)+\beta_2(H_2\cdot\ell_1) =\beta_2a.
  \end{align*}
  By symmetry, we obtain $\mu_2=\beta_1a$. Hence (4) holds.  
  Finally, (5) follows from  (4) and $(-K_X) \cdot c_2(X) =24$ (Lemma~\ref{firstlem}(3)). 
\end{proof}

\begin{lem}\label{l-degB}
We use Notation \ref{n-pic-2}. 
Assume that $R_1$ is of type $E_1$. 
Then 
\[
\deg B_1 \leq \left(r_1 - \frac{\mu_2}{a}\right)^2 \cdot L_1^3
\]
for $a :=  |\Pic(X) / (\mathbb{Z}H_1 \oplus \mathbb{Z}H_2)| \in \Z_{>0}$ (cf. Lemma \ref{l-basis}(1)). 
\end{lem}

\begin{proof}
Since $H_1$ and $H_2$ are nef, we get $H_1 \cdot H_2^2 \geq 0$. 
By $-aK_X \sim \mu_2 H_1 + H_2$ (Lemma \ref{l-basis}(4)) and $K_X -D_1 = f_1^*K_{Y_1} \sim -r_1 H_1$, the following holds: 
\begin{eqnarray*}
0 \leq H_1 \cdot H_2^2 
&=& H_1 \cdot (-aK_X -\mu_2H_1)^2\\
&=& H_1 \cdot ( (ar_1-\mu_2)H_1-aD_1)^2\\
&=& (ar_1-\mu_2)^2 H^3_1 -2a(ar_1-\mu_2) H^2_1 \cdot D_1 +a^2 H_1 \cdot D_1^2. 
\end{eqnarray*}
Then the required inequality 
$\deg B_1 \leq \left(r_1 - \frac{\mu_2}{a}\right)^2 \cdot L_1^3$ 
follows from $H_1^3 = L_1^3, H_1^2 \cdot D_1 =0,$ and $H_1 \cdot D_1^2 = -\deg B_1$.
\end{proof}

\begin{prop}\label{p-not-r1}
Let $X$ be a Fano threefold with $\rho(X)=2$. 
Let $R$ be an extremral ray of type $E_1$ and let $f: X \to Y$ be the contraction of $R$. 
Then $Y$ is a Fano threefold of index $r_Y \geq 2$. 
\end{prop}

\begin{proof}
Set $R_1 :=R$ and  $a :=  |\Pic(X) / (\mathbb{Z}H_1 \oplus \mathbb{Z}H_2)| \in \Z_{>0}$ (Lemma \ref{l-basis}(1)).
We use Notation \ref{n-pic-2}. 
By $-K_Y = f_*(-K_X)$, $-K_Y$ is big. 
This, together with $\rho(Y)=1$, implies that $-K_Y$ is ample. 
Suppose ${\cred r_1} = r_Y=1$. 
It is enough to derive a contradiction. 
Recall that $|-K_Y|$ is not very ample \cite[Theorem 1.1{\cred (3)}]{TanII}. 
By \cite[Theorem 1.1]{TanI}, we obtain $(-K_Y)^3 \in \{ 2, 4\}$. 
{\cred We have  $(-K_X)^3 \in 2\Z$ by the Riemann-Roch theorem $\chi(X, -K_X) = \frac{(-K_X)^3}{2}+3$ \cite[Corollary 2.6]{TanI}.}
{\cred These, together with 
$0 < (-K_X)^3 < (-K_Y)^3$ (Lemma \ref{l-(-K)^3-drops}), imply} that $(-K_Y)^3 =4$ and $(-K_X)^3 =2$. 
Let $\varphi_{|-K_X|} : X \dashrightarrow \P^{h^0(X, -K_X)-1}$ be  the rational map induced by the complete linear system $|-K_X|$. 
If $\dim (\Im\,\varphi_{|-K_X|}) \geq 2$, then our condition $\rho(X)=2$ contradicts \cite[Proposition 6.8]{TanI}. 
Hence we may assume that $\dim (\Im\,\varphi_{|-K_X|}) =1$. 
In this case, the mobile part $M$ of $|-K_X|$ is base point free and 
the induced morphism $\varphi_{|M|} : X \to \P^1$ is 
the contraction of the extremal ray $R_2$ of type D 
with $d_2 := (-K_X)^2 \cdot H_2 =1$ \cite[Proposition 4.1]{TanI} 
{\cred (note that the equality $-K_X^2 \cdot E=1$ claimed in \cite[Proposition 4.1(3)]{TanI} is a typo; the correct equality is $(-K_X)^2 \cdot E=1$, as given in the last line of its proof)}. 
In particular, $\mu_2 =1$. 
Then Lemma \ref{c_times_H} and Lemma \ref{l-basis} imply 
\[
24 a = \frac{24}{r_1} + \deg B_1 + (12 - d_2) = 35+ \deg B_1. 
\]
By $\deg B_1 \leq (r_1 - \frac{\mu_2}{a})^2 \cdot L_1^3 = (1 -\frac{1}{a})^2 \cdot  4<4$ (Lemma \ref{l-degB}), 
we obtain $24a =35 + \deg B_1 < 35 +4$, which implies $a=1$. 
However, this leads to the following contradiction: 
$\deg B_1 \leq (1 -\frac{1}{a})^2 \cdot  4=0$. 
\end{proof}

\begin{prop}\label{basis}
We use Notation \ref{n-pic-2}. 
Then the following hold. 
\begin{enumerate}
\item  $\{H_1, H_2\}$ is a $\mathbb{Z}$-linear basis of $\Pic\,X \simeq \Z^{\oplus 2}$. 
\item $\{\ell_1, \ell_2\}$ is a dual basis of $\mathrm{N}_1(X)_\mathbb{Z}$, i.e., 
$H_1 \cdot \ell_1 = H_2 \cdot \ell_2 =0$ and $H_1 \cdot \ell_2 = H_2 \cdot \ell_1=1$. 
Here we set $\mathrm{N}_1(X)_\mathbb{Z} := {\rm Z}_1(X)/ \equiv$, 
where $ {\rm Z}_1(X) := \bigoplus_{C} \Z C$ is the free $\Z$-module that is freely generated by all the curves $C$ on $X$ and $\equiv$ denotes the numerical equivalence. 
\item  $-K_X \sim \mu_2 H_1 + \mu_1 H_2.$
\item $24 = \mu_2 H_1 \cdot c_2(X) + \mu_1 H_2 \cdot c_2(X)$. 
\end{enumerate}
\end{prop}

\begin{proof}[Proof of Proposition \ref{basis} for the case when $X$ is primitive] 
Set 
$a :=  |\Pic(X) / (\mathbb{Z}H_1 \oplus \mathbb{Z}H_2)| \in \Z_{>0}$ (Lemma \ref{l-basis}(1)). 
Since $X$ is primitive, there is an extremal ray of type $C$ (Proposition \ref{hasC}). 
In what follows, we assume that $R_1$ is of type $C$.  
The proof consists of 7 steps. 

\setcounter{step}{0}

\begin{step}\label{s4-basis}
In order to complete the proof of Proposition \ref{basis}, 
it is enough to show that $a=1$. 
\end{step}

\begin{proof}[Proof of Step \ref{s4-basis}]
Assume $a=1$. 
Then the assertion (1) follows from Lemma \ref{l-basis}(1)(2). 
The assertions (3) and (4) hold by Lemma \ref{l-basis}(4)(5). 
The assertion (2) follows from Lemma \ref{l-basis}(3) and $H_1 \cdot \ell_1 = H_2 \cdot \ell_2=0$. 
This completes the proof of Step \ref{s4-basis}. 
\end{proof}




\begin{step}\label{s5-basis}
If $(\mu_1, \mu_2)=(1, 1)$, then $a=1$. 
\end{step}

\begin{proof}[Proof of Step \ref{s5-basis}]
By Lemma \ref{l-basis}(5), we obtain 
\[
24a= H_1\cdot c_2(X) + H_2\cdot c_2(X).
\]
  Since $R_1$ is of type $C_1$, the following holds (Lemma \ref{c_times_H2}(1)): 
  \[7\le c_2(X)\cdot H_1 \le 17.\]
  If $c_2(X)\cdot H_2=45$, then $52 \le 24a\le 62$,  which is a contradiction.
  Hence $1 \leq c_2(X)\cdot H_2 \le 24$  (Lemma \ref{c_times_H2}(2)), which implies 
  \[
  24a = H_1\cdot c_2(X) + H_2\cdot c_2(X) \leq 17+24 = 41.
  \]
  Thus $a=1$.
This completes the proof of Step \ref{s5-basis}. 
\end{proof}

\begin{step}\label{s6-basis}
If $(\mu_1, \mu_2)=(1, 2)$, then $a=1$. 
\end{step}

\begin{proof}[Proof of Step \ref{s6-basis}]
By Lemma \ref{l-basis}(5), we obtain 
\[
24a= 2 H_1\cdot c_2(X) + H_2\cdot c_2(X).
\]  
  In this case, $a$ is odd.
  Indeed, if $a$ is even, then 
  the linear equivalence 
  \[
  a(-K_X)\sim 2H_1+H_2
  \]
  would imply that $H_2$ is not a primitive element, which is a contradiction (Lemma \ref{l-Hi-prim}). 
Since $R_1$ is of type $C_1$ and $R_2$ is not of type $E_5$ by $\mu_2 =2$, 
we get $c_2(X)\cdot H_1 \le 17$ and $c_2(X)\cdot H_2 \le 24$ (Lemma \ref{c_times_H2}), which imply 
  \[
  24a= 2 H_1\cdot c_2(X) + H_2\cdot c_2(X)\le 2\cdot 17+24=58. 
  \]
  As $a$ is odd, we get $a=1$. 
This completes the proof of Step \ref{s6-basis}. 
\end{proof}

\begin{step}\label{s7-basis}
If $(\mu_1, \mu_2)=(1, 3)$, then $a=1$. 
\end{step}

\begin{proof}[Proof of Step \ref{s7-basis}]
In this case, $R_1$ is of type $C_1$ and $R_2$ is of type $D_3$. 
By Lemma \ref{l-basis}(5), we obtain 
\[
24a= 3 H_1\cdot c_2(X) + H_2\cdot c_2(X).
\]

We now show that $a$ is odd.
  By $-aK_X \sim 3H_1+H_2$ (Lemma \ref{l-basis}(4)),  we get 
  \begin{align*}
    a^3(-K_X)^3 &=27H_1^3+27 H_1^2\cdot H_2+9H_1\cdot H_2^2+H_2^3 \\
    &\overset{{\rm (a)}}{=} 27\cdot 2(H_2\cdot \ell_1) \\
    &\overset{{\rm (b)}}{=} 27\cdot 2a, 
  \end{align*}
where (b) holds by Lemma \ref{l-basis}(3) and (a) follows from 
$H^3_1=0$, $H_2^2\equiv 0$, and $H_1^2 \equiv 2\ell_1$ (Lemma~\ref{pic2-numequiv}). 
Hence we obtain $a^2 (-K_X)^3 = 2 \cdot 27$. 
By $(-K_X)^3 \in 2\Z$, we get $a\not\in 2\Z$. 

We have $c_2(X)\cdot H_1\le 17$ (Lemma \ref{c_times_H2}) and $c_2(X)\cdot H_2 =3$  (Lemma \ref{c_times_H}), which imply 
  \[
  24a = 3 H_1\cdot c_2(X) + H_2\cdot c_2(X) \le 3\cdot 17+3=54. 
  \]
  As $a$ is odd, we get $a=1$. This completes the proof of Step \ref{s7-basis}. 
\end{proof}

In what follows, we treat the case when $\mu_1=2$. 
In this case, $R_1$ is of type $C_2$ and $H_1 \cdot c_2(X)=6$ (Lemma \ref{c_times_H}).

\begin{step}\label{s8-basis}
If $(\mu_1, \mu_2)=(2, 1)$, then $a=1$. 
\end{step}

\begin{proof}[Proof of Step \ref{s9-basis}]
By Lemma \ref{l-basis}(5), we obtain 
\[
24a= H_1\cdot c_2(X) + 2H_2\cdot c_2(X) =6 + 2H_2\cdot c_2(X).
\]  
By the same argument as in 
the case $(\mu_1, \mu_2)=(1, 2)$ (Step \ref{s6-basis}), 
we see that $a$ is odd. 
  If $c_2(X)\cdot H_2=45$, then 
  \[24a=6+2\cdot 45=96,\]
  which contradicts the fact that $a$ is odd.
  Hence $c_2(X)\cdot H_2 \le 24$ (Lemma \ref{c_times_H2}), which implies  
  \[24a=6 + 2H_2\cdot c_2(X) \le 6+2\cdot 24=54.\]
  As $a$ is odd, we get $a=1$. 
This completes the proof of Step \ref{s8-basis}. 
\end{proof}

\begin{step}\label{s9-basis}
If $(\mu_1, \mu_2)=(2, 2)$, then $a=1$. 
\end{step}

\begin{proof}[Proof of Step \ref{s9-basis}]
By Lemma \ref{l-basis}(5), we obtain 
\[
24a= 2 H_1\cdot c_2(X) + 2H_2\cdot c_2(X) =12+ 2H_2\cdot c_2(X).
\]  
Note that the type of $R_2$ is $E_2, C_2$, or $D_2$ (Remark \ref{r-length}). 
 If $R_2$ is of type $E_2$, then $c_2(X)\cdot H_2=24/r$ for some $r\in \mathbb{Z}_{> 0}$ (Lemma \ref{c_times_H}).
  Hence 
  \[12a=6+ H_2\cdot c_2(X)=6+\frac{24}{r}.\]
By $24/r \in 6 \Z$, we get $r \in \{1, 2, 4\}$. 
Then it is easy to see that $(a, r)=(1, 4)$ is a unique solution of $12a = 6+\frac{24}{r}$.

We may assume that $R_2$ is of type $C_2$ or $D_2$. 
In this case, we have $c_2(X)\cdot H_2 \le 6$ (Lemma \ref{c_times_H}), which implies  
  \[
  12a = 6+ H_2\cdot c_2(X)\le 12.
  \]
Thus $a=1$. 
This completes the proof of Step \ref{s9-basis}. 
\end{proof}

\begin{step}\label{s10-basis}
If $(\mu_1, \mu_2)=(2, 3)$, then $a=1$. 
\end{step}

\begin{proof}[Proof of Step \ref{s10-basis}]
In this case, $R_2$ is of type $D_3$. 
By Lemma \ref{l-basis}(5) and Lemma \ref{c_times_H}, we obtain 
\[
24a= 3 H_1\cdot c_2(X) + 2H_2\cdot c_2(X) =18+ 6 = 24.
\]  
Hence $a=1$. 
This completes the proof of Step \ref{s10-basis}. 
\end{proof}
Step \ref{s4-basis}--Step \ref{s10-basis} complete the proof of Proposition \ref{basis} for the case when $X$ is primitive. 
\end{proof}

\begin{proof}[Proof of  Proposition \ref{basis} for the case when $X$ is imprimitive] 
Since $X$ is imprimitive, $X$ has an extremal ray of type $E_1$. 
In what follows, we assume that $R_1$ is of type $E_1$.  
Set $a :=  |\Pic(X) / (\mathbb{Z}H_1 \oplus \mathbb{Z}H_2)| \in \Z_{>0}$ (Lemma \ref{l-basis}(1)). 
It is enough to show $a=1$ (cf. Step \ref{s4-basis} of the primitive case). 
{\cred Recall that $2 \leq r_1 \leq 4$ (Proposition \ref{p-not-r1}).} 

\begin{claim}\label{cl-imprim}
\begin{enumerate}
\item $H_1 \cdot c_2(X) \leq 31$. 
\item $H_1 \cdot c_2(X) \leq 
\begin{cases}
25 \qquad (\text{if } a \leq 3)\\ 
23 \qquad (\text{if } a \leq 2)\\
17 \qquad (\text{if } a=1).
\end{cases}$
\item $H_1 \cdot c_2(X) \leq 20$ if $a \leq 3$ and $\mu_2 =2$. 
\end{enumerate}
\end{claim}

\begin{proof}[Proof of Claim \ref{cl-imprim}]
By Lemma \ref{c_times_H} and Lemma \ref{l-degB}, 
we have 
\begin{equation}\label{e-imprim-Hc2-bound}
H_1 \cdot c_2(X) =\frac{24}{r_1} +\deg B_1 \leq \frac{24}{r_1} + \left( r_1 - \frac{\mu_2}{a}\right)^2 \cdot L_1^3. 
\end{equation}
In what follows, we shall separate the argument depending on the value of $r_1$. 


Assume $r_1 =4$. In this case, we have 
{\cred $Y_1 \simeq \P^3$ and} $L_1^3 = 1$ {\cred \cite[Theorem 2.18]{TanI}}. 
This, together with (\ref{e-imprim-Hc2-bound}), implies 
\[
H_1 \cdot c_2(X) \leq 6 +  \left( 4 - \frac{\mu_2}{a}\right)^2< 6 + 4^2 =22. 
\]
If $a=1$, then $H_1 \cdot c_2(X) \leq 6 +  (4-1)^2 =15$. 
If $a \leq 3$ and $\mu_2 =2$, then $H_1 \cdot c_2(X) \leq 6 + (4-\frac{2}{3})^2 = 6+ \frac{100}{9} <6+12 =18$.

Assume $r_1 =3$. In this case, 
{\cred $Y_1$ is a smooth quadric hypersurace in $\P^4$ and} 
we have $L_1^3 = 2$ {\cred \cite[Theorem 2.18]{TanI}}. 
This, together with (\ref{e-imprim-Hc2-bound}), implies 
\[
H_1 \cdot c_2(X) \leq 8 +  2\left( 3 - \frac{\mu_2}{a}\right)^2< 8 + 2 \cdot 3^2 =28. 
\]
If $a \leq 3$, then $H_1 \cdot c_2(X) \leq 8 +  2\left( 3 - \frac{1}{3}\right)^2 = 8+ \frac{128}{9} < 8+\frac{135}{9}=8+15=23$. 
If $a \leq 2$, then $H_1 \cdot c_2(X) \leq 8 +  2\left( 3 - \frac{1}{2}\right)^2 = 8+ \frac{25}{2} <21$. 
If $a = 1$, then $H_1 \cdot c_2(X) \leq 8 +  2\left( 3 - 1\right)^2 =16$.  
If $a \leq 3$ and $\mu_2 =2$, then $H_1 \cdot c_2(X) \leq 8 +  2\left( 3 - \frac{2}{3}\right)^2 = 8 + \frac{98}{9} < 8+ \frac{99}{9} =19$.

Assume $r_1 =2$. In this case, we have $1 \leq L_1^3 \leq 5$ 
{\cred \cite[Remark 4.9]{TanII}}. 
This, together with (\ref{e-imprim-Hc2-bound}), implies 
\[
H_1 \cdot c_2(X) \leq 12 +  5 \left( 2 - \frac{\mu_2}{a}\right)^2< 12 +5 \cdot 4 =32. 
\]
If $a \leq 3$, then $H_1 \cdot c_2(X) \leq 12 +  5\left( 2 - \frac{1}{3}\right)^2 =12 + \frac{125}{9}
<12 +\frac{126}{9} = 26$. 
If $a \leq 2$, then $H_1 \cdot c_2(X) \leq 12 +  5\left( 2 - \frac{1}{2}\right)^2 = 12 + \frac{45}{4} <24$. 
If $a  = 1$, then $H_1 \cdot c_2(X) \leq 12 +  5\left( 2 - 1\right)^2  = 17$. 
If $a \leq 3$ and $\mu_2 =2$, then $H_1 \cdot c_2(X) \leq 12 +  5\left( 2 - \frac{2}{3}\right)^2 = 12 + \frac{80}{9} <21$. 
This completes the proof of Claim \ref{cl-imprim}. 
\qedhere



\end{proof}

We first treat the case when $R_2$ is of type $E_1$. 
In this case, the same conclusion of Claim \ref{cl-imprim} holds for $H_2 \cdot c_2(X)$. 
Claim \ref{cl-imprim}(1) implies 
\[
24 a =  H_1 \cdot c_2(X) + H_2 \cdot c_2(x) \leq 31+31 =62. 
\]
Hence we may assume that $a = 2$. 
By Claim \ref{cl-imprim}(2), we get 
\[
48= 24 a  =  H_1 \cdot c_2(X) + H_2 \cdot c_2(x) \leq 23+23 =46. 
\]
This is absurd. Thus $R_2$ is not of type ${\cred E_1}$. 

\medskip

We now consider the case when $H_2 \cdot c_2(X) =45$. 
Then $R_2$ is of type $E_5$ (Lemma \ref{c_times_H}, Lemma \ref{c_times_H2}), and hence $\mu_2 =1$. Then Claim \ref{cl-imprim}(1) implies 
\[
24a =  H_1 \cdot c_2(X) + H_2 \cdot c_2(X) = H_1 \cdot c_2(X) +45 \leq 31+45=76. 
\]
Hence $a \leq 3$. 
By Claim \ref{cl-imprim}(2), we obtain 
\[
24a =  H_1 \cdot c_2(X) + H_2 \cdot c_2(X) = H_1 \cdot c_2(X) +45 \leq 25+45=70. 
\]
Hence we may assume that $a =2$. 
We then get $48 =24a = H_1 \cdot c_2(X) +45$, which implies 
\[
\frac{24}{r_1}<\frac{24}{r_1} +\deg B = H_1 \cdot c_2(X) =3.
\]
This contradicts  $r_1 \in \{2, 3, 4\}$. 
Thus $H_2 \cdot c_2(X) \neq 45$, and hence $1 \leq H_2 \cdot c_2(X) \leq 24$ (Lemma \ref{c_times_H2}). 

\medskip

Assume that $\mu_2 =1$. 
Then Claim \ref{cl-imprim}(1) implies 
\[
24 a = H_1 \cdot c_2(X) + H_2 \cdot c_2(X) \leq 31 + 24 = 55. 
\]
Hence we may assume that $a= 2$. By Claim \ref{cl-imprim}(2), we get  
\[
48 = 24 a = H_1 \cdot c_2(X) + H_2 \cdot c_2(X) \leq 23 + 24=47, 
\]
which is absurd.

\medskip

Assume that $\mu_2 =2$. 
Then Claim \ref{cl-imprim}(1) implies 
\[
24 a = 2H_1 \cdot c_2(X) + H_2 \cdot c_2(X) \leq 2 \cdot 31 + 24 = 86. 
\]
Hence $a \leq 3$. By Claim \ref{cl-imprim}(3), we get 
\[
24 a = 2H_1 \cdot c_2(X) + H_2 \cdot c_2(X) \leq 2 \cdot 20 + 24 = 64. 
\]
Thus we may assume that $a=2$. 
We then get 
\[
-2K_X \sim 2H_1 + H_2, 
\]
which contradicts the fact that $H_2$ is primitive (Lemma \ref{l-Hi-prim}).

\medskip

Assume that $\mu_2 =3$. 
Then $R_2$ is of type $D_3$ and $H_2 \cdot c_2(X) =3$ (Lemma \ref{c_times_H}). 
Therefore, $24 a = 3H_1 \cdot c_2(X) + H_2 \cdot c_2(X)$ 
(Lemma \ref{l-basis}) implies the following (Claim \ref{cl-imprim}(1)): 
\[
8a = H_1 \cdot c_2(X) + 1 \leq 32.
\]
By $-aK_X \sim 3H_1 + H_2$, we obtain $a \not\in 3\Z$. 
Hence we may assume that $a \in \{2, 4\}$.

Suppose $a=4$. Then Lemma \ref{c_times_H} and Lemma \ref{l-degB} imply 
\[
31 = H_1 \cdot c_2(X) =  \frac{24}{r_1} + \deg B_1 \leq \frac{24}{r_1} +
\left(r_1-\frac{3}{4}\right)^2 L_1^3. 
\]
If $r_1 =4$, then $L_1^3 =1$,  and hence
$31 \leq 6 + (4-\frac{3}{{\cred 4}})^2<6+16=22$, which is absurd. 
If $r_1 =3$, then $L_1^3 =2$,  and hence
$31 \leq 8 + (3-\frac{3}{{\cred 4}})^2 \cdot 2<8+9 \cdot 2 =26$, which is impossible. 
If $r_1 =2$, then $1 \leq L_1^3 \leq 5$, and 
hence $31 \leq 12 + (2-\frac{3}{4})^2 \cdot L_1^3 
\leq 12 + \frac{25}{16} \cdot 5 <12 +8$, which is a contradiction. 
We then obtain $a \neq 4$. 

Suppose $a=2$. 
It suffices to derive a contradiction. 
We get 
$16 =  H_1 \cdot c_2(X) + 1$. 
Hence Lemma \ref{c_times_H} and Lemma \ref{l-degB} imply 
\[
15 =H_1 \cdot c_2(X) = \frac{24}{r_1} + \deg B_1 \leq \frac{24}{r_1} + 
\left(r_1 -\frac{3}{2}\right)^2 \cdot L_1^3. 
\]
If $r_1 =4$, then $L_1^3=1$ and $15 \leq 6 + \frac{25}{4} \cdot 1 <6+7$, which is absurd. 
If $r_1 =3$, then $L_1^3 =2$ and $15 \leq 8 + \frac{9}{4} \cdot 2 < 8+5$, which is a contradiction. 
If $r_1 =2$, then $1 \leq L_1^3 \leq 5$ and 
$15 \leq 12 + \frac{1}{4} \cdot L_1^3 \leq 12 + \frac{5}{4} <12 +2$, which is absurd. 
This completes the proof of Proposition \ref{basis}. 
\end{proof}

\begin{lem}\label{l-E1-H1H2}
We use Notation \ref{n-pic-2}. 
Assume that $R_1$ is of type $E_1$. 
Then the following hold. 
\begin{enumerate}
\item $H_1^2 \cdot H_2 = (r_1-\mu_2)L_1^3$. 
\item $H_1 \cdot H_2^2 = (r_1-\mu_2)^2 L_1^3 -\deg B_1$. 
\end{enumerate}
\end{lem}

\begin{proof}
Let us show (1). 
Set $D_1 := \Ex(f_1)$. 
We obtain 
\[
H_1^2 \cdot H_2 \overset{{\rm (i)}}{=} H_1^2 \cdot ( -K_X-\mu_2 H_1) \overset{{\rm (ii)}}{=}   H_1^2 \cdot ( (r_1-\mu_2)H_1 -D_1) 
= (r_1 -\mu_2) L_1^3, 
\]
where (i) follows from $-K_X \sim \mu_2 H_1 + H_2$ (Proposition \ref{basis}(3)) 
and (ii) holds by $K_X =f_1^*K_{X_1} +D_1 = -r_1H_1 + D_1$. 
Thus (1) holds. 
The assertion (2) follows from the displayed equation 
in the proof of Lemma \ref{l-degB} by using $H_1^2 \cdot D_1 =0$ and $a=1$ (Proposition \ref{basis}). 
\end{proof}


\subsection{Case $C-C$}\label{ss-rho2-CC}

\begin{nota}\label{n-pic-2-C}
We use Notation \ref{n-pic-2}. 
Assume that $R_1$ is of type $C$ and $R_2$ is of type $C$. 
In particular, for each $i \in \{1, 2\}$, 
$f_i : X \to Y_i = \mathbb P^2$ is of type $C$ 
{\cred (Theorem \ref{pic-over2})}. 
Set 
\[
f:=f_1 \times f_2 \colon X\to \mathbb{P}^2\times \mathbb{P}^2.
\]
Let $\widetilde f : X \to f(X)$ be the induced morphism, where $f(X)$ denotes the scheme-theoretic image. 
For the $i$-th projection $\pi_i : \mathbb{P}^2\times \mathbb{P}^2 \to \P^2$, set $M_i := \pi_i^*L_i$. 
Note that $\{M_1, M_2\}$ is a free $\Z$-linear basis of $\Pic (\mathbb{P}^2\times \mathbb{P}^2)$. 
\end{nota}

\begin{lem}\label{f_push_div}
We use Notation \ref{n-pic-2-C}. 
Then 
\[
\deg(\widetilde{f})\cdot f(X)=f_*(X) \sim \frac{2}{\mu_2}M_1+\frac{2}{\mu_1}M_2. 
\]
In particular, $\deg(\widetilde{f})=1$ or $\deg(\widetilde{f})=2$. 
\end{lem}

\begin{proof}
Recall that $f : X \to Y_1 \times Y_2 = \P^2 \times \P^2$ is a finite morphism (Lemma \ref{l-finite-morph}). 
  We can write
  \[f_*(X)\sim a_1M_1+a_2M_2\]
  for some $a_1,a_2\in \mathbb{Z}$.
  Pick a curve $\Gamma_1=(\mathrm{point}\times \mathrm{line})$ on $\mathbb{P}^2\times \mathbb{P}^2$.
  Since $M_1\sim (\mathrm{line}\times \mathbb{P}^2)$ and $M_2 \sim (\mathbb{P}^2\times \mathrm{line})$, we have 
  \[ 
  M_1\cdot \Gamma_1=0\qquad\text{and} \qquad M_2\cdot \Gamma_1=1.
  \]
  Hence we obtain $f_*(X)\cdot \Gamma_1 =a_2$.
  Now, 
  \begin{align*}
    M_1^2\cdot M_2 
    &\equiv (\mathrm{point}\times \mathbb{P}^2)\cdot (\mathbb{P}^2\times \mathrm{line}) \\
    &\equiv (\mathrm{point} \times \mathrm{line}) \\
    &\equiv \Gamma_1.
  \end{align*}
We have $H_2\cdot \ell_1=1$ (Proposition~\ref{basis}(2)) 
and $H_1^2\equiv \frac{2}{\mu_1}\ell_1$ (Lemma \ref{pic2-numequiv}).
  We obtain 
  \begin{align*}
    f_*(X)\cdot \Gamma_1 &= f_*(X)\cdot M_1^2\cdot M_2\\
    &= (f^*M_1)^2\cdot f^*M_2\\
    &= H_1^2\cdot H_2 \\
    &= \frac{2}{\mu_1}H_2\cdot \ell_1 \\
    &= \frac{2}{\mu_1}.
  \end{align*}
  Hence $a_2=\frac{2}{\mu_1}$.
  Similarly, we have $a_1=\frac{2}{\mu_2}$.
\end{proof}

\begin{lem}\label{pic2-r2C-K^3}
We use Notation \ref{n-pic-2-C}. Then   $(-K_X)^3=6(\mu_1^2+\mu_2^2)$.
\end{lem}

\begin{proof}
  By Proposition~\ref{basis}, we have $-K_X\sim \mu_2 H_1 + \mu_1 H_2$.
  Since $L_i$ is a line on $Y_i = \mathbb{P}^2$, 
  \[H_i^3=f_i^*L_i^3=0.\]
  By $H_i^2 \equiv \frac{2}{\mu_i}\ell_i$ (Lemma~\ref{pic2-numequiv}) and $H_1 \cdot \ell_2 = H_2 \cdot \ell_1 = 1$ (Proposition \ref{basis}), we obtain  
  \begin{align*}
    (-K_X)^3 &= (\mu_2 H_1 + \mu_1  H_2)^3 \\
    &= \mu_2 ^3 H_1^3 + 3\mu_2^2\mu_1 H_1^2\cdot H_2 + 3 \mu_2\mu_1^2 H_1 \cdot H_2^2 + \mu_1^3 H_2^3 \\
    &= 3\mu_2^2\mu_1\cdot \frac{2}{\mu_1}(H_2\cdot \ell_1) + 3\mu_2\mu_1^2\cdot \frac{2}{\mu_2}(H_1\cdot \ell_2) \\
    &= 6(\mu_1^2 + \mu_2^2).
  \end{align*}
\end{proof}

\subsubsection{Case $\deg (\widetilde f)=1$}

\begin{lem}\label{deg1-isom}
We use Notation \ref{n-pic-2-C}. 
Assume that $\deg \widetilde{f} =1$. 
Then the following hold. 
\begin{enumerate}
\item 
$\widetilde f : X\to f(X)$ is an isomorphism.
\item 
If $(\mu_1, \mu_2)=(1, 1)$, then $X \subset \mathbb P^2 \times \mathbb P^2$ is of bidegree $(2, 2)$. 
\item 
If $(\mu_1, \mu_2)=(1, 2)$, then $X \subset \mathbb P^2 \times \mathbb P^2$ is of bidegree $(1, 2)$. 
\item
If $(\mu_1, \mu_2)=(2, 2)$, then $X \subset \mathbb P^2 \times \mathbb P^2$ is of bidegree $(1, 1)$. 
\end{enumerate}
\end{lem}

\begin{proof}
If (1) holds, then (2), (3), and (4) hold. 
Thus it suffices to show (1). 

  Set $Y := f(X)$. 
Since $\widetilde f : X \to Y$ is a finite birational morphism from a smooth variety $X$ to a variety $Y$, $\widetilde f$ is nothing but the normalisation of $Y$. 
As $Y$ is a prime divisor on $\P^2 \times \P^2$, $Y$ is Gorenstein. 
By \cite[Proposition 2.3]{Rei94}, we have 
\[
\omega_X \otimes \MO_X(C) \simeq \widetilde f^*\omega_Y 
\]
for the conductor $C$, which is an effective Cartier divisor on $X$ whose support coincides with the non-isomorphic 
locus $\Ex(\widetilde f)$ of $\widetilde f$. 
By the adjunction formula, 
we have 
  \begin{align*}
    \omega_Y &\simeq 
(\omega_{\mathbb{P}^2\times \mathbb{P}^2} \otimes \MO_{\mathbb{P}^2\times \mathbb{P}^2}(Y))\otimes \MO_Y\\
    &\simeq 
\left( \MO_{\mathbb{P}^2\times \mathbb{P}^2}(-3,-3)\otimes \MO_{\mathbb{P}^2\times \mathbb{P}^2}\left(\frac{2}{\mu_2},\frac{2}{\mu_1}\right)\right)\otimes \MO_Y \\
    &\simeq \MO_{\mathbb{P}^2\times \mathbb{P}^2}\left(\frac{2}{\mu_2}-3,\frac{2}{\mu_1}-3\right)\otimes \MO_Y, 
  \end{align*}
  where the second isomorphims holds by Lemma \ref{f_push_div}. 
Therefore, we obtain 
  \begin{align*}
K_X+ C \sim \left(\frac{2}{\mu_2}-3\right)H_1+\left(\frac{2}{\mu_1}-3\right)H_2.
  \end{align*}
It follows from $-K_X\sim \mu_2 H_1+\mu_1 H_2$ (Proposition~\ref{basis}) that 
  \[
C \sim \left(\mu_2 + \frac{2}{\mu_2}-3\right)H_1+\left(\mu_1 + \frac{2}{\mu_1}-3\right)H_2.\]
As we have $\mu_i\in\{1,2\}$ for each $i$, we obtain $C\sim 0$, which implies $C=0$.
Hence $\tilde{f}$ is an isomorphism. 
\end{proof}

\subsubsection{Case $\deg (\widetilde f)=2$} 

\begin{lem}\label{l-bideg11}
Let $W$ be a nonzero effective Cartier divisor on $\P^2 \times \P^2$ 
of bidegree $(1, 1)$. 
Then the following hold. 
\begin{enumerate}
\item $W$ is isomorphic to one of the following: 
\begin{itemize}
\item $W_1:=\{x_0y_0 =0\} \subset \P^2 \times \P^2$. 
\item  $W_2 :=\{x_0y_0 +x_1y_1=0\} \subset \P^2 \times \P^2$. 
\item  $W_3 :=\{x_0y_0 +x_1y_1+x_2y_2=0\} \subset \P^2 \times \P^2$. 
\end{itemize}
Note that $W_1$ is not irreducible, $W_2$ is not smooth but irreducible, 
and $W_3$ is smooth. 
\item 
If $W \simeq W_2$, then there exists a closed point $Q \in \P^2$ such that 
$\dim q_1^{-1}(Q) =2$, where 
$q_1 : W \hookrightarrow \P^2 \times \P^2 \xrightarrow{{\rm pr}_1} \P^2$. 
\end{enumerate}
\end{lem}

\begin{proof}
Let us show (1). 
Since $W$ is defined by 
\[
f \in H^0(\P^2 \times \P^2, \MO_{\P^2 \times \P^2}(1, 1)) \simeq 
H^0(\P^2, \MO_{\P^2}(1)) \otimes_k H^0(\P^2, \MO_{\P^2}(1)),
\]
we have 
\[
W = \left\{ \sum_{0 \leq i,j\leq 2} a_{ij}x_iy_j =0 \right\} \subset \P^2 \times \P^2. 
\]
Applying the elementary operations to the $3 \times 3$-matrix $(a_{ij})$, 
we may assume that $W$ is one of $W_1, W_2$, and $W_3$ (note that we can apply row and column elementary operations). 
Thus (1) holds.

Let us show (2). 
By the proof of (1), $W \simeq W_i$ is obtained by applying 
a suitable ${\rm Aut}(\P^2) \times {\rm Aut}(\P^2)$-action 
on $\P^2 \times \P^2$. 
Hence we may assume that $W=W_2$. 
Then the assertion (2) holds 
for $Q := [0:0:1]$. 
\end{proof}

\begin{lem}\label{l-CC-deg2}
We use Notation \ref{n-pic-2-C}. 
Assume that $\deg \widetilde{f} \neq 1$. 
Set $W :=f(X)$. 
Then the following hold. 
\begin{enumerate}
\item 
$\deg \widetilde{f} =2$ and $\mu_1=\mu_2=1$. 
\item $(-K_X)^3 = 12$. 
\item 
$W$ is a smooth prime divisor on $\mathbb P^2 \times \mathbb P^2$ 
of bidegree $(1, 1)$. 
\item 
For $\mathcal L := (\Coker( \MO_{W} \to \widetilde{f}_*\MO_X))^{-1}$, it holds that 
$\mathcal L \simeq \MO_{\mathbb P^2 \times \mathbb P^2}(1, 1)|_{W}$. 
\item $\widetilde f: X \to W$ is a split double cover. 
\end{enumerate}
\end{lem}

\begin{proof}
By Lemma \ref{f_push_div}, we have 
\[
\deg (\tilde{f})\cdot f(X) = f_*(X) \sim \frac{2}{\mu_2} M_1+\frac{2}{\mu_1}M_2.\] 
We then get $\deg (\tilde{f}) =2$ and $\mu_1 = \mu_2 =1$. 
Hence (1) holds. 
By Lemma~\ref{pic2-r2C-K^3}, $(-K_X)^3=6(\mu_1^2+\mu_2^2)=12.$ 
Thus (2) holds. 
Furthermore, we have 
\[
\MO_{\mathbb P^2 \times \mathbb P^2}(W) \simeq \MO_{\P^2 \times \P^2}(1, 1). 
\]
Then the assertion (3) follows from Lemma \ref{l-bideg11}(2), 
{\cred because $f_1$ commutes with the projection ${\rm pr}_1$, 
i.e., $f_1 : X \xrightarrow{\wt{f}} f(X) =W 
\hookrightarrow \P^2 \times \P^2 \xrightarrow{{\rm pr}_1}\P^2$ (Notation \ref{n-pic-2-C})}. 

Let us show (4). 
Pick a Cartier divisor $L$ on $W$ with $\mathcal L \simeq \MO_W(L)$. 
By Lemma \ref{l-dc-omega}, we have 
\begin{equation}\label{C1_C1_1}
  K_X\sim \widetilde{f}^*(K_{W}+L). 
\end{equation}
By the adjunction formula, we obtain 
\begin{align*}
  K_W &\sim (K_{\mathbb{P}^2\times \mathbb{P}^2}+W)|_W \\
  &\sim (-3M_1-3M_2+M_1+M_2)|_W \\
  &\sim -2(M_1+M_2)|_W.
\end{align*}
By (\ref{C1_C1_1}) and $K_X\sim -H_1-H_2$ (Lemma~\ref{basis}),  we obtain
\[
\widetilde{f}^*L \sim K_X - 
\widetilde{f}^*K_W \sim ( -H_1-H_2)-\widetilde{f}^*( -2(M_1+M_2)|_W)
=H_1+H_2. 
\]
It holds that $L \equiv c_1( \MO_{\mathbb P^2 \times \mathbb P^2}(1, 1)|_W)$. 
Since $W$ is a smooth Fano threefold, we obtain 
$L \sim c_1( \MO_{\mathbb P^2 \times \mathbb P^2}(1, 1)|_W)$. 
Thus (4) holds. 
The assertion (5) follows from (3), (4),  and Lemma \ref{l-split-criterion}(2). 
\end{proof}

\begin{thm}\label{t-pic2-C}
Let $X$ be a  Fano threefold with $\rho(X)=2$. 
Set $W := \{x_0y_0+x_1y_1 + x_2y_2 =0\} \subset \P^2 \times \P^2$, which is a smooth hypersurface on $\P^2 \times \P^2$ of bidegree $(1, 1)$. 
Let $R_1$ and $R_2$ be the distinct extremal rays of $\NE(X)$. 
Assume that each of $R_1$ and $R_2$ is of type $C$. 
Then, possibly after permuting $R_1$ and $R_2$, one and only one of the following holds. 
\begin{enumerate}
\item 
$R_1$ is of type $C_1$, $R_2$ is of type $C_1$, and one of the following holds {\cred (No.\ \hyperref[table-2-6]{2-6})}. 
\begin{enumerate}
\item 
$X$ is isomorphic to a hypersurface of $\P^2 \times \P^2$ of bidegree $(2, 2)$ and   $(-K_X)^3=12$.
\item 
There exists a split double cover $f: X \to W$ such that $(f_*\MO_X/\MO_W)^{-1} \simeq \MO_{\P^2 \times \P^2}(1, 1)|_W$ and   $(-K_X)^3=12$. 
\end{enumerate}
\item 
$R_1$ is of type $C_1$, $R_2$ is of type $C_2$, $(-K_X)^3=30$, and 
$X$ is isomorphic to a hypersurface of $\P^2 \times \P^2$ of bidegree $(1, 2)$ {\cred (No.\ \hyperref[table-2-24]{2-24})}. 
\item 
$R_1$ is of type $C_2$, $R_2$ is of type $C_2$, $(-K_X)^3=48$, and 
$X$ is isomorphic to $W$ {\cred (No.\ \hyperref[table-2-32]{2-32})}. 
\end{enumerate}
\end{thm}

\begin{proof}
If one of the contractions is not generically smooth, 
then (2) holds by {\cred Theorem \ref{t-wild-cb}(2)}. 
Hence we may assume that each contraction is generically smooth. 
Then the assertion follows from 
Lemma \ref{pic2-r2C-K^3}, Lemma \ref{deg1-isom}, and Lemma \ref{l-CC-deg2}. 
\end{proof}

\subsection{Case $C-D$}\label{ss-rho2-CD}

\begin{nota}\label{n-pic-2-D}
We use Notation \ref{n-pic-2}. 
Assume that $R_1$ is of type $C$ and $R_2$ is of type $D$. 
In particular, $f_1 : X \to Y_1 = \mathbb P^2$ is of type $C$
and 
$f_2 : X \to Y_2 =\mathbb P^1$ is of type $D$. 
If $f_2$ is of type $D_i$ for $i \in \{1, 2, 3\}$, then we have $\mu_2 = i$ (Remark \ref{r-length}). 
Set 
\[
f:=f_1 \times f_2 \colon X\to \mathbb{P}^2\times \mathbb{P}^1.
\]
\end{nota}

\begin{lem}\label{pic2-r2D-K^3}\label{l-pic2-D-1}
We use Notation \ref{n-pic-2-D}. 
Then $(-K_X)^3=6\mu_2^2$. 
\end{lem}

\begin{proof}
  By Proposition~\ref{basis}, $-K_X\sim \mu_2H_1+\mu_1H_2$. 
By $L_1 \simeq \MO_{\mathbb P^2}(1)$ and $L_2 \simeq \MO_{\mathbb P^1}(1)$, 
we have $H_1^3=f_1^*L_1^3=0$ and $H_2^2=f_2^*L_2^2\equiv 0$.
Therefore, we obtain 
  \begin{align*}
    (-K_X)^3 &= (\mu_2H_1+\mu_1H_2)^3 \\
    &= \mu_2^3 H_1^3 + 3\mu_2^2\mu_1 H_1^2\cdot H_2 + 3\mu_2\mu_1^2 H_1\cdot H_2^2 + \mu_2^3H_2^3 \\
    &= 3\mu_2^2\mu_1\cdot \left( \frac{2}{\mu_1} \ell_1 \right) \cdot H_2\\
    &= 6\mu_2^2,
  \end{align*}
where the third equality holds by 
$H_1^2 \equiv \frac{2}{\mu_1}\ell_1$ (Lemma~\ref{pic2-numequiv}) 
and the fourth one follows from 
$H_2 \cdot \ell_1 =1$ (Proposition~\ref{basis}). 
\end{proof}

\begin{lem}\label{l-pic2-D-2}
We use Notation \ref{n-pic-2-D}. 
  \begin{enumerate}
    \item If $R_1$ is of type $C_1$, 
then $f : X \to \mathbb{P}^2\times \mathbb{P}^1$ is a double cover. 
    \item If $R_1$ is of type $C_2$, 
then $f : X \to \mathbb{P}^2\times \mathbb{P}^1$ is an isomorphism 
and $R_2$ is of type $D_3$. 
  \end{enumerate}
\end{lem}

\begin{proof}
By Lemma \ref{l-finite-morph}(1), $f$ is a finite morphism. 
It follows from $L_1 \simeq \MO_{\mathbb P^2}(1)$ and $L_2 \simeq \MO_{\mathbb P^1}(1)$ 
that $\deg f= (f_1^*L_1)^2 \cdot f_2^*L_2$. 
By Corollary~\ref{fib_of_C}, {\cred every}  fibre of $f_1 : X \to Y_1 = \mathbb P^2$ is numerically equivalent to $\frac{2}{\mu_1}\ell_1$, i.e., $(f_1^*L_1)^2 \equiv \frac{2}{\mu_1}\ell_1$. 
We then obtain  
  \[\deg f = (f_1^*L_1)^2 \cdot f_2^*L_2 
= \frac{2}{\mu_1}\ell_1 \cdot H_2=\frac{2}{\mu_1},\]
where the last equality holds by Proposition \ref{basis}. 

(1) 
Assume that $R_1$ is of type $C_1$. 
We then have $\mu_1=1$, and hence $f: X \to \mathbb{P}^2\times \mathbb{P}^1$ is a double cover. 
%

(2) 
Assume that ${\cred R_1}$ is of type $C_2$. 
Then we have $\mu_1={\cred 2}$ and $\deg f=1$. 
Hence, $f: X \to \mathbb{P}^2\times \mathbb{P}^1$ is a finite birational morphism of normal varieties, which is automatically an isomorphism. 
In particular, $R_2$ is of type $D_3$. 
\end{proof}

\begin{lem}\label{l-pic2-D-3}
We use Notation \ref{n-pic-2-D}. 
Assume that $R_1$ is of type $C_1$. 
Set $\mathcal L := (f_*\MO_X/\MO_{\mathbb P^2 \times \mathbb P^1})^{-1}$. 
Then the following hold. 
\begin{enumerate}
\item $R_2$ is of type $D_1$ or $D_2$. 
\item $\mathcal L \simeq \MO_{\mathbb P^2 \times \mathbb P^1}(3-\mu_2, 1)$. 
\item 
$f: X \to \mathbb P^2 \times \mathbb P^1$ is a split double cover. 
\end{enumerate}
\end{lem}

\begin{proof}
By $\mu_1=1$, we have $K_X \sim -\mu_2 H_1 -H_2$, i.e., 
$\omega_X \simeq f^*\MO_{\P^2 \times \P^1}(-\mu_2, -1)$. 
Recall that we have $\omega_X \simeq f^*(\omega_{\mathbb P^2 \times \mathbb P^1} \otimes \mathcal L)$ (Lemma \ref{l-dc-omega}). 
It holds that $\mathcal L \simeq \MO_{\mathbb{P}^2\times \mathbb{P}^1}(b_1, b_2)$ for some $b_1, b_2 \in \Z$. 
Then 
\begin{align*}
\omega_X 
  &\simeq f^*\left(\omega_{\mathbb{P}^2\times \mathbb{P}^1} \otimes \mathcal L\right) \\
  &\simeq f^*\left( \MO_{\mathbb{P}^2\times \mathbb{P}^1}(-3, -2) \otimes \MO_{\mathbb{P}^2\times \mathbb{P}^1}(b_1, b_2)\right) \\
  &\simeq  f^* \MO_{\mathbb{P}^2\times \mathbb{P}^1}(b_1-3, b_2-2). 
\end{align*}
Hence $b_1 = 3- \mu_2$ and $b_2 = 1$. 
Thus (2) holds. 

Let us show (3), i.e.,  
\[
0 \to \MO_{\mathbb P^2 \times \mathbb P^1} \to f_*\MO_X \to \mathcal L \to 0 
\]
splits. 
The extension class lies in 
\[
{\rm Ext}^1_{\mathbb P^2 \times \mathbb P^1}(\mathcal L, \MO_{\mathbb P^2 \times \mathbb P^1}) 
\simeq H^1(\mathbb P^2 \times \mathbb P^1, \mathcal L^{-1})\simeq H^1(\mathbb P^2 \times \mathbb P^1, \MO_{\mathbb P^2 \times \mathbb P^1}(\mu_2-3, -1))=0, 
\]
where the last equality holds by the K\"{u}nneth formula. 
Thus (3) holds. 

Let us show (1). 
Suppose that $R_2$ is of type $D_3$. 
In this case, we have that $X \simeq S \times \mathbb P^1$ for some  smooth projective surface $S$ (Lemma \ref{l-trivial-double}). 
Then $X$ has an extremal ray of type $C_2$, 
which is a contraction. 
Thus (1) holds. 
\end{proof}


\begin{thm}\label{t-pic2-D}
Let $X$ be a  Fano threefold with $\rho(X)=2$. 
Let $R_1$ and $R_2$ be {\cred the} extremal rays of $\NE(X)$. 
Assume that $R_1$ is of type $C$ and $R_2$ is of type $D$. 
Then one and only one of the following holds.  
\begin{enumerate}
\item 
$R_1$ is of type $C_1$, $R_2$ is of type $D_1$, $(-K_X)^3 =6$, 
and 
there exists a split double cover $f: X \to \P^2 \times \P^1$ such that 
$(f_*\MO_X/\MO_{\P^2 \times \P^1})^{-1} \simeq \MO_{\P^2 \times \P^1}(2, 1)$ 
{\cred (No.\ \hyperref[table-2-2]{2-2})}. 
\item 
$R_1$ is of type $C_1$, $R_2$ is of type $D_2$, $(-K_X)^3 =24$, 
and 
there exists a split double cover $f: X \to \P^2 \times \P^1$ such that 
$(f_*\MO_X/\MO_{\P^2 \times \P^1})^{-1}  \simeq \MO_{\P^2 \times \P^1}(1, 1)$ {\cred (No.\ \hyperref[table-2-18]{2-18})}. 
\item 
$R_1$ is of type $C_2$, $R_2$ is of type $D_3$, $(-K_X)^3 =54$, 
and $X \simeq \P^2 \times \P^1$ {\cred (No.\ \hyperref[table-2-34]{2-34})}. 
\end{enumerate}
\end{thm}

\begin{proof}
The assertions follow from 
Lemma \ref{l-pic2-D-1}, 
Lemma \ref{l-pic2-D-2}, and 
Lemma \ref{l-pic2-D-3}. 
\end{proof}

\subsection{Case $C-E$ (primitive)}\label{ss-rho2-CE}





\begin{lem} \label{r1-c*}
We use Notation \ref{n-pic-2}. 
Assume that $R_1$ is of type $C$ and $R_2$ is of type $E_2, E_3, E_4$, or $E_5$. 
Set $D := \Ex(f_2 : X \to Y_2)$, which is a prime divisor on $X$. 
Then the following hold. 
  \begin{enumerate} 
    \item If $R_1$ is of type $C_1$, then 
$f_1|_D : D \to Y_1 =\mathbb P^2$ is a double cover and 
$R_2$ is of type $E_3$ or $E_4$.
    \item If $R_1$ is of type $C_2$, then 
$f_1|_D : D \to Y_1 =\mathbb P^2$ is an isomorphism and 
$R_2$ is of type $E_2$ or $E_5$.
  \end{enumerate}
\end{lem}

\begin{proof}
By Lemma \ref{l-finite-morph}(2), $f_1|_D : D \to Y_1 =\mathbb P^2$ is a finite surjective morphism. 
\setcounter{step}{0}

\begin{step}\label{s1-r1-c*}
It holds that 
\[
\deg(f_1|_D)=\frac{2}{\mu_1} D\cdot\ell_1=\frac{1}{\mu_2^2}(\omega_D \otimes \MO_D(-D))^2, 
\]
where $\MO_D(-D) := \MO_X(-D)|_D$. 
\end{step}

\begin{proof}[Proof of Step \ref{s1-r1-c*}]
By Corollary~\ref{fib_of_C}, a fibre  of $f_1$ is numerically equivalent to $(2/\mu_1)\ell_1$. 
We then obtain  
  \[
\deg(f_1|_D)=\frac{2}{\mu_1} D\cdot\ell_1 = D \cdot H_1^2, 
\]
where the latter equality follows from Lemma~\ref{pic2-numequiv}. 
  By Proposition~\ref{basis}(3), we have
  \begin{align*}
    H_1^2 &\equiv \frac{1}{\mu_2^2}(-K_X-\mu_1H_2)^2 \\
    &= \frac{1}{\mu_2^2}((-K_X)^2-2\mu_1(-K_X)\cdot H_2+\mu_1^2 H_2^2).
  \end{align*}
  Since $R_2$ is of type $E_2, E_3, E_4$ or $E_5$, $f_2(D)$ is a point on ${\cred Y_2}$. 
Therefore, we get 
  \[ H_2\cdot D = f^*_2L_2\cdot D \equiv 0.\]
Hence we obtain 
\[
\deg(f_1|_D) = D \cdot H_1^2
=
\frac{1}{\mu_2^2}D \cdot ((-K_X)^2-2\mu_1(-K_X)\cdot H_2+\mu_1^2 H_2^2)
\]
\[
= \frac{1}{\mu_2^2}D \cdot (-K_X)^2 = \frac{1}{\mu_2^2}(\omega_D \otimes \MO_D(-D))^2,
\]
where the last equality holds by the adjunction formula. 
This completes the proof of Step \ref{s1-r1-c*}. 
\end{proof}

\begin{step}\label{s2-r1-c*}
The assertion (1) holds. 
\end{step}

\begin{proof}[Proof of Step \ref{s2-r1-c*}]
Assume that $R_1$ is of type $C_1$. 
We have $\mu_1=1$.  
Recall that $R_2$ is of type $E_2, E_3, E_4,$ or $E_5$.

Assume that $R_2$ is of type $E_2$. 
Then we have $\mu_2=2, D\simeq \mathbb{P}^2$, and $\MO_D(D)\simeq \MO_D(-1)$.
By Step \ref{s1-r1-c*},  we get 
  \[
  2 D\cdot\ell_1 =
  \frac{1}{\mu_2^2}(\omega_D \otimes \MO_D(-D))^2 = 
  \frac{1}{4}(-3+1)^2=1.
  \]
This is absurd. 


Assume that $R_2$ is of type $E_5$. 
Then $\mu_2=1, D\simeq \mathbb{P}^2,$ and $\MO_D(D)\simeq \MO_{\mathbb{P}^2}(-2)$. 
By Step \ref{s1-r1-c*},  we get 
  \[
  2D\cdot\ell_1  =
  \frac{1}{\mu_2^2}(\omega_D \otimes \MO_D(-D))^2 = (\MO_{\mathbb{P}^2}(-3+2))^2=1.\]
This is a contradiction.

Assume that $R_2$ is of type $E_3$ or $E_4$. 
Then we have $\mu_2=1$, 
$D$ is isomorphic to a possibly singular quadric surface in $\mathbb P^3$, 
and $\MO_D(D)\simeq \MO_{\mathbb{P}^3}(-1)|_D$. 
We have 
\[
\omega_D \simeq 
(\omega_{\mathbb{P}^3}\otimes \MO_{\mathbb{P}^3}(D))|_D 
\simeq \MO_{\mathbb{P}^3}(-4+2)|_D \\
\simeq \MO_{\mathbb{P}^3}(-2)|_D. 
\]
By Step \ref{s1-r1-c*},  we get 
\[
\deg(f_1|_D)=    2 D\cdot \ell_1 
=
  \frac{1}{\mu_2^2}(\omega_D \otimes \MO_D(-D))^2 
= (\MO_{\mathbb{P}^3}(-2+1)|_D)^2 =2. 
\]
This completes the proof of Step \ref{s2-r1-c*}. 
\end{proof}

\begin{step}\label{s3-r1-c*}
The assertion (2) holds. 
\end{step}

\begin{proof}[Proof of Step \ref{s3-r1-c*}]
Assume that $R_1$ is of type $C_2$. 
We have $\mu_1=2$. 

Suppose that $R_2$ is of type $E_3$ or $E_4$. Then $\mu_2=1$, and 
hence  Proposition \ref{basis} implies 
  \[
  -K_X \sim H_1+2H_2.
  \]
We then get
  \begin{align*}
    24&=c_2(X)\cdot (-K_X) \\
    &=c_2(X)\cdot H_1 + 2c_2(X)\cdot H_2  \\
    &=6+\frac{48}{r}
  \end{align*}
  for some $r\in \mathbb{Z}_{>0}$ (Lemma \ref{c_times_H}).
  This is impossible. 

Hence $R_2$ is of type $E_2$ or $E_5$. 
It holds that 
\[
\deg(f_1|_D)= D\cdot\ell_1 
=\frac{1}{\mu_2^2}(\omega_D \otimes \MO_D(-D))^2 =1, 
\]
where the first two equalities are guaranteed by Step \ref{s1-r1-c*} 
and the last one follows from the same argument as in 
Step \ref{s2-r1-c*}. 
  Then $f_1|_D : D \to \mathbb P^2$ is a finite birational morphism 
of normal varieties, which is automatically an isomorphism.
This completes the proof of Step \ref{s3-r1-c*}. 
\end{proof}
Step \ref{s2-r1-c*} and Step \ref{s3-r1-c*} 
complete the proof of Lemma \ref{r1-c*}. 
\end{proof}


\subsubsection{Case $C_1-E$  {\cred (primitive)}}

\begin{nota}\label{n-pic-2-C_1E}
We use Notation \ref{n-pic-2}. 
Assume that $R_1$ is of type $C_1$ and $R_2$ is of type $E_2, E_3, E_4,$ or $E_5$. 
Set $D := \Ex(f_2 : X \to Y_2)$. 
In particular, $f_1 : X \to Y_1 = \mathbb P^2$ is of type $C_1$ and 
$f_2 : X \to Y_2$ is of type $E_3$ or $E_4$ (Lemma \ref{r1-c*}). 
Recall that $D$ is a (possibly singular) quadric surface in $\P^3$ and set $\MO_D(n) := \MO_{\P^3}(n)|_D$ for $n \in \Z$. 
By Lemma~\ref{conic-embedding}, we have a diagram:
\begin{center}
  \begin{tikzcd}
    X \arrow[r,hook, "\iota"] & P:=\mathbb{P}((f_1)_*\MO_X(-K_X)) \arrow[d,"\pi"]\\
    & \mathbb{P}^2.
  \end{tikzcd}
\end{center}
Since $\iota$ is a closed immersion, we identify $X$ with the smooth prime divisor $\iota(X)$ 
on the $\mathbb{P}^2$-bundle $P$ over $\mathbb{P}^2$.
\end{nota}

\begin{lem}\label{l-pic2E-ext-seq}
We use Notation \ref{n-pic-2-C_1E}. 
Then   $R^1(f_1)_*\MO_X(-K_X-D)=0$ and there exists an exact sequence: 
\[ 
0 \longrightarrow (f_1)_*\MO_X(-K_X-D) \longrightarrow (f_1)_*\MO_X(-K_X) \longrightarrow 
(f_1)_*\MO_D(-K_X) \longrightarrow 0. \]
\end{lem}

\begin{proof}
Since $(-K_X -D)-K_X$ is $f_1$-ample (Lemma \ref{r1-c*}), 
the assertion follows from 
$R^1(f_1)_*\MO_X(-K_X-D)=0$, 
which is guaranteed by \cite[Theorem 0.5]{Tan15}. 
\end{proof}



\begin{lem}\label{l-pic2E-D-pol}
We use Notation \ref{n-pic-2-C_1E}. 
Then  $(f_1|_D)^*\MO_{\mathbb{P}^2}(n)\simeq \MO_D(n)$ for all $n\in \mathbb{Z}$. 
\end{lem}

\begin{proof}
We may assume $n=1$.   

  Assume that  $R_2$ is of type $E_3$. 
  Then $D\simeq \mathbb{P}^1\times \mathbb{P}^1$ and we can write
  \[ (f_1|_D)^*\MO_{\mathbb{P}^2}(1)\simeq \MO_D(a,b) \]
  for some $a, b\in \mathbb{Z}$. Since $f_1|_D : D \to \P^2$ is a finite morphism, 
  this invertible sheaf is ample, which implies $a>0$ and $b>0$. 
  By the projection formula, we get 
  \begin{align*}
    2ab &= (\MO_D(a,b))^2 \\
    &= \deg(f_1|_D) \cdot (\MO_{\mathbb{P}^2}(1))^2 \\
    &=2.
  \end{align*}
  Hence we have $ab=1$. 
  We then get $(a, b) = (1, 1)$ and 
  \[
  (f_1|_D)^*\MO_{\mathbb{P}^2}(1)\simeq \MO_{\P^1 \times \P^1}(1, 1) \simeq \MO_{\P^3}(1)|_D = \MO_D(1).
  \]
  

  \medskip

  Assume that $R_2$ is of type $E_4$. 
  By \cite[Ch. II, Exercise 6.3 and 6.5]{Har77}, the restriction homomorphism $\Pic(\mathbb{P}^3) \to \Pic(D)$ is an isomorphism. 
  Hence we have 
  \[ (f_1|_D)^*\MO_{\mathbb{P}^2}(1)\simeq \MO_D(a) \]
  for some $a\in \mathbb{Z}$. 
  Since this invertible sheaf is ample, we get $a>0$. 
  By the projection formula,
  \[ 
  a^2 c_1(\MO_D(1))^2 =
  c_1(\MO_D(a))^2 = \deg(f_1|_D) \cdot c_1(\MO_{\mathbb{P}^2}(1))^2=2.
  \]
  Hence $a=1$. 
\end{proof}

\begin{prop}\label{p-pic2E-vb}
We use Notation \ref{n-pic-2-C_1E}. 
Then the following hold. 
\begin{enumerate}
\item $(f_1)_*\MO_X(-K_X-D)\simeq \MO_{\mathbb{P}^2}(2).$
\item $(f_1)_*\MO_D(-K_X) \simeq \MO_{\mathbb{P}^2}\oplus \MO_{\mathbb{P}^2}(1)$. 
\item 
$(f_1)_*\MO_X(-K_X) \simeq \MO_{\mathbb{P}^2}\oplus \MO_{\mathbb{P}^2}(1)\oplus \MO_{\mathbb{P}^2}(2)$. 
\end{enumerate}
\end{prop}

\begin{proof}
Let us show (1). 
For a closed point $t \in \P^2$ and the fibre $X_t$ of $f_1 : X \to Y_1 = \P^2$ over $t$, we have $-K_X \cdot X_t =D \cdot X_t =2$. 
Hence $h^0(X_t,\MO_{X_t}(-K_X-D))=1$ for all $t\in Y_1 = \mathbb{P}^2$.
By Grauert's theorem \cite[Ch. III, Theorem 12.9]{Har77}, 
$(f_1)_*\MO_X(-K_X-D)$ is an invertible sheaf on $\P^2$. 
Hence we can write 
\[
(f_1)_*\MO_X(-K_X-D)\simeq \MO_{\mathbb{P}^2}(a)\]
for some $a\in \mathbb{Z}$.
Again by \cite[Ch. III, Theorem 12.9]{Har77}, we get 
\[\MO_X(-K_X-D)\simeq f^*_1(f_1)_*\MO_X(-K_X-D) \simeq f^*_1 \MO_{\mathbb{P}^2}(a).\]
By Lemma \ref{l-pic2E-D-pol} and adjunction formula, we obtain
\[
\MO_D(a)\simeq (f_1|_D)^*\MO_{\mathbb{P}^2}(a) \simeq \MO_D(-K_X-D) \simeq \omega_D^{-1} \simeq \MO_{\P^3}(2)|_D = \MO_D(2), 
\]
which implies $a=2$. Thus (1) holds. 


Let us show (2).  
By $\MO_D(-K_X-D)\simeq \MO_D(2)$ and $\MO_D(D)\simeq \MO_D(-1)$, 
we obtain $\MO_D(-K_X) \simeq \MO_D(1)$. 
Then 
\begin{align*}
  (f_1|_D)_*\MO_D(-K_X) &\simeq (f_1|_D)_*\MO_D(1) \\
  &\simeq (f_1|_D)_*(\MO_D \otimes (f_1|_D)^*\MO_{\mathbb{P}^2}(1)) \\
  &\simeq (f_1|_D)_*\MO_D \otimes \MO_{\mathbb{P}^2}(1).
\end{align*}
Since $f_1|_D \colon D\to \mathbb{P}^2$ is a double cover (Lemma \ref{r1-c*}), 
we have an exact sequence 
\[
0 \to \MO_{\mathbb P^2} \to (f_1|_D)_*\MO_D \to \MO_{\mathbb P^2}(-b)\to 0
\]
for some $b \in \Z$. 
Since this sequence splits, we get $(f_1|_D)_*\MO_D \simeq \MO_{\mathbb{P}^2}\oplus \MO_{\mathbb{P}^2}(-b)$. 
Furthermore, $\omega_D \simeq (f_1|_D)^*(\omega_{\mathbb{P}^2}\otimes \MO_{\mathbb{P}^2}(b))$ (Lemma \ref{l-dc-omega}). 
By $\omega_D \simeq \MO_D(-2)$, we obtain 
\begin{align*}
  \MO_D(-2) &\simeq \omega_D \\
  &\simeq (f_1|_D)^*(\MO_{\mathbb{P}^2}(-3)\otimes \MO_{\mathbb{P}^2}(b)) \\
  &\simeq (f_1|_D)^*(\MO_{\P^2}(b-3)) \\
  &\simeq \MO_D(b-3), 
\end{align*}
where the last isomorphism holds by Lemma \ref{l-pic2E-D-pol}. 
Thus $b=1$ and (2) holds.  

Let us show (3). 
By (1), (2), and Lemma \ref{l-pic2E-ext-seq}, we have the following exact sequence: 
\[0 \longrightarrow \MO_{\mathbb{P}^2}(2) \longrightarrow (f_1)_*\MO_X(-K_X) \longrightarrow \MO_{\mathbb{P}^2}\oplus \MO_{\mathbb{P}^2}(1) \longrightarrow 0.\]
By $\mathrm{Ext}^1(\MO_{\mathbb{P}^2}\oplus \MO_{\mathbb{P}^2}(1),\MO_{\mathbb{P}^2}(2)) =0$, 
this exact sequence splits, and hence 
\[
(f_1)_*\MO_X(-K_X) \simeq \MO_{\mathbb{P}^2}\oplus \MO_{\mathbb{P}^2}(1)\oplus \MO_{\mathbb{P}^2}(2). 
\]
Thus (3) holds. 
\end{proof}

\begin{lem}\label{l-pic2E-XinP}
We use Notation \ref{n-pic-2-C_1E}. 
Then 
{\cred $\MO_X(-K_X) \simeq \MO_P(1)|_X$ and} 
$\MO_P(X) \simeq \MO_P(2)$, 
where $\MO_P(1)$ denotes the tautological line bundle 
of the $\P^2$-bundle structure $ \pi : P =\P( (f_1)_*\MO_X(-K_X)) \to \P^2$ and we set $\MO_P(2) := \MO_P(1)^{\otimes 2}$. 
\end{lem}

\begin{proof}
Recall that we have the following diagram (Notation \ref{n-pic-2-C_1E}, Proposition \ref{p-pic2E-vb}): 
\begin{center}
  \begin{tikzcd}
    X \arrow[r, hook, "\iota"] & P=\mathbb{P}(\MO_{\mathbb{P}^2}\oplus \MO_{\mathbb{P}^2}(1)\oplus \MO_{\mathbb{P}^2}(2) ) \arrow[d,"\pi"]\\
    & \mathbb{P}^2.
  \end{tikzcd}
\end{center}
By $\Pic\,P\simeq \pi^*\Pic\,\mathbb{P}^2 \oplus \mathbb{Z}\MO_P(1)$, we can write 
\begin{equation}\label{rho2-a-n}
  \MO_P(X)\simeq \pi^*\MO_{\mathbb{P}^2}(a)\otimes \MO_P(n)
\end{equation}
for some $a,n\in \mathbb{Z}$.
It follows from the adjunction formula that  
\begin{equation}\label{rho2-P-X-adjunction}
  \MO_X(K_X)\simeq \MO_P(K_P+X)|_X.
\end{equation}
By Proposition~\ref{standard_2}(2), it holds that 
\begin{equation}\label{KP=-3}
\MO_P(K_P)\simeq \MO_P(-3)\otimes \pi^*(\omega_{\mathbb{P}^2}\otimes \det(\MO_{\mathbb{P}^2}\oplus \MO_{\mathbb{P}^2}(1) \oplus \MO_{\mathbb{P}^2}(2))) \simeq \MO_P(-3).
\end{equation}
Since the closed immersion $\iota \colon X \to P$ is induced by the surjection $(f_1)^*(f_1)_*\MO_X(-K_X) \to  \MO_X(-K_X)$,  
it follows from \cite[Ch. II, the proof of Proposition 7.12]{Har77} that 
\begin{equation}\label{rho2-X-K_X}
  \MO_X(-K_X) \simeq \iota^*\MO_P(1) \simeq \MO_P(1)|_X. 
\end{equation}
By (\ref{rho2-a-n}), (\ref{rho2-P-X-adjunction}), (\ref{KP=-3}) 
and (\ref{rho2-X-K_X}), we have
\[
\MO_P(1)|_X \simeq \MO_X(-K_X)  \simeq \MO_P(-K_P-X)|_X \simeq (\pi^*\MO_{\mathbb{P}^2}(-a)\oplus \MO_P(3-n))|_X.
\]
Taking the intersection number of this equation with a general fibre of $\pi|_X \colon X \to \mathbb{P}^2$, 
we have $1=3-n$, i.e., $n=2$. 
Hence we obtain $\MO_X\simeq \pi^*\MO_{\mathbb{P}^2}({\cred -}a)|_X$. 
We see that $a=0$ by taking the intersection number of this equation with a curve $C$ on $X$  such that $\pi(C)$ is a curve. 
\end{proof}

\begin{lem}\label{l-pic2E-vol}
We use Notation \ref{n-pic-2-C_1E}. 
Then $(-K_X)^3 =14$. 
\end{lem}

\begin{proof}
{\cred Set $\mathcal E := \MO_{\P^2} \oplus \MO_{\P^2}(1) \oplus \MO_{\P^2}(2)$.} 
It follows from Proposition~\ref{standard_2}(5) that 
\begin{align*}
  (-K_X)^3 = &2c_1(\mathcal{E})^2-2c_2(\mathcal{E})+4c_1(\mathcal{E})\cdot c_1(\MO_{\mathbb{P}^2}) +6 c_1(\mathcal{E})\cdot K_{\mathbb{P}^2} \\
  & + 9c_1(\MO_{\mathbb{P}^2})\cdot K_{\mathbb{P}^2} + 6K_{\mathbb{P}^2}^2+3c_1(\MO_{\mathbb{P}^2})^2.
\end{align*}
Moreover, we have  
\begin{align*}
  &c_1(\mathcal{E})=c_1(\MO_{\mathbb{P}^2})+c_1(\MO_{\mathbb{P}^2}(1))+c_1(\MO_{\mathbb{P}^2}(2))=c_1(\MO_{\mathbb{P}^2}(3)) \\
  &c_2(\mathcal{E})=c_1(\MO_{\mathbb{P}^2}(1))\cdot c_1(\MO_{\mathbb{P}^2}(2))=2 \\
  &c_1(\MO_{\mathbb{P}^2})=0
\end{align*}
by \cite[Appendix A, \S3, C3 and C.5]{Har77}.
Hence we obtain
\[
  (-K_X)^3 = 2\cdot 3^2-2\cdot 2+  4\cdot 0+6\cdot 3 \cdot(-3) + 9\cdot 0+6\cdot(-3)^2 + 3\cdot 0 \\
=18-4  -54+54
=14.
\]
\end{proof}


\begin{prop}\label{p-pic2-C1E-final} 
We use Notation \ref{n-pic-2-C_1E}. 
Then there exist a split double cover 
\[
f: X \to V_7
\]
such that $K_X \sim f^*( \frac{1}{2} K_{V_7})$, where $V_7$ is the blowup of $\P^3$ at a point and 
$\frac{1}{2} K_{V_7}$ is a Weil divisor on $V_7$ satisfying $2 \cdot (\frac{1}{2} K_{V_7}) \sim K_{V_7}$ (note that $\frac{1}{2} K_{V_7}$ is unique up to linear equivalence). 
In particular, for $\mathcal L := (f_*\MO_X/\MO_{V_7})^{-1}$, we have $\mathcal L \simeq \MO_{V_7}(-\frac{1}{2} K_{V_7})$. 
\end{prop}

\begin{proof}
Set $S \subset P$ to be the section of $\pi : P = \P_{\P^2}(\MO_{\P^2} \oplus \MO_{\P^2}(1) \oplus \MO_{\P^2}(2)) \to Z := \P^2$ corresponding 
to the natural projection: 
\[
\MO_{\P^2} \oplus \MO_{\P^2}(1) \oplus \MO_{\P^2}(2) \to \MO_{\P^2}. 
\]
Let 
$\mu : \widetilde P \to P$ be the blowup along $S$.

\setcounter{step}{0}

\begin{step}\label{s1-pic2-C1E-final}
There exists the following commutative diagram 
\begin{equation}\label{e1-pic2-C1E-final}
\begin{tikzcd}
	& \widetilde P \\
	\P(\MO_{\P^2} \oplus \MO_{\P^2}(1) \oplus \MO_{\P^2}(2))=P && 
\P(\MO_{\P^2}(1) \oplus \MO_{\P^2}(2))=V_7 \\
	& \P^2=Z
	\arrow["\mu"', from=1-2, to=2-1]
	\arrow["\pi"', from=2-1, to=3-2]
	\arrow["\alpha", from=1-2, to=2-3]
	\arrow["\beta", from=2-3, to=3-2]
\end{tikzcd}
\end{equation}
such that for {\cred every}  closed point $z \in \P^2 =Z$, the base change of (\ref{e1-pic2-C1E-final}) by $z \hookrightarrow \P^2=Z$ can be written as follows
\begin{equation}\label{e2-pic2-C1E-final}
\begin{tikzcd}
	& \F_1:=\P(\MO_{\P^1} \oplus \MO_{\P^1}(1)) =\widetilde P_z \\
	\P^2=P_z && 
\P^1 =(V_7)_z\\
	& z
	\arrow["\mu_z"', from=1-2, to=2-1]
	\arrow["\pi_z"', from=2-1, to=3-2]
	\arrow["\alpha_z", from=1-2, to=2-3]
	\arrow["\beta_z", from=2-3, to=3-2]
\end{tikzcd}
\end{equation}
where $\mu_z$ is the blowup at the closed point $S_z$ and $\alpha_z : \F_1 \to \P^1$ is the $\P^1$-bundle. 

\end{step}

\begin{proof}[Proof of Step \ref{s1-pic2-C1E-final}]
We have the canonical injective graded $\MO_{\P^2}$-algebra homomorhism
\[
\iota : {\rm Sym} (\MO_{\P^2}(1) \oplus \MO_{\P^2}(2)) \hookrightarrow {\rm Sym} (\MO_{\P^2} \oplus \MO_{\P^2}(1) \oplus \MO_{\P^2}(2)). 
\]
Then $\iota$ induces a dominant rational map $\varphi$ over $\P^2$ as follows: 
\[
\begin{tikzcd}
	P= \P_{\P^2}(\MO_{\P^2} \oplus \MO_{\P^2}(1) \oplus \MO_{\P^2}(2))  && 
    \P_{\P^2}(\MO_{\P^2}(1) \oplus \MO_{\P^2}(2))=V_7 \\
	& Z = \P^2
	\arrow["\pi"', from=1-1, to=2-2]
	\arrow["\beta", from=1-3, to=2-2]
	\arrow["\varphi", dashed, from=1-1, to=1-3]
\end{tikzcd}
\]
Note that this construction commutes with 
base changes by open immersions   $Z' \hookrightarrow Z = \P^2$. 
Then the indeterminacies of $\varphi$ are resolved 
by the blowup along $S$, 
because this can be checked after taking an affine open cover of $Z =\P^2$ trivialising $\MO_{\P^2}(n)$. 
Therefore, we get a commutative diagram (\ref{e1-pic2-C1E-final}). 
This diagram fibrewisely induces the diagram (\ref{e2-pic2-C1E-final}), because 
we have the corresponding diagram over an arbitrary open subset of $Z= \P^2$ which trivialises $\MO_{\P^2}(n)$. 
This completes the proof of Step \ref{s1-pic2-C1E-final}. 
\end{proof}

\begin{step}\label{s2-pic2-C1E-final}
It holds that $X \cap S = \emptyset$. 
In particular, $\mu^{-1}(X) \simeq X$. 
\end{step}

\begin{proof}[Proof of Step \ref{s2-pic2-C1E-final}]
Suppose $X \cap S \neq  \emptyset$. 
{\cred 
Since $P = \P_{\P^2}(\MO_{\P^2} \oplus \MO_{\P^2}(1) \oplus \MO_{\P^2}(2))$ is smooth, 
we have ${\rm codim}_P\,T \leq {\rm codim}_P\,X + {\rm codim}_P\,S =1+2=3$ for an irreducible component $T$ of $X \cap S$ \cite[Theorem 0.2]{EH16}.} 
Hence there exists a curve $C$ with $C \subset X\cap S$. 
By $\MO_X(-K_X)\simeq \MO_P(1)|_X$ {\cred (Lemma \ref{l-pic2E-XinP})} and $\MO_P(1)|_S \simeq \MO_S$, we obtain 
\[
-K_X \cdot C = (\MO_P(1)|_X)\cdot C =\MO_P(1) \cdot C = (\MO_P(1)|_S) \cdot C = \MO_S \cdot C=0. 
\]
This contradicts the ampleness of $-K_X$. 
This completes the proof of Step \ref{s2-pic2-C1E-final}. 
\end{proof}

\begin{step}\label{s3-pic2-C1E-final}
The induced morphism 
\[
f : X \xrightarrow{\simeq} \mu^{-1}(X) \xrightarrow{\alpha|_{\mu^{-1}(X)}} V_7
\]
is a double cover. 
\end{step}

\begin{proof}[Proof of Step \ref{s3-pic2-C1E-final}]
Fix a closed point $v \in V_7$. 
Set $z := \beta(v) \in Z =\P^2$. 
Take the base change by $z \hookrightarrow \P^2$. 
We then obtain the above diagram (\ref{e2-pic2-C1E-final}) and the fibre $X_z$ is a conic in $\P^2  =P_z$. 
It holds by Step \ref{s2-pic2-C1E-final} that  
$X_z$ is disjoint from the $(-1)$-curve on 
the blowup $\widetilde P_z = \F_1$ of $P_z = \P^2$ lying over the point $S_z$. 
We then have $X_z \cdot \zeta = 2$ for a fibre $\zeta$ of $\alpha_z : \widetilde P_z \to (V_7)_z$, 
because $\zeta$ is nothing but the proper transform of a line on $\P^2 = P_z$ passing through $S_z$. 
Since $f^{-1}(v)$ consists of at most two points, $f : X \to V_7$ is a finite surjective morphism. 
By $X_z \cdot \zeta =2$, $f$ is a double cover. 
This completes the proof of Step \ref{s3-pic2-C1E-final}. 
\end{proof}

\begin{step}\label{s4-pic2-C1E-final}
It holds that $K_X \sim f^*( \frac{1}{2} K_{V_7} + \beta^*M)$
 for some Cartier divisor $M$ on $\P^2$. 
\end{step}

\begin{proof}[Proof of Step \ref{s4-pic2-C1E-final}]
The above diagram (\ref{e1-pic2-C1E-final}) consists of smooth projective toric varieties and toric morphisms. 
Furthermore, each morphism is of relative Picard number one. 

Set 
\[
N := -2(K_{\widetilde P} +\mu^*X) + \alpha^*K_{V_7}. 
\]
Fix a closed point $z \in Z  =\P^2$. 
We then have the diagram (\ref{e2-pic2-C1E-final}). 
Take a fibre $\zeta \simeq \P^1$ of 
the $\P^1$-bundle $\alpha_z : \F_1 =\widetilde P_z \to \P^1 =(V_7)_z$ over a closed point. 
Note that $\mu_z^*X_z \cdot \zeta =2$ (cf. the proof of Step \ref{s3-pic2-C1E-final}), and hence we get 
\begin{equation}\label{e6-pic2-C1E-final}
N \cdot \zeta = (-2(K_{\widetilde P} +\mu^*X) + \alpha^*K_{V_7}) \cdot \zeta = -2(-2+2) +0 =0 
\end{equation}
and 
\begin{equation}\label{e7-pic2-C1E-final}
N \cdot \mu^{-1}(X_z) = (-2(K_{\widetilde P} +\mu^*X) + \alpha^*K_{V_7}) \cdot \mu^{-1}(X_z) 
= -2 \deg \omega_{X_z} + 2 \deg K_{(V_7)_z} = 4-4=0. 
\end{equation}
By (\ref{e6-pic2-C1E-final}), there exists a Cartier divisor $N'$ on $V_7$ such that $N \sim \alpha^*N'$. 
It follows from (\ref{e7-pic2-C1E-final}) that we can find a Cartier divisor $N''$ on $\P^2 =Z$ such that $N' \sim \beta^*N''$. 
For $N''':= -N''$, 
we get $2(K_{\widetilde P} +\mu^*X) \sim \alpha^*(K_{V_7} + \beta^*N''')$. 
Taking the restriction to $\mu^{-1}(X) (\simeq X)$, we have $2K_X \sim f^*(K_{V_7} + \beta^*N''')$. 
Then 
$K_X - f^*(\frac{1}{2}K_{V_7})$ is numerically trivial over $Z$. 
Therefore, we can find a Cartier divisor $M$ on $Z$ such that $K_X =  f^*(\frac{1}{2}K_{V_7} + \beta^*M)$. 
This completes the proof of Step \ref{s4-pic2-C1E-final}. 
\end{proof}

\begin{step}\label{s5-pic2-C1E-final}
It holds that $K_X \sim f^*( \frac{1}{2} K_{V_7})$. 
\end{step}

\begin{proof}[Proof of Step \ref{s5-pic2-C1E-final}]
By Step \ref{s4-pic2-C1E-final}, we have $K_X \sim f^*(\frac{1}{2}K_{V_7} + \beta^*M)$.

We now show that $\deg M \leq 0$, i.e., $-M$ is nef.  
For the exceptional divisor $E \simeq \P^2$ of the blowup $V_7 \to \P^3$, it holds that 
\[
\MO_{V_7}\left( \frac{1}{2}K_{V_7}\right)\Big|_E \simeq \MO_E(-1). 
\]
Pick a line $\ell \subset E$. 
Since $-(\frac{1}{2}K_{V_7} + \beta^*M)$ is ample, we obtain 
\[
0<-\left( \frac{1}{2}K_{V_7} + \beta^*M \right) \cdot \ell = 1-\beta^*M \cdot \ell. 
\]
Since $E \to Z = \P^2$ is a finite surjective morphism, we get $\deg M \leq 0$.

By $K_X \sim f^*(\frac{1}{2}K_{V_7} + \beta^*M)$, we have 
\[
8K_X^3 = ( f^*(K_{V_7} + 2\beta^*M))^3 =  
2(K_{V_7}^3 + 6 K_{V_7}^2 \cdot \beta^*M + 12 K_{V_7} \cdot (\beta^*M)^2), 
\]
where the last equality holds by $(\beta^*M)^3=0$. 
It follows from $K_X^3 = -14$ (Lemma \ref{l-pic2E-vol}) and $K_{V_7}^3 = K_{\P^3}^3 +8 = -56$ that 
\[
K_{V_7}^2 \cdot \beta^*M + 2 K_{V_7} \cdot (\beta^*M)^2 = 0. 
\]
Suppose that $M \not\sim 0$. 
By $\deg M \leq 0$, we have $\deg M <0$, i.e., $-M$ is ample. 
We then have 
\[
K_{V_7}^2 \cdot \beta^*M = -(-K_{V_7})^2 \cdot (-\beta^*M)<0 \quad {\rm and} 
\quad K_{V_7} \cdot (\beta^*M)^2 = - (-K_{V_7}) \cdot (-\beta^*M)^2<0, 
\]
which is a contraction. 
Therefore, we obtain $M \sim 0$. 
This completes the proof of Step \ref{s5-pic2-C1E-final}. 
\end{proof}

\begin{step}\label{s6-pic2-C1E-final}
The induced morphism 
$f : X \to V_7$ is a split double cover. 
\end{step}

\begin{proof}[Proof of Step \ref{s6-pic2-C1E-final}]
For a Cartier divisor $L$ satisying 
\[
0 \to \MO_{V_7} \to f_*\MO_X \to \MO_{V_7}(-L) \to 0, 
\]
it suffices to show that $\Ext^1(\MO_{V_7}(-L), \MO_{V_7})=0$, 
which is equivalent to $H^1(V_7, \MO_{V_7}(L))=0$. 
We have $K_X \sim f^*(K_{V_7} + L)$ (Lemma \ref{l-dc-omega}). 
By Step \ref{s5-pic2-C1E-final}, we get $L \sim -\frac{1}{2}K_{V_7}$. 
Since $V_7$ is toric and $L$ is ample, 
we get $H^1(V_7, \MO_{V_7}(L))=0$ \cite[Theorem 1.6]{Fuj07}. 
This completes the proof of Step \ref{s6-pic2-C1E-final}. 
\end{proof}
Step \ref{s5-pic2-C1E-final} and Step \ref{s6-pic2-C1E-final} complete 
the proof of Proposition \ref{p-pic2-C1E-final}. 
\end{proof}

\subsubsection{Case $C_2-E$ {\cred (primitive)}}

\begin{lem}\label{l-pic2-C2E}
We use Notation \ref{n-pic-2}. 
Assume that $R_1$ is of type $C_2$ and $R_2$ is of type $E_2, E_3, E_4$, or $E_5$. 
Then the following hold. 
\begin{enumerate}
\item $R_2$ is of type $E_2$ or $E_5$. 
\item If $R_2$ is of type $E_2$, then $X \simeq \mathbb P(\MO_{\mathbb P^2} \oplus \MO_{\mathbb P^2}(1))$ and $(-K_X)^3=56$. 
\item 
If $R_2$ is of type $E_5$, then $X \simeq \mathbb P(\MO_{\mathbb P^2} \oplus \MO_{\mathbb P^2}(2))$ and $(-K_X)^3=62$. 
\end{enumerate}
\end{lem}

\begin{proof}
Set $D := \Ex(f_2 : X \to Y_2)$. 
Since $f_1\colon X\to \mathbb{P}^2$ is a $\mathbb{P}^1$-bundle, we can write 
\[f_1\colon X = \mathbb{P}_{\P^2}(\mathcal{E})\to \mathbb{P}^2 \]
for some locally free sheaf $\mathcal{E}$ of rank 2 on $\mathbb{P}^2$.
By Lemma~\ref{r1-c*}, $f_1|_D \colon D \to \mathbb{P}^2$ is an isomorphism and $R_2$ is of type $E_2$ or $E_5$. 
In particular, (1) holds. 
We have 
$\MO_X(D) \simeq \MO_X(n) \otimes (f_1)^*\mathcal L$ for some $n \in \Z$ and $\mathcal L \in \Pic\,\P^2$. 
Since $D \cdot F=1$ for a fibre $F \simeq \P^1$ of $f_1\colon X = \mathbb{P}_{\P^2}(\mathcal{E})\to \mathbb{P}^2$, 
we obtain $n=1$, i.e., 
\begin{equation}\label{e1-pic2-C2E}
\MO_X(D) \simeq \MO_X(1) \otimes (f_1)^*\mathcal L.
\end{equation}

We have the following exact sequence
\[0\longrightarrow \MO_X \longrightarrow \MO_X(D) \longrightarrow j_*(\MO_X(D)|_D) \longrightarrow 0,\]
where $j : D \hookrightarrow X$ denotes the induced closed immersion. 
By $R^1(f_1)_*\MO_X=0$ (Proposition \ref{p-cont-ex2}),  we obtain the exact sequence
\begin{equation}\label{e2-pic2-C2E}
0\longrightarrow (f_1)_*\MO_X \longrightarrow  (f_1)_*\MO_X(D) \longrightarrow (f_1)_*j_*(\MO_X(D)|_D) \longrightarrow 0.
\end{equation}
We have $(f_1)_*\MO_X = \MO_{\mathbb P^2}$. 
By (\ref{e1-pic2-C2E}) and $(f_1)_*\MO_X(1) \simeq \mathcal E$ \cite[Ch. II, Proposition 7.11]{Har77}, it holds that 
$(f_1)_*\MO_X(D) \simeq \mathcal E \otimes \mathcal L$. 
If $R_2$ is of type $E_2$ (resp. $E_5$), then we set $e :=1$ (resp. $e:=2$). 
By Theorem \ref{dimy=3}, we have $\MO_X(D)|_D \simeq \MO_{\mathbb P^2}(-e)$. 
Since $f_1 \circ j: D \to \mathbb P^2$ is an isomorphism, we get $(f_1)_*j_*(\MO_X(D)|_D) \simeq \MO_{\mathbb P^2}(-e)$. 
Hence, the exact sequence (\ref{e2-pic2-C2E}) can be written as follows: 
\begin{equation*}\label{e3-pic2-C2E}
0\longrightarrow \MO_{\P^2} \longrightarrow   \mathcal E \otimes \mathcal L
\longrightarrow \MO_{\mathbb P^2}(-e) \longrightarrow 0.
\end{equation*}
By 
$\mathrm{Ext}^1(\MO_{\mathbb{P}^2}(-e),\MO_{\mathbb{P}^2}) \simeq H^1(\P^2, \MO_{\P^2}(e))=0$, 
we obtain $\mathcal E \otimes \mathcal L \simeq \MO_{\P^2} \oplus \MO_{\mathbb P^2}(-e)$. 
Therefore, it hold that 
\[
X = \mathbb P_{\P^2}(\mathcal E) \simeq \mathbb P_{\P^2}(\mathcal E \otimes \mathcal L \otimes \MO_{\P^2}(e)) 
\simeq 
 \mathbb P_{\P^2}(\MO_{\P^2} \otimes \MO_{\P^2}(e)).
\]

What is remaining is to compute  $(-K_X)^3$. 
By Proposition~\ref{standard_2}(3), we have 
\[
(-K_X)^3=2c_1(\MO_{\mathbb{P}^2}\oplus \MO_{\mathbb{P}^2}(e))^2-8c_2(\MO_{\mathbb{P}^2}\oplus \MO_{\mathbb{P}^2}(e)) + 6K_{\mathbb{P}^2}^2.
\]
By \cite[Appendix A]{Har77}, it holds that
\begin{align*}
  c_1(\MO_{\mathbb{P}^2}\oplus \MO_{\mathbb{P}^2}(e)) = c_1(\MO_{\mathbb{P}^2}) + c_1(\MO_{\mathbb{P}^2}(e)) = c_1(\MO_{\mathbb{P}^2}(e)),
\end{align*}
\begin{align*}
  c_2(\MO_{\mathbb{P}^2}\oplus \MO_{\mathbb{P}^2}(e)) = c_1(\MO_{\mathbb{P}^2})\cdot c_1(\MO_{\mathbb{P}^2}(e))
=0.
\end{align*}
Therefore, we obtain 
\begin{align*}
  (-K_X)^3 &= 2c_1(\MO_{\mathbb{P}^2}\oplus \MO_{\mathbb{P}^2}(e))^2-8c_2(\MO_{\mathbb{P}^2}\oplus \MO_{\mathbb{P}^2}(e)) + 6K_{\mathbb{P}^2}^2 \\
  &= 2e^2-8\cdot 0+6\cdot 9 \\
  &= 2e^2+54.
\end{align*}
\end{proof}

\begin{thm}\label{t-pic2-E}
Let $X$ be a  Fano threefold with $\rho(X)=2$. 
Let $R_1$ and $R_2$ be {\cred the} extremal rays of $\NE(X)$. 
Assume that $R_1$ is of type $C$ and $R_2$ is of type $E_2, E_3, E_4$, or $E_5$. 
Then one and only one of the following holds.  
\begin{enumerate}
\item 
$R_1$ is of type $C_1$, $R_2$ is of type $E_3$ or $E_4$, 
$(-K_X)^3 =14$, 
and there exists a split double cover $f: X \to V_7$ such that 
$f_*\MO_X/\MO_{V_7} \simeq \MO_{V_7}(-\frac{1}{2}K_{V_7})^{-1}$ 
{\cred (No.\ \hyperref[table-2-8]{2-8})}. 
\item 
$R_1$ is of type $C_2$, $R_2$ is of type $E_2$, $(-K_X)^3 =56$, 
and $X \simeq \P(\MO_{\P^2} \oplus \MO_{\P^2}(1))$ {\cred (No.\ \hyperref[table-2-35]{2-35})}. 
\item 
$R_1$ is of type $C_2$, $R_2$ is of type $E_5$, $(-K_X)^3 =62$, 
and $X \simeq \P(\MO_{\P^2} \oplus \MO_{\P^2}(2))$ {\cred (No.\ \hyperref[table-2-36]{2-36})}. 
\end{enumerate}
\end{thm}

\begin{proof}
The assertion follows from 
Lemma \ref{r1-c*}, 
{\cred Lemma \ref{l-pic2E-vol},} 
Proposition \ref{p-pic2-C1E-final}, and 
Lemma \ref{l-pic2-C2E}.  
\end{proof}

\subsection{Case $E_1-C$}\label{ss-rho2-E1C}

\begin{lem}\label{l-genus-E1}
We use Notation \ref{n-pic-2}. 
Assume that $R_1$ is of type $E_1$. 
Then 
\[
g(B_1) = \frac{(-K_X)^3}{2} - \frac{(-K_{Y_1})^3}{2} + r_1 \deg B_1 +1
\]
\end{lem}

\begin{proof}
See, e.g.,  \cite[Lemma 3.21(2)(a)]{TanII}. 
\end{proof}

\begin{lem}\label{l-E1C}
We use Notation \ref{n-pic-2}. 
Assume that $R_1$ is of type $E_1$ and $R_2$ is of type $C$. 
Then the following hold. 
\begin{enumerate}
\item[(A)] $18 = \mu_2\left(\frac{24}{r_1} + \deg B_1\right)  + \deg \Delta_{f_2}$. 
\item[(B)] $H_1^2 \cdot H_2 = (r_1-\mu_2)L_1^3 = \frac{1}{\mu_2^2}(8-\deg \Delta_{f_2})$. 
\item[(C)] $H_1\cdot H_2^2 = (r_1-\mu_2)^2 L_1^3 -  \deg B_1 = \frac{2}{\mu_2}$. 
\end{enumerate}
\end{lem}

\begin{proof}
The assertion (A) follows from Lemma \ref{c_times_H} and Proposition \ref{basis}(4). 
The assertion (C) holds by Lemma \ref{l-E1-H1H2} and 
\[
2 =(-K_X) \cdot H_2^2  = (\mu_2 H_1 +H_2) \cdot H_2^2 =\mu_2 H_1 \cdot H_2^2. 
\]
The assertion (B) follows from  Lemma \ref{l-E1-H1H2} and 
\[
12 -\deg \Delta_{f_2} = (-K_X)^2 \cdot H_2  = (\mu_2 H_1+H_2)^2 \cdot H_2 
= \mu_2^2 H_1^2 \cdot H_2 + 2\mu_2 H_1 \cdot H_2^2 
= \mu_2^2 H_1^2 \cdot H_2 + 4, 
\]
where the first equality is guaranteed by Proposition \ref{conic_disc-num}(3).  
\end{proof}

\begin{prop}\label{p-E1C1}
We use Notation \ref{n-pic-2}. 
Assume that $R_1$ is of type $E_1$ and $R_2$ is of type $C_1$. 
Then one of the following holds. 
\begin{enumerate}
\item $(-K_X)^3 = 16, r_1 =4, \deg B_1 = 7, g(B_1) =5, \deg \Delta_{f_2}=5$ (No.\ \hyperref[table-2-9]{2-9}). 
\item $(-K_X)^3 = 20, r_1 =3, \deg B_1 = 6, g(B_1) =2, \deg \Delta_{f_2}=4$ (No.\ \hyperref[table-2-13]{2-13}). 
\item $(-K_X)^3 = 18, r_1 =2, (-K_{{\cred Y_1}})^3 =24, \deg B_1 = 1, g(B_1) =0, \deg \Delta_{f_2}=5$ (No.\ \hyperref[table-2-11]{2-11}).  
\item $(-K_X)^3 = 22, r_1 =2, (-K_{{\cred Y_1}})^3 =32, \deg B_1 = 2, g(B_1) =0, \deg \Delta_{f_2}=4$ (No.\ \hyperref[table-2-16]{2-16}). 
\item $(-K_X)^3 = 26, r_1 =2, (-K_{{\cred Y_1}})^3 =40, \deg B_1 = 3, g(B_1) =0, \deg \Delta_{f_2}=3$ (No.\ \hyperref[table-2-20]{2-20}). 
\end{enumerate}

\end{prop}

\begin{proof}
We have $\mu_2 =1$. 
By Lemma \ref{l-E1C}, the following hold. 
\begin{enumerate}
\item[(A)] $18 = \frac{24}{r_1} + \deg B_1 + \deg \Delta_{f_2}$. 
\item[(B)] $H_1^2 \cdot H_2 = (r_1-1)L_1^3 = 8 - \deg \Delta_{f_2}$. 
\item[(C)] $H_1\cdot H_2^2 = (r_1-1)^2 L_1^3 -  \deg B_1 = 2$. 
\end{enumerate}
{\cred By (C) (or Proposition \ref{p-not-r1}), we obtain $r_1 \in \{2, 3, 4\}$.}  

Assume $r_1 =4$, and hence $L_1^3 =1$. 
Then (B) and (C) imply $\deg \Delta_{f_2}=5$ and $\deg B_1 = 7$, respectively. 
Hence the following hold (Lemma \ref{l-genus-E1}): 
\[
(-K_X)^3 = (H_1+H_2)^3 = H_1^3 + 3H_1^2 \cdot H_2 +3H_1 \cdot H_2^2 +H_2^3 =  L_1^3 + 9 + 6+0  =16
\]
\[
g(B_1) = 8 - 32 + 4 \cdot 7 +1 = 5. 
\]

Assume $r_1 =3$, and hence $L_1^3 =2$. 
Then (B) and (C) imply $\deg \Delta_{f_2}=4$ and $\deg B_1 = 6$, respectively. 
Hence the following hold (Lemma \ref{l-genus-E1}): 
\[
(-K_X)^3 = ( H_1+H_2)^3  = H_1^3 + 3H_1^2 \cdot H_2 +3H_1 \cdot H_2^2 +H_2^3 = 2 + 12 + 6 +0 =20
\]
\[
g(B_1) = 10 - 27 + 3 \cdot 6 +1 = 2. 
\]

Assume $r_1 =2$, and hence $1 \leq L_1^3 \leq 5$. Then  
\begin{enumerate}
\item[(A)] $6 = \deg B_1 + \deg \Delta_{f_2}$. 
\item[(B)] $H_1^2 \cdot H_2 = L_1^3 = 8 - \deg \Delta_{f_2}$. 
\item[(C)] $H_1\cdot H_2^2 =  L_1^3 -  \deg B_1 = 2$. 
\end{enumerate}
By (C), we get $L_1^3 =2 + \deg B_1 \geq 3$. 

Assume $(r_1, L_1^3)=(2, 3)$. 
Then (B) and (C) imply $\deg \Delta_{f_2}=5$ and $\deg B_1 = 1$, respectively. 
In particular, $g(B_1)=0$. 
Hence the following hold: 
\[
(-K_X)^3 = ( H_1+H_2)^3 = H_1^3 + 3H_1^2 \cdot H_2 +3H_1 \cdot H_2^2 +H_2^3  = 3 + 9 + 6 +0 =18
\]

Assume $(r_1, L_1^3)=(2, 4)$. 
Then (B) and (C) imply $\deg \Delta_{f_2}=4$ and $\deg B_1 = 2$, respectively. 
In particular, $g(B_1)=0$. 
Hence the following hold: 
\[
(-K_X)^3 = ( H_1+H_2)^3  = H_1^3 + 3H_1^2 \cdot H_2 +3H_1 \cdot H_2^2 +H_2^3 = 4 + 12 + 6 +0 =22.
\]

Assume $(r_1, L_1^3)=(2, 5)$. 
Then (B) and (C) imply $\deg \Delta_{f_2}=3$ and $\deg B_1 = 3$, respectively. 
Hence the following hold (Lemma \ref{l-genus-E1}): 
\[
(-K_X)^3 = ( H_1+H_2)^3 = H_1^3 + 3H_1^2 \cdot H_2 +3H_1 \cdot H_2^2 +H_2^3= 5 + 15 + 6 +0 =26.
\]
\[
g(B_1) = 13 - 20 +2 \cdot 3+1 = 0. 
\]
\end{proof}

\begin{prop}\label{p-E1C2}
We use Notation \ref{n-pic-2}. 
Assume that $R_1$ is of type $E_1$ and $R_2$ is of type $C_2$. 
Then one of the following holds. 
\begin{enumerate}
\item $(-K_X)^3 = 38, r_1 =4, \deg B_1 = 3, g(B_1) =0$ (No.\ \hyperref[table-2-27]{2-27}). 
\item $(-K_X)^3 = 46, r_1 =3, \deg B_1 = 1, g(B_1) =0$ (No.\ \hyperref[table-2-31]{2-31}). 
\end{enumerate}

\end{prop}

\begin{proof}
We have $\mu_2 =2$ and $\Delta_{f_2}=0$. 
By Lemma \ref{l-E1C}, the following hold. 
\begin{enumerate}
\item[(A)] $9 = \frac{24}{r_1} + \deg B_1$. 
\item[(B)] $H_1^2 \cdot H_2 = (r_1-2)L_1^3 =2$. 
\item[(C)] $H_1\cdot H_2^2 = (r_1-2)^2 L_1^3 -  \deg B_1 = 1$. 
\end{enumerate}
By (B), we obtain $r_1 \in \{3, 4\}$. 
Moreover, (B) and (C) imply
\[
(-K_X)^3 = ( 2H_1+H_2)^3 = 8H_1^3 + 12H_1^2\cdot H_2 + 6H_1 \cdot H_2^2 +H_2^3 = 8L_1^3 + 12 \cdot 2 + 6 \cdot 1 +0 = 8L_1^3 +30. 
\]

Assume $r_1 =4$. 
Then $L_1^3=1$ and $(-K_X)^3 = 8L_1^3 + 30 =38$. 
Moreover, (A) implies $\deg B_1 = 3$. 
Hence the following hold (Lemma \ref{l-genus-E1}): 
\[
g(B_1) = 19 - 32 + 4 \cdot 3 +1 = 0. 
\]

Assume $r_1 =3$. 
Then $L_1^3=2$ and $(-K_X)^3 = 8L_1^3 + 30 =46$. 
Moreover, 
 (A) implies $\deg B_1 = 1$. 
Hence the following hold (Lemma \ref{l-genus-E1}): 
\[
g(B_1) = 23 - 27 + 3 \cdot 1 +1 = 0. 
\]
\end{proof}

\subsection{Case $E_1-D$}\label{ss-rho2-E1D}

\begin{lem}\label{l-E1D}
We use Notation \ref{n-pic-2}. 
Assume that $R_1$ is of type $E_1$ and $R_2$ is of type $D$. 
Set $d_2 := (-K_X)^2 \cdot H_2$. 
Then the following hold. 
\begin{enumerate}
\item[(A)] $12 = \mu_2\left(\frac{24}{r_1} + \deg B_1\right)  -d_2$. 
\item[(B)] $H_1^2 \cdot H_2 = (r_1-\mu_2)L_1^3 = \frac{d_2}{\mu_2^2}$. 
\item[(C)] $H_1\cdot H_2^2 = (r_1-\mu_2)^2 L_1^3 -  \deg B_1 =0$. 
\end{enumerate}
\end{lem}

\begin{proof}
The assertion (A) follows from Lemma \ref{c_times_H} and Proposition \ref{basis}(4). 
The assertion (C) holds by Lemma \ref{l-E1-H1H2} and $H_2^2 \equiv 0$. 
The assertion (B) follows from  Lemma \ref{l-E1-H1H2} and 
\[
d_2 = (-K_X)^2 \cdot H_2  = (\mu_2 H_1+H_2)^2 \cdot H_2 
= \mu_2^2 H_1^2 \cdot H_2. 
\]
\end{proof}

\begin{prop}\label{p-E1D3}
We use Notation \ref{n-pic-2}. 
Assume that $R_1$ is of type $E_1$ and $R_2$ is of type $D_3$. 
Then $(-K_X)^3 = 54, r_1 =4, \deg B_1 = 1, g(B_1) =0$ (No.\ \hyperref[table-2-6]{2-33}). 
\end{prop}

\begin{proof}
We have $d_2 =(-K_X)^2 \cdot H_2 =9$. 
By Lemma \ref{l-E1D}, the following hold: 
\begin{enumerate}
\item[(A)] $7 = \frac{24}{r_1} + \deg B_1$. 
\item[(B)] $H_1^2 \cdot H_2 = (r_1-3)L_1^3 = 1$. 
\item[(C)] $H_1\cdot H_2^2 = (r_1-3)^2 L_1^3 -  \deg B_1 =0$. 
\end{enumerate}
By (B), we get $r_1 =4$ and $L_1^3 =1$. 
Then (C) implies $\deg B_1 = 1$. 
Hence the following holds: 
\[
(-K_X)^3 = ( 3H_1+H_2)^3 = 
27H_1^3 + 27H_1^2 \cdot H_2 + 9 H_1 \cdot H_2^2 + H_{{\cred 2}}^3 = 27 + 27 + 0 +0 =54.
\]
\end{proof}

\begin{prop}\label{p-E1D2}
We use Notation \ref{n-pic-2}. 
Assume that $R_1$ is of type $E_1$ and $R_2$ is of type $D_2$. 
Then one of the following holds. 
\begin{enumerate}
\item $(-K_X)^3 = 32, r_1 =4, \deg B_1 = 4, g(B_1) =1$ (No.\ \hyperref[table-2-6]{2-25}). 
\item $(-K_X)^3 = 40, r_1 =3, \deg B_1 = 2, g(B_1) =0$ (No.\ \hyperref[table-2-6]{2-29}). 
\end{enumerate}

\end{prop}

\begin{proof}
We have $d_2 =(-K_X)^2 \cdot H_2=8$. 
By Lemma \ref{l-E1D}, the following hold: 
\begin{enumerate}
\item[(A)] $10 = \frac{24}{r_1} + \deg B_1$. 
\item[(B)] $H_1^2 \cdot H_2 = (r_1-2)L_1^3 = 2$. 
\item[(C)] $H_1\cdot H_2^2 = (r_1-2)^2 L_1^3 -  \deg B_1 =0$. 
\end{enumerate}
By (B), ${\cred r_1 \in \{3, 4\}}$. Moreover, we get 
\[
(-K_X)^3 = ( 2H_1+H_2)^3 = 8H_1^3 + 12 H_1^2 \cdot H_2 + 6H_1\cdot H_2^2
+H_2^3 = 8L_1^3 + 24. 
\]

Assume $r_1 =4$. 
We then get $L_1^3 =1$ and $(-K_X)^3 = 8+24=32$. 
Moreover, (A) implies $\deg B_1 = 4$. 
Hence the following hold (Lemma \ref{l-genus-E1}): 
\[
g(B_1) = 16 - 32 +4 \cdot 4+1 = 1. 
\]

Assume $r_1 =3$. 
We then get $L_1^3 =2$ and $(-K_X)^3 = 16+24=40$. 
Moreover, (A) implies $\deg B_1 = 2$. 
In particular, $g(B_1)=0$. 
\end{proof}

\begin{prop}\label{p-E1D1}
We use Notation \ref{n-pic-2}. 
Assume that $R_1$ is of type $E_1$ and $R_2$ is of type $D_1$. 
Then one of the following holds. 
\begin{enumerate}
\item $(-K_X)^3 = 10, r_1 =4, \deg B_1 = 9, g(B_1) =10, (-K_X)^2 \cdot H_2 =3$ (No.\ \hyperref[table-2-6]{2-4}). 
\item $(-K_X)^3 = 14, r_1 =3, \deg B_1 = 8, g(B_1) =5, (-K_X)^2 \cdot H_2 =4$ (No.\ \hyperref[table-2-6]{2-7}). 
\item $(-K_X)^3 = 4, r_1 =2, (-K_{Y_1})^3 = 8, \deg B_1 = 1, g(B_1) =1,  (-K_X)^2 \cdot H_2 =1$ (No.\ \hyperref[table-2-1]{2-1}). 
\item $(-K_X)^3 = 8, r_1 =2, (-K_{Y_1})^3 = 16, \deg B_1 = 2, g(B_1) =1,  (-K_X)^2 \cdot H_2 =2$ (No.\ \hyperref[table-2-3]{2-3}). 
\item $(-K_X)^3 = 12, r_1 =2, (-K_{Y_1})^3 = 24, \deg B_1 = 3, g(B_1) =1,  (-K_X)^2 \cdot H_2 =3$ (No.\ \hyperref[table-2-5]{2-5}). 
\item $(-K_X)^3 = 16, r_1 =2, (-K_{Y_1})^3 = 32, \deg B_1 = 4, g(B_1) =1,  (-K_X)^2 \cdot H_2 =4$ (No.\ \hyperref[table-2-10]{2-10}). 
\item $(-K_X)^3 = 20, r_1 =2, (-K_{Y_1})^3 = 40, \deg B_1 = 5, g(B_1) =1,  (-K_X)^2 \cdot H_2 =5$ (No.\ \hyperref[table-2-14]{2-14}). 
\end{enumerate}
\end{prop}

\begin{proof}
By Lemma \ref{l-E1D}, the following hold: 
\begin{enumerate}
\item[(A)] $12 +d_2 = \frac{24}{r_1} + \deg B_1$. 
\item[(B)] $H_1^2 \cdot H_2 = (r_1-1)L_1^3 = d_2$. 
\item[(C)] $H_1\cdot H_2^2 = (r_1-1)^2 L_1^3 -  \deg B_1 =0$. 
\end{enumerate}
{\cred By (C) (or Proposition \ref{p-not-r1}), we obtain $r_1 \in \{2, 3, 4\}$.}  
It holds that 
\[
(-K_X)^3 = ( H_1+H_2)^3 = 
H_1^3 +3H_1^2 \cdot H_2 + 3H_1 \cdot H_2^2 + H_3^3 = L_1^3 + 3 d_2. 
\]

Assume $r_1 =4$, and hence $L_1^3 =1$. 
Then (B) and (C) imply $d_2 = 3$ and $\deg B_1 =9$, respectively. 
In particular, $(-K_X)^3 = 1 + 3 \cdot 3 = 10$. 
Hence the following hold (Lemma \ref{l-genus-E1}): 
\[
g(B_1) = 5 - 32 +4 \cdot 9+1 = 10. 
\]

Assume $r_1 =3$, and hence $L_1^3 =2$. 
Then (B) and (C) imply $d_2 = 4$ and $\deg B_1 =8$, respectively. 
In particular, $(-K_X)^3 = 2 + 3 \cdot 4 = 14$. 
Hence the following hold (Lemma \ref{l-genus-E1}): 
\[
g(B_1) = 7 - 27 +3 \cdot 8+1 = 5.  
\]

Assume $r_1 =2$, and hence $1 \leq L_1^3 \leq 5$. 
Then 
\begin{enumerate}
\item[(A)] $d_2 =  \deg B_1$. 
\item[(B)] $H_1^2 \cdot H_2 = L_1^3 = d_2$. 
\item[(C)] $H_1\cdot H_2^2 =  L_1^3 -  \deg B_1 =0$. 
\end{enumerate}
We get  $(-K_X)^3 = L_1^3 + 3d_2 = 4L_1^3$. 
Then the following hold (Lemma \ref{l-genus-E1}): 
\[
g(B_1) = 2L_1^3 - 4L_1^3 +2 \cdot L_1^3 +1 = 1.  
\]

\end{proof}


\begin{lem}\label{l-dP-fib-CI}
Let $\sigma: X \to Y$ be a blowup along a smooth curve $\Gamma$ on a Fano threefold $Y$. 
Let $\pi:X \to \P^1$ be a morphism with $\pi_*\MO_X = \MO_{\P^1}$. 
Take a Cartier divisor $D$ on $Y$. 
Assume that {\rm (i)} and {\rm (ii)} hold. 
\begin{enumerate}
\renewcommand{\labelenumi}{(\roman{enumi})}
\item $(-K_Y) \cdot D^2= (-K_Y) \cdot \Gamma$. 
\item 
$D \sim \sigma_*F$ for a fibre $F$ of $\pi : X \to \P^1$. 
\end{enumerate}
Then $\Gamma$ is a complete intersection of two members of $|D|$. 
\end{lem}


\begin{proof}
Take two fibres $F_1$ and $F_2$ of $\pi : X \to \P^1$. 
Set $D_{1} := \sigma_*F_1$ and $D_{2} := \sigma_*F_2$. 
For each $i \in \{1, 2\}$, 
we get 
$D_{i} = \sigma_*F_i \sim \sigma_* F \overset{{\rm (ii)}}{\sim} D$. 
We then obtain $\Gamma \supset (D_{1} \cap D_{2})_{\red}$. 
Since $D_{1} \cap D_{2}$ is of pure {\cred dimension one}, 
we get $(D_{1} \cap D_{2})_{\red} = \Gamma$ or $(D_{1} \cap D_{2})_{\red} = \emptyset$. 
We then obtain $(D_{1} \cap D_{2})_{\red} = \Gamma$, 
as otherwise the equality $(D_{1} \cap D_{2})_{\red} = \emptyset$ would lead to  the following contradiction: 
\[
0 = (-K_Y) \cdot (D_{1} \cap D_{2}) = (-K_Y) \cdot D_{1} \cdot D_{2} 
= (-K_Y) \cdot D^2 \overset{{\rm (i)}}{=} (-K_Y) \cdot \Gamma >0. 
\]

It is enough to show the scheme-theoretic equality  $D_{1} \cap D_{2} = \Gamma$. 
Recall that we have the equality 
$(D_{1} \cap D_{2})_{\red} = \Gamma$. 
Since $D_{1} \cap D_{2}$ is Cohen-Macaulay, 
it suffices to prove that $D_{1} \cap D_{2}$ is $R_0$, i.e., 
$\MO_{D_{1} \cap D_{2}, \xi}$ is reduced for the generic point $\xi$ of $D_{1} \cap D_{2}$. 
This follows from 
$(-K_Y) \cdot (D_{1} \cap D_{2}) = (-K_Y) \cdot D^2 \overset{{\rm (i)}}{=} -K_Y \cdot \Gamma$, 
because the effective $1$-cycle associated with $D_{1} \cap D_{2}$ is determined by its length at the generic point \cite[Lemma 1.18]{Bad01}. 
\end{proof}

\begin{prop}\label{p-pic2-dP-CI}
Let $X$ be a Fano threefold with $\rho(X)=2$ which has extremal rays of type $D$ and $E_1$. 
Let $\sigma: X \to Y$ be the contraction of the extremal ray of type $E_1$. 
Then the blowup centre $\sigma(\Ex(\sigma))$ of $\sigma$ is a complete intersection $D_1 \cap D_2$ of two members $D_1$ and $D_2$ of $|D|$ 
for some  Cartier divisor $D$ on $Y$. 
\end{prop}

In the following proof, 
we use No.\ of $X$ 
given in Proposition \ref{p-E1D3}, Proposition \ref{p-E1D2}, and Proposition \ref{p-E1D1}. 
In particular, No.\ of $X$ is one of 
2-1, 2-3, 2-4, 2-5, 2-7, 2-10, 2-14, 
2-25, 2-29, 2-33.

\begin{proof}
Set $\Gamma := \sigma(\Ex(\sigma))$. 
Let $\pi : X \to \P^1$ be the contraction of the extremal ray of type $D$. 
Take a fibre $F$ of $\pi$. 
Set $\mu$ to be the length of $\pi$. 
Then $\pi : X \to \P^1$ is of type  $D_{\mu}$ (Remark \ref{r-length}(2)). 
For each case, it is enough to find a Cartier divisor $D$ on $Y$ 
satisfying  (i) and (ii) of Lemma \ref{l-dP-fib-CI}. 
In what follows, we treat the following three cases separately. 
\begin{enumerate}
\item 2-4, 2-25, 2-33. 
\item 2-7, 2-29.
\item 2-1, 2-3, 2-5, 2-10, 2-14.
\end{enumerate}

(1) 
In these cases, we have $Y=\P^3$ and $\deg \Gamma = s^2$ for some $1 \leq s \leq 3$ 
(Proposition \ref{p-E1D3}, Proposition \ref{p-E1D2}, Proposition \ref{p-E1D1}). 
Let $D$ be a Cartier divisor satisfying $\MO_{\P^3}(D) \simeq \MO_{\P^3}(s)$. 
Then Lemma \ref{l-dP-fib-CI}(i) holds by 
$(-K_Y) \cdot D^2 = 4s^2 = (-K_Y) \cdot \Gamma$. 
It holds that $-K_X \sim \sigma^*\MO_{\P^3}(\mu) + F$ (Proposition \ref{basis}(3)). 
Hence $\sigma_*F \sim -K_Y -\MO_{\P^3}(\mu) \sim \MO_{\P^3}(4-\mu)$. 
In order to prove Lemma \ref{l-dP-fib-CI}(ii), it suffices to show $s+\mu=4$. 
By using the fact that $\pi$ is of type $D_{\mu}$, 
this can be checked for all the cases 2-4, 2-25, 2-33 (Proposition \ref{p-E1D3}, Proposition \ref{p-E1D2}, Proposition \ref{p-E1D1}).

(2) 
In these cases, we have $Y=Q$ and $\deg \Gamma = 2s^2$ for some $1 \leq s \leq 2$ (Proposition \ref{p-E1D2}, Proposition \ref{p-E1D1}). 
Let $D$ be a Cartier divisor satisfying $\MO_{Q}(D) \simeq \MO_{Q}(s)$. 
Then Lemma \ref{l-dP-fib-CI}(i) holds by 
$(-K_Y) \cdot D^2 = 6s^2 = (-K_Y) \cdot \Gamma$. 
It holds that $-K_X \sim \sigma^*\MO_{Q}(\mu) + F$ (Proposition \ref{basis}(3)). 
Hence $\sigma_*F \sim -K_Y -\MO_{Q}(\mu) \sim \MO_{Q}(3-\mu)$. 
In order to prove Lemma \ref{l-dP-fib-CI}(ii), it suffices to show $s+\mu=3$. 
By using the fact that $\pi$ is of type $D_{\mu}$, 
this can be checked for both cases 2-7, 2-29  (Proposition \ref{p-E1D2}, Proposition \ref{p-E1D1}).

(3) 
In these cases, 
$Y$ is a Fano threefold of index $r_Y=2$ and 
$\Gamma$ is an elliptic curve such that $d = D^3$ and $1 \leq d \leq 5$ for $d := \deg \Gamma$, 
where $D$ is a Cartier divisor satisfying $-K_Y \sim 2D$ 
(Proposition \ref{p-E1D1}). 
Recall that the degree $\deg \Gamma$  is defined as $\frac{1}{2}(-K_Y) \cdot \Gamma$. 
Then Lemma \ref{l-dP-fib-CI}(i) holds by 
$(-K_Y) \cdot D^2 = 2D^3 =2d = (-K_Y) \cdot \Gamma$. 
Since $\pi$ is of type $D_1$ for all the cases 2-1, 2-3, 2-5, 2-10, 2-14  
 (Proposition \ref{p-E1D1}), 
it holds that $-K_X \sim \sigma^*(\mu D) + F 
\sim \sigma^*D+F$  (Proposition \ref{basis}(3)). 
Hence $\sigma_*F \sim -K_Y -D \sim D$. 
Thus  Lemma \ref{l-dP-fib-CI}(ii) holds. 
\qedhere 





\end{proof}

\subsection{Case $E_1-E$}\label{ss-rho2-E1E}

\begin{lem}\label{l-E1E}
We use Notation \ref{n-pic-2}. 
Assume that $R_1$ is of type $E_1$ and 
$R_2$ is of type $E_2, E_3,$ or $E_4$. 
Then the following hold. 
\begin{enumerate}
\item[(A)] $24 = \mu_2( \frac{24}{r_1} + \deg B_1) + \frac{24}{r_2}$. 
\item[(B)] $H_1^2 \cdot H_2 = (\frac{r_2-1}{\mu_2})^2 L_2^3 = 
(r_1-\mu_2)L_1^3$. 
\item[(C)] $H_1 \cdot H_2^2 = \frac{r_2-1}{\mu_2}L_2^3 = (r_1-\mu_2)^2 L_1^3 - \deg B_1$. 
\end{enumerate}
\end{lem}

\begin{proof}
The assertion (A) follows from Lemma \ref{c_times_H}(3) and Proposition \ref{basis}(4). 
For each $i \in \{1, 2\}$, 
let $a_i$ be the positive integer satisfying 
\[
K_X \sim f_i^*K_{Y_i} +a_iD_i. 
\]
Specifically, we set 
$a_1 :=1$ and $a_2 := 2$ (resp. $a_2 :=1)$ if $R_2$ is of type $E_2$ (resp. $E_3$ or $E_4$) \cite[Proposition 3.22]{TanII}.
For each $i \in \{1, 2\}$. we have 
\[
\mu_2 H_1 + H_2 \sim -K_X \sim -f^*_iK_{Y_i} -a_i D_i =r_i H_i -a_iD_i. 
\]
Then the assertion (B) holds by $H_1^2 \cdot H_2 = (r_1-\mu_2)L_1^3$ (Lemma \ref{l-E1-H1H2}) and 
\[
\mu_2^2H_1^2 \cdot H_2 = ( (r_2-1)H_2 -a_2D_2)^2 \cdot H_2 
= (r_2-1)^2 L_2^3. 
\]
The assertion (C) follows from 
$H_1 \cdot H_2^2 = (r_1-\mu_2)^2 L_1^3 - \deg {\cred B_1}$ (Lemma \ref{l-E1-H1H2}) 
and 
\[
\mu_2 H_1 \cdot  H_2^2 = ((r_2 -1)H_2 -a_2D_2) \cdot H_2^2 = (r_2-1)L_2^3. 
\]
\end{proof}

\subsubsection{Case $E_1-E_2$}

\begin{prop}\label{p-E1E2}
    We use Notation \ref{n-pic-2}. 
Assume that $R_1$ is of type $E_1$ and $R_2$ is of type $E_2$. 
Then $(-K_X)^3 = 46, r_1 =4, \deg B_1 = 2, g(B_1) =0, r_2 =3$ (No.\ \hyperref[table-2-30]{2-30}). 
\end{prop}

\begin{proof}
For each $i \in \{1, 2\}$, $Y_i$ is smooth and $r_i$ is the index of the Fano threefold $Y_i$. 
In particular, $r_i \in \{1, 2, 3, 4\}$. 
We have $\mu_2 =2$.  
Hence Lemma \ref{l-E1E} implies the following: 
\begin{enumerate}
\item[(A)] $24 =  \frac{48}{r_1} + 2\deg B_1 + \frac{24}{r_2}$. 
\item[(B)] $H_1^2 \cdot H_2 = (\frac{r_2-1}{2})^2 L_2^3 = 
(r_1-2)L_1^3$. 
\item[(C)] $H_1 \cdot H_2^2 = \frac{r_2-1}{2}L_2^3 = (r_1-2)^2 L_1^3 - \deg B_1$. 
\end{enumerate}
Since each of $H_1$ and $H_2$ is nef and big, we have $H_1^2 \cdot H_2 >0$ and $H_1 \cdot H_2^2 >0$. By (B) and (C), we obtain $r_1 >2$ and $r_2 >1$, respectively. 
In particular, $r_1 \in \{ 3, 4\}$ and $r_2 \in \{2, 3, 4\}$. 

Assume $r_1 =4$, and hence $L_1^3 =1$. 
By (B) and (C), the following hold: 
\[
r_2 =3,\quad
L_2^3=2, \quad 
H_1^2 \cdot H_2 = 2, \quad H_1 \cdot H_2^2 = 2, 
\quad \deg B_1 = 2. 
\]
Then $g(B_1)=0$ and 
\[
(-K_X)^3 = (2H_1+H_2)^3 = 8 + 24 + 12 + 2 = 46. 
\]

Assume $r_1 =3$, and hence $L_1^3 =2$. 
By (B) and (C), the following hold: 
\[
r_2 =3,\quad
L_2^3=2, \quad 
H_1^2 \cdot H_2 = 2, \quad 
H_1 \cdot H_2^2 = 2, 
\quad \deg B_1 = 0. 
\]
This is absurd. 
\end{proof}

\subsubsection{Cases $E_1-E_3$ and $E_1-E_4$}

\begin{prop}\label{p-E1E3E4}
    We use Notation \ref{n-pic-2}. 
Assume that $R_1$ is of type $E_1$ and $R_2$ is of type $E_3$ or $E_4$. 
Then one of the following holds. 
\begin{enumerate}
\item 
$(-K_X)^3 = 22, r_1 =4, \deg B_1 = 6, g(B_1) =4, r_2 =2, L_2^3 = 3$ 
(No.\ \hyperref[table-2-15]{2-15}). 
\item 
$(-K_X)^3 =30, r_1 =3, \deg B_1 = 4, g(B_1) =1, r_2 =2, L_2^3 = 4$ 
(No.\ \hyperref[table-2-23]{2-23}). 
\end{enumerate}
\end{prop}

\begin{proof}
Note that $Y_2$ is not smooth. 
We have $\mu_2 =1$. 
Hence Lemma \ref{l-E1E} implies the following: 
\begin{enumerate}
\item[(A)] $24 =  \frac{24}{r_1} + \deg B_1 + \frac{24}{r_2}$. 
\item[(B)] $H_1^2 \cdot H_2 = (r_2-1)^2 L_2^3 = 
(r_1-1)L_1^3$. 
\item[(C)] $H_1 \cdot H_2^2 = (r_2-1)L_2^3 = (r_1-1)^2 L_1^3 - \deg B_1$. 
\end{enumerate}
Since each of $H_1$ and $H_2$ is nef and big, we have $H_1^2 \cdot H_2 >0$ and $H_1 \cdot H_2^2 >0$. By (B) and (C), 
we obtain $r_1 \geq 2$ and $r_2 \geq 2$, respectively. 

Assume $r_1 =4$, and hence $L_1^3 =1$. 
By (B) and (C), the following hold: 
\[
r_2 =2,\quad
L_2^3=3, \quad 
H_1^2 \cdot H_2 = 3, \quad H_1 \cdot H_2^2 = 3, 
\quad \deg B_1 = 6. 
\]
It follows from Lemma \ref{l-genus-E1} that  
\[
(-K_X)^3 = (H_1+H_2)^3 = 1 + 9 + 9 + 3 = 22, 
\]
\[
g(B_1) = 11 - 32 + 24 + 1= 4. 
\]

Assume $r_1 =3$, and hence $L_1^3 =2$. 
By (B), we get $(r_2-1)^2 L_2^3 = 4$. 
Thus $(r_2, L_2^3) \in \{(2, 4), (3, 1)\}$. 
Suppose that $(r_2, L_2^3) = (3, 1)$. 
Then (B) and (C) imply $H_1^2 \cdot H_2 \in 2\Z$ and $H_1 \cdot H_2^2 \in 2\Z$. 
We get the following contradiction: 
$(-K_X)^3 =(H_1+H_2)^3 \equiv L_1^3 +L_2^3 \not\equiv 0 \mod 2$. 
Thus $(r_2, L_2^3) = (2, 4)$. 
Then (B) and (C) imply the following: 
\[
H_1^2 \cdot H_2 = 4, \quad 
H_1 \cdot H_2^2 = 4, 
\quad \deg B_1 = 4. 
\]
Then the following hold (Lemma \ref{l-genus-E1}):  
\[
(-K_X)^3 = (H_1+H_2)^3 = 2 + 12 + 12 + 4 = 30, 
\]
\[
g(B_1) = 15 - 27 + 12 +1 = 1. 
\]

Assume $r_1 =2$, and hence $1 \leq L_1^3 \leq 5$. 
Then (A)-(C) can be rewritten as follows: 
\begin{enumerate}
\item[(A)] $12 =  \deg B_1 + \frac{24}{r_2}$. 
\item[(B)] $H_1^2 \cdot H_2 = (r_2 -1)^2L_2^3 = L_1^3$. 
\item[(C)] $H_1 \cdot H_2^2 = (r_2-1)L_2^3 = L_1^3 - \deg B_1$. 
\end{enumerate}
By $1 \leq L_1^3 \leq 5$ and (B), 
we obtain $r_2 \in \{2, 3\}$. 
Since the case $r_2 =2$ contradicts (A), 
we get $r_2 =3$. 
Then (B) and $1 \leq L_1^3 \leq 5$ imply $L_2^3 =1$ and $L_1^3 =4$. 
By (B) and (C), we get $H_1^2 \cdot H_2 \in 2\Z$ and $H_1 \cdot H_2^2 \in 2\Z$, 
respectively. 
We then obtain the following contradiction: 
$(-K_X)^3 = (H_1+H_2)^3 \equiv L_1^3 +L_2^3 \not\equiv 0 \mod 2$. 
\qedhere 


\end{proof}

Although the following result is not used in this paper, 
we will need it in the classification of Fano threefolds with $\rho \geq 3$. 
{\cred For the definition of No.\ 2-xx, we refer to Definition \ref{d No pic2}.} 

\begin{prop}\label{p-2-15-2-23}
\begin{enumerate}
    \item Let $X$ be a Fano threefold of No.\ 2-15. Let $f: X \to \P^3$ be a blowup along a smooth curve $B$ of degree $6$. Then $B$ is contained in a (possibly singular) quadric surface $S$ on $\P^3$. 
    \item 
    Let $X$ be a Fano threefold of No.\ 2-23. Let $f: X \to Q$ be a blowup along a smooth curve $B$ of degree $4$. Then $B$ is contained in a (possibly singular) quadric surface $S$ on $Q$, i.e., $S$ is a hyperplane section of $Q \subset \P^4$. 
\end{enumerate}
\end{prop}

\begin{proof}
We prove (1) and (2) simultaneously. 
If $X$ is 2-15 (resp. 2-23), then 
we set $Y := \P^3$ (resp. $Y:=Q$) and $s :=4$ (resp. $s :=3$). 
Then we have a blowup $f: X \to Y$ and $Y$ is a Fano threefold of index $s$. 
We have the contraction $g: X \to Z$ of the other extremal ray, which is of type $E_3$ or $E_4$. 
In particular, $T :=\Ex(g)$ is a (possibly singular) quadric surface. 
Set $S := f(T)$. 
It is enough to show that $B \subset S$ and 
the induced 
morphism 
\[
f_T :T \to S
\]
is an isomorphism. 
Set $H_f$ (resp. $H_g$) to be the pullback on $X$ of the ample generator of $Y$ (resp. $Z$). 
It holds that $-K_X \sim H_f + H_g$ (Proposition \ref{basis}). 
Let $r$ be the positive integer satisfying $-g^*K_{{\cred Z}} \sim r H_g$. 

\setcounter{step}{0}

\begin{step}\label{s1-2-15-2-23}
$B \subset S$. 
\end{step}

\begin{proof}[Proof of Step \ref{s1-2-15-2-23}]
Suppose $B \not\subset S$. 
Then $T = f^*S$. 
Pick a curve $C$ on $S$ disjoint from $S \cap B$. 
Set $C_X :=f^{-1}(C) \xrightarrow{\simeq} C$. 
Since $S$ is ample, we have $S \cdot C>0$. 
On the other hand, $\MO_X(-T)|_T$ is ample, and hence $T \cdot C_X = \MO_X(T)|_T \cdot C <0$. 
By $T = f^*S$, we get the following contradiction: 
\[
0 < S \cdot C = f^*S \cdot C_X = T \cdot C_X <0. 
\]
This completes the proof of Step \ref{s1-2-15-2-23}. 
\end{proof}

\begin{step}\label{s2-2-15-2-23}
$T \cdot \zeta =r-1$ for  a one-dimensional fibre $\zeta$ of $f: X \to Y$. 
\end{step}

\begin{proof}[Proof of Step \ref{s2-2-15-2-23}]
{\cred For $E := \Ex(f)$}, we have that 
\[
-K_X \sim H_f + H_g \equiv \frac{-K_X+E}{s} + \frac{-K_X+T}{r}, 
\]
which implies 
\[
-(sr-r-s)K_X \sim rE + sT. 
\]
Taking the intersection with $\zeta$, we get 
$sr-r-s  = -r +s T \cdot \zeta$, 
which implies $T \cdot \zeta = r-1$. 
This completes the proof of Step \ref{s2-2-15-2-23}. 
\end{proof}

\begin{step}\label{s3-2-15-2-23}
It holds that 
$r=2$ and $T \cdot \zeta =1$. 
\end{step}

\begin{proof}[Proof of Step \ref{s3-2-15-2-23}]
By \cite[Proposition 3.22]{TanII}, we get 
\[
K_Z^3 = (g^*K_Z)^3 =(K_X-T)^3 = K_X^3 -3 K_X^2 \cdot T +3 K_X \cdot T^2 -T^3 = 
K_X^3 -3 \cdot 2 +3 \cdot 2 -2, 
\]
and hence $r^3 H_g^3 = (-K_Z)^3 = (-K_X)^3 +2$. 
If $X$ is of No.\ 2-15, then we get $r^3 H_g^3 = 22 + 2 = 24$.
If $X$ is of No.\ 2-23, then we get $r^3 H_g^3 = 30 + 2 = 32$. 
Since the inequality $r \geq 2$ holds by Step \ref{s2-2-15-2-23}, 
we obtain $r=2$. 
Again by Step \ref{s2-2-15-2-23}, we get $T \cdot \zeta =r-1=1$. 
This completes the proof of Step \ref{s3-2-15-2-23}. 
\end{proof}

\begin{step}\label{s4-2-15-2-23}
$f_T: T \to S$ is an isomorphism. 
\end{step}

\begin{proof}[Proof of Step \ref{s4-2-15-2-23}]
Let us show that $S$ is normal. Consider an exact sequence 
\[
0 \to \MO_X(-T) \to \MO_X \to \MO_T \to 0. 
\]
By {\cred $E \cdot \zeta <0$ and}
$(-T-K_X) \cdot \zeta =-1+1=0$ (Step \ref{s3-2-15-2-23}), 
we get $R^1f_*\MO_X(-T)=0$, 
{\cred because \cite[Theorem 0.5]{Tan15} is applicable 
for $A := -T-K_X - \frac{1}{2}E$}. 
Then $\MO_S = f_*\MO_T$, i.e., $S$ is normal.

Note that any curve $\zeta$ on $T$ is not contracted by $f$, 
as otherwise we would get a contradiction: $\zeta \in R_1 \cap R_2 =\{0\}$. 
Hence $T \neq \Ex(f)$, which implies that ${\cred f_T}: T \to S$ is birational. 
Since $T$ is a (possibly singular) quadric surface, ${\cred f_T} : T \to S$ is an isomorphism. 
This completes the proof of Step \ref{s4-2-15-2-23}. 
\end{proof}
Step \ref{s1-2-15-2-23} and 
Step \ref{s4-2-15-2-23} complete the proof of Proposition \ref{p-2-15-2-23}. 
\end{proof}

\subsubsection{Case $E_1-E_5$}

\begin{prop}\label{p-E1E5}
We use Notation \ref{n-pic-2}. 
Assume that $R_1$ is of type $E_1$ and $R_2$ is of type $E_5$. 
Then 
$(-K_X)^3 = 40, r_1 =4, \deg B_1 = 3, g(B_1) =1, r_2 =3, L_2^3 = 12$ 
(No.\ \hyperref[table-2-28]{2-28}). 
\end{prop}

\begin{proof}
By Lemma \ref{c_times_H}, we obtain $\mu_2 =1$ and $H_2 \cdot c_2(X) = 45/r_2$. 
It holds that 
\[
H_1 + H_2 \sim -K_X = -{\cred f_1^*}K_{{\cred Y_1}} -D_1 =r_1 H_1 -D_1, 
\]
\[
H_1 + H_2 \sim -K_X = -{\cred f_2^*}K_{{\cred Y_2}} -\frac{1}{2} D_2 
=\frac{r_2}{2} H_2 -\frac{1}{2}D_2. 
\]
We then get the following: 
\[
H_1^2 \cdot H_2 =  H_1^2  \cdot( (r_1-1)H_1 -D_1) = (r_1-1)L_1^3. 
\]
\[
H_1^2 \cdot H_2 = ( (\frac{r_2}{2}-1)H_2 -\frac{1}{2}D_2)^2 \cdot H_2 = (\frac{r_2}{2}-1)^2 L_2^3. 
\]
\[
H_1 \cdot H_2^2 = ((\frac{r_2}{2} -1)H_2 -\frac{1}{2}D_2) \cdot H_2^2 
= (\frac{r_2}{2}-1)L_2^3. 
\]
\[
H_1 \cdot H_2^2 = H_1\cdot ( (r_1-1) H_1 -D_1)^2 
=(r_1-1)^2 L_1^3 +H_1 \cdot D_1^2 =(r_1-1)^2 L_1^3 - \deg {\cred B_1}. 
\]
Therefore, the following hold, 
where (A) is guaranteed by Lemma \ref{c_times_H}(3) and Proposition \ref{basis}(4): 
\begin{enumerate}
\item[(A)] $24 =  \frac{24}{r_1} + \deg B_1 + \frac{45}{r_2}$. 
\item[(B)] $H_1^2 \cdot H_2= (\frac{r_2}{2}-1)^2 L_2^3 =  (r_1-1)L_1^3 $. 
In particular, $L_2^3 = \frac{4(r_1-1)}{(r_2-2)^2} L_1^3$. 
\item[(C)] $H_1 \cdot H_2^2 = (\frac{r_2}{2}-1)L_2^3 = (r_1-1)^2 L_1^3 - \deg B_1$. 
\end{enumerate}
By (A), we have $r_1 \geq 2$  and $r_2 \geq 3$. 
{\cred In particular, $r_1 \in \{2, 3, 4\}$.} 

Assume $r_1 =4$, and hence $L_1^3 =1$. 
By (B), we get $\Z \ni L_2^3 = \frac{4(r_1-1)}{(r_2-2)^2} L_1^3 
= \frac{12}{(r_2-2)^2}$. Hence $r_2 \in \{3, 4\}$. 
Since $r_2 \neq 4$ by (A), 
we obtain $r_2 =3$. 
Then (A)-(C) imply the following: 
\[
L_2^3 =12, \quad 
\deg B_1 = 3, \quad 
H_1^2 \cdot H_2 = 3, \quad H_1\cdot H_2^2 = 6. 
\]
Then the following hold (Lemma \ref{l-genus-E1}): 
\[
(-K_X)^3 = (H_1+H_2)^3 = 1 + 9 + 18 + 12 = 40, 
\]
\[
g(B_1) = 20 -32 + 12 +1 = 1. 
\]

Assume $r_1 =3$, and hence $L_1^3 =2$. 
By (A), we have $r_2 \in \{3, 5, 9, 15, 45\}$. 
This, together with $\Z_{>0} \ni L_2^3 =\frac{4(r_1-1)}{(r_2-2)^2} L_1^3 = \frac{16}{(r_2-2)^2}$, 
implies $r_2 =3$ and $L_2^3 =16$. 
By (C), we get the following contradiction: 
\[
\deg B_1 = (r_1-1)^2 L_1^3 - (\frac{r_2}{2}-1)L_2^3  =8 - 8 =0. 
\]

Assume $r_1 =2$, and hence $1 \leq L_1^3 \leq 5$. 
Then (A)-(C) can be rewritten as follows: 
\begin{enumerate}
\item[(A)] $12 =  \deg B_1 + \frac{45}{r_2}$. 
\item[(B)] $H_1^2 \cdot H_2= (\frac{r_2}{2}-1)^2 L_2^3 =  L_1^3 $. 
In particular, $L_2^3 = \frac{4}{(r_2-2)^2} L_1^3$. 
\item[(C)] $H_1 \cdot H_2^2 = (\frac{r_2}{2}-1)L_2^3 =L_1^3 - \deg B_1$. 
\end{enumerate}
By (A), we have $r_2 \in \{5, 9, 15, 45\}$. 
However, this contradicts $\Z_{>0} \ni L_2^3 =\frac{4}{(r_2-2)^2} L_1^3$ 
and $1 \leq L_1^3 \leq 5$. 
\end{proof}

\subsubsection{Case $E_1-E_1$}

\begin{prop}\label{p-E1E1}
We use Notation \ref{n-pic-2}. 
Assume that $R_1$ is of type $E_1$ and $R_2$ is of type $E_1$. 
Moreover, suppose $r_1 \geq r_2$. 
Then one of the following holds possibly after permuting 
$R_1$ and $R_2$. 
\begin{enumerate}
\item $(-K_X)^3 = 20, Y_1 \simeq Y_2 \simeq  \P^3, g(B_1) = g(B_2) =3, \deg B_1 = \deg B_2  =6$ (No.\ \hyperref[table-2-12]{2-12}). 
\item $(-K_X)^3 = 24, Y_1 \simeq \P^3, g(B_1) = 1, \deg B_1 =5, 
Y_2 \simeq Q, g(B_2) = 1, \deg B_2 =5$ (No.\ \hyperref[table-2-17]{2-17}).  
\item $(-K_X)^3 = 26, Y_1 \simeq \P^3, g(B_1) = 2, \deg B_1 =5, 
r_2 =2, L_2^3 =4, g(B_2) = 0, \deg B_2 =1$ (No.\ \hyperref[table-2-19]{2-19}).  
\item $(-K_X)^3 = 30, Y_1 \simeq \P^3, g(B_1) = 0, \deg B_1 =4, 
r_2 =2, L_2^3 =5, g(B_2) = 0, \deg B_2 =2$ (No.\ \hyperref[table-2-22]{2-22}). 
\item $(-K_X)^3 = 28, Y_1 \simeq Y_2 \simeq Q, 
g(B_1) =g(B_2) =0, \deg B_1 =\deg B_2 = 4$ (No.\ \hyperref[table-2-21]{2-21}). 
\item $(-K_X)^3 = 34, Y_1 \simeq Q, g(B_1) = 0, \deg B_1 =3, 
r_2 =2, L_2^3 =5, g(B_2) = 0, \deg B_2 =1$ (No.\ \hyperref[table-2-26]{2-26}). 
\end{enumerate}
\end{prop}

\begin{proof}
We may apply Lemma \ref{l-E1-H1H2} for both the extremal rays $R_1$ and $R_2$, 
although we need to permute the indices when we apply it for $R_2$. 
Then the following hold (Lemma \ref{c_times_H}(3), Proposition \ref{basis}(4)): 
\begin{enumerate}
\item[(A)] $24 = \frac{24}{r_1} + \deg B_1  + \frac{24}{r_2} + \deg B_2$. 
\item[(B)] $H_1^2 \cdot H_2  = (r_1-1)L_1^3 = (r_2-1)^2 L_2^3-\deg B_2$. 
\item[(C)] $H_1\cdot H_2^2 = (r_2-1)L_2^3 = (r_1-1)^2 L_1^3 -  \deg B_1$. 
\end{enumerate}
Since each of $H_1$ and $H_2$ is nef and big, 
we have $H_1^2 \cdot H_2 >0$ and $H_1 \cdot H_2^2 >0$. 
By (B) and (C), we get $r_1 \geq 2$ and $r_2 \geq 2$, respectively. 
By $r_1 \geq r_2$ and (A), 
we obtain $r_1 \geq 3$. 
In particular, 
\[
(r_1, r_2) \in \{ (4, 4), (4, 3), (4, 2), (3, 3), (3, 2)\}. 
\]
In what follows, we shall use 
\begin{itemize}
\item $(-K_X)^3 = (H_1+H_2)^3 = L_1^3 + 3 H_1^2 \cdot H_2 + 3 H_1 \cdot H_2^2 + L_2^3$, and 
\item $g(B_i) = \frac{(-K_X)^3}{2} - \frac{(-K_{Y_i})^3}{2} + r_i \deg B_i +1$ (Lemma \ref{l-genus-E1}). 
\end{itemize}

\medskip

(1) 
Assume $(r_1, r_2) = (4, 4)$. Then $(L_1^3, L_2^3)=(1, 1)$ and the following hold: 
\begin{enumerate}
\item[(A)] $12 =  \deg B_1 + \deg B_2$. 
\item[(B)] $H_1^2 \cdot H_2 = 3 = 9 -  \deg B_2$. 
\item[(C)] $H_1 \cdot H_2^2= 3 = 9 - \deg B_1$. 
\end{enumerate}
Hence $\deg B_1 = \deg B_2 =6$, 
\[
(-K_X)^3 = 1 + 9+9+1=20, \qquad g(B_1) =g(B_2) = 10 -32 + 24 + 1 =3. 
\]

\medskip

(2) 
Assume $(r_1, r_2) = (4, 3)$. Then $(L_1^3, L_2^3) = (1 , 2)$ and 
the following hold: 
\begin{enumerate}
\item[(A)] $10 =  \deg B_1 + \deg B_2$. 
\item[(B)] $H_1^2\cdot H_2 = 3 = 8 -  \deg B_2$. 
\item[(C)] $H_1 \cdot H_2^2 = 4 = 9 - \deg B_1$. 
\end{enumerate}
Hence $(\deg B_1, \deg B_2, H_1^2\cdot H_2, H_1 \cdot H_2^2) = (5, 5, 3, 4)$, 
\[
(-K_X)^3 = 1 + 9+12+2=24,
\]
\[
g(B_1) = 12 -32 + 20 + 1 =1, \qquad 
g(B_2) = 12 -27 + 15 + 1 =1. 
\]

\medskip

(3), (4) 
Assume $(r_1, r_2) = (4, 2)$. Then $L_1^3=1$, $1 \leq L_2^3 \leq 5$,  and 
the following hold: 
\begin{enumerate}
\item[(A)] $6 =  \deg B_1 + \deg B_2$. 
\item[(B)] $H_1^2\cdot H_2^2 = 3 = L_2^3 -  \deg B_2$. 
\item[(C)] $H_1 \cdot H_2^2 = L_2^3 = 9 - \deg B_1$. 
\end{enumerate}
By $1 \leq L_2^3 \leq 5$, $L_2^3 - \deg B_2 = 3$, and $\deg B_2 >0$, 
we have two solutions: $(\deg B_2, L_2^3) \in \{ (1, 4), (2, 5)\}$. 

Assume $(\deg B_2, L_2^3) =(1, 4)$. 
Then we get $(\deg B_1, H_1^2 \cdot H_2, H_1 \cdot H_2^2) = (5, 3, 4)$. Then 
\[
(-K_X)^3 = 1 + 9 + 12 + 4 = 26, 
\]
\[
g(B_1) = 13 -32 + 20 +1 = 2, \qquad 
g(B_2) = 0.
\]

Assume $(\deg B_2, L_2^3) =(2, 5)$. 
Then we get $(\deg B_1, H_1^2 \cdot H_2, H_1 \cdot H_2^2) = (4, 3, 5)$. Then 
\[
(-K_X)^3 = 1 + 9 + 15 + 5 = 30,
\]
\[
g(B_1) = 15 -32 + 16 +1 = 0, \qquad 
g(B_2) = 0.
\]

\medskip

(5) 
Assume $(r_1, r_2) = (3, 3)$. Then $L_1^3 =L_2^3 =2$,  and 
the following hold: 
\begin{enumerate}
\item[(A)] $8 =  \deg B_1 + \deg B_2$. 
\item[(B)] $H_1^2\cdot H_2^2 = 4 = 8 -  \deg B_2$. 
\item[(C)] $H_1 \cdot H_2^2 = 4 = 8 - \deg B_1$. 
\end{enumerate}
Hence $\deg B_1 = \deg B_2 =H_1^2 \cdot H_2 = H_1\cdot H_2^2 =4$, 
\[
(-K_X)^3 = 2 + 12 + 12 + 2 = 28, \qquad 
g(B_1) = g(B_2) = 14 -27 + 12 +1 = 0. 
\]

\medskip

(6)
Assume $(r_1, r_2) = (3, 2)$. Then $L_1^3=2$, $1 \leq L_2^3 \leq 5$,  and 
the following hold: 
\begin{enumerate}
\item[(A)] $4 =  \deg B_1 + \deg B_2$. 
\item[(B)] $H_1^2\cdot H_2^2 = 4 = L_2^3 - \deg B_2$. 
\item[(C)] $H_1 \cdot H_2^2 = L_2^3 = 8 - \deg B_1$. 
\end{enumerate}
By $1 \leq L_2^3 \leq 5$, $L_2^3 - \deg B_2 = 4$, and $\deg B_2 >0$, 
we obtain $\deg B_2 = 1$ and $L_2^3 =5$. 
Hence we get 
$(\deg B_1, H_1^2 \cdot H_2, H_1 \cdot H_2^2) = (3, 4, 5)$. Then 
\[
(-K_X)^3 = 2 + 12 + 15 + 5 = 34,\]
\[
g(B_1) = 17 - 27 + 9 +1 =0, \qquad 
g(B_2) = 17 - 20 + 2 +1 =0. 
\]
\end{proof}

\section{Primitive Fano threefolds with $\rho(X)=3$}

In this section, we classify primitive Fano threefolds with $\rho(X)=3$.
The goal of this section is to prove the following theorem. 

\begin{thm}[Theorem \ref{t-pic3-CC}, Theorem \ref{t-pic3-CE}]\label{t-pic3}
Let $X$ be a primitive Fano threefold with $\rho(X)=3$. 
Then one and only one of the following holds.  
\begin{enumerate}
\item 
$(-K_X)^3=12$ and 
there exists a split double cover $f: X \to \P^1 \times \P^1 \times \P^1$ such that 
$f_*\MO_X/\MO_{\P^1 \times \P^1 \times \P^1} \simeq \MO_{\P^1 \times \P^1 \times \P^1}(1, 1, 1)^{-1}$. 
\item 
$(-K_X)^3=14$ 
and $X$ is isomorphic to a prime divisor $X'$ on 
$P := \P_{\P^1 \times \P^1}(\MO_{\P^1 \times \P^1} \oplus \MO_{\P^1 \times \P^1}(-1, -1)^{\oplus 2})$ 
such that 
$\MO_P(X') \simeq \MO_P(2) \otimes \pi^*\MO_{\P^1 \times \P^1}(2, 3)$, 
where $\pi : P= \P_{\P^1 \times \P^1}(\MO_{\P^1 \times \P^1} \oplus \MO_{\P^1 \times \P^1}(-1, -1)^{\oplus 2}) \to \P^1 \times \P^1$ denotes the natural projection. 
\item 
$(-K_X)^3=48$ and $X \simeq \P^1 \times \P^1 \times \P^1$. 
\item 
$(-K_X)^3=52$ and 
$X \simeq \P_{\P^1 \times \P^1}(\MO_{\P^1 \times \P^1} \oplus \MO_{\P^1 \times \P^1}(1, 1))$. 
\end{enumerate}
\end{thm}

By Theorem~\ref{pic-over2}, there are the following two cases, which we shall treat separately:  
\begin{itemize} 
  \item All the extremal rays are of type $C$ (Subsection \ref{ss-rho3-CC}). 
  \item There exist extremal rays $R_1$ and $R_2$ such that 
  $R_1$ is of type $C$ and $R_2$ is of type $E_1$ (Subsection \ref{ss-rho3-CE}). 
\end{itemize}

\subsection{Case $C-C-C$}\label{ss-rho3-CC}

\begin{lem}
Let $X$ be a primitive Fano threefold with $\rho(X)=3$. 
Assume that all the extremal rays are of type $C$. 
Let $F$ be a two-dimensional extremal face generated by extremal rays $R_1$ and $R_2$. 
Let $f_1 : X \to Y_1$, $f_2 : X \to Y_2$, and 
$\pi : X \to B$ be the contractions of 
$R_1$, $R_2$, and $F$, respectively. 
Then we have the following commutative diagram 
\[
\begin{tikzcd}
& X \arrow[dl, "f_1"'] \arrow[dr, "f_2"] \arrow[dd, "\pi"]\\
Y_1 \arrow[rd, "g_1"'] && Y_2 \arrow[dl, "g_2"]\\
& B
\end{tikzcd}
\]
such that the following hold. 
\begin{enumerate}
    \item $Y_1 = B \times C_1$. 
    \item $Y_2 = B \times C_2$. 
    \item $B = C_1 = C_2 = \P^1$. 
    \item $g_1: Y_1= B \times C_1 \to B$ and $g_2: Y_2= B \times C_2 \to B$ are the first projections. 
\end{enumerate}
\end{lem}

\begin{proof}
Each extremal ray $R$ is an intersection of two two-dimensional extremal faces $F$ and $F'$. 
This implies that there are exactly two non-trivial contractions $X \to B_1$ and $X \to B_2$ that factor through the contraction $X \to Y$ of $R$. 
By $Y = \P^1 \times \P^1$ {\cred (Theorem \ref{pic-over2})}, these must {\cred coincide  with} the ones induced by the projections. 
\end{proof}

\begin{rem}\label{r-rho3-CC-face}
Let $X$ be a primitive Fano threefold with $\rho(X)=3$. 
Assume that all the extremal rays are of type $C$. 
Let $F$ be an extremal face of $\NE(X)$ and 
let $h : X \to Z$ be the contraction of $F$. 
Then the following hold. 
\begin{enumerate}
\item  The following are equivalent. 
\begin{itemize}
    \item $\dim F =1$, i.e., $F$ is an extremal ray. \item $\dim Z =2$. 
    \item $Z \simeq \P^1 \times \P^1$. 
\end{itemize}
\item The following are equivalent. 
\begin{itemize}
    \item $\dim F =2$. 
    \item $\dim Z =1$. 
    \item $Z \simeq \P^1$. 
\end{itemize}
\end{enumerate}
\end{rem}

\begin{prop}\label{p-rho3-CC-ext=3}
Let $X$ be a primitive Fano threefold with $\rho(X)=3$. 
Assume that all the extremal rays are of type $C$. 
Then the number of the extremal rays of $X$ is three. 
\end{prop}

\begin{proof}
Suppose that $X$ has at least four {\cred extremal} rays. 
Then we can find two two-dimensional extremal faces $F$ and $F'$ such that $F \cap F' = \{0 \}$. 
Let $\pi : X \to B$ and $\pi' : X \to B'$ be the contractions 
of $F$ and $F'$, respectively. 
By Remark \ref{r-rho3-CC-face}, we get $\dim B = \dim B' =1$. 
Then there exists a curve $C$ on $X$ such that $(\pi \times \pi')(C)$ is a point for $\pi \times \pi' : X \to B \times B'$. 
This implies that $\pi(C)$ and $\pi'(C)$ are points, i.e., 
$[C] \in F \cap F' = \{0\}$, which is absurd. 
\end{proof}

\begin{nota}\label{n-pic3-CC}
Let $X$ be a primitive Fano threefold with $\rho(X)=3$. 
Assume that all the extremal rays are of type $C$. 
Then there exist exactly three two-dimensional extremal faces $F_1, F_2, F_3$ (Proposition \ref{p-rho3-CC-ext=3}). 
We have also exactly three extremal rays $R_{12}, R_{23}, R_{13}$, 
which are given as follows: 
\[
R_{12} := F_1 \cap F_2, \qquad 
R_{23} := F_2 \cap F_3, \qquad 
R_{13} := F_1 \cap F_3. 
\]
Corresponding to these extremal faces, we have the following contraction morphisms: 
\begin{itemize}
\item $\pi^X_i : X \to B_i := \P^1$ for all $1 \leq i \leq 3$. 
\item $f_{ij}:X \to Y_{ij} :=B_i \times B_j = \P^1 \times \P^1$ 
for all $1 \leq i < j \leq 3$. 
\end{itemize}
Note that $f_{ij}$ and $\pi_i$ are compatible, i.e.,
\[
\pi^X_i : X \xrightarrow{f_{ij}} Y_{ij} = B_i \times B_j \xrightarrow{{\rm pr}_1} B_i, \qquad 
\pi^X_j : X \xrightarrow{f_{ij}} Y_{ij} = B_i \times B_j \xrightarrow{{\rm pr}_2} B_j. 
\]
Set 
\[
f := \pi_1^X \times \pi_2^X \times \pi_3^X : X \to Z := B_1 \times B_2 \times B_3 = \P^1 \times \P^1 \times \P^1. 
\]
Let $\pi_i^Z : Z \to B_i =\P^1$ be the $i$-th projection. 
For each $i \in \{1, 2, 3\}$, we set 
\begin{itemize}
\item $H_i:=(\pi^X_i)^*\MO_{\mathbb P^1}(1)$, and
\item $H^Z_i:=(\pi^Z_i)^*\MO_{\mathbb P^1}(1)$. 
\end{itemize}
For $1 \leq i < j \leq 3$, let $\ell_{ij}$ be an extremal rational curve of the extremral ray $R_{ij}$. 
Note that we have 
$H_i \cdot \ell_{ij} = H_j \cdot \ell_{ij} =0$. 
\end{nota}



\begin{lem}\label{l-rho3-CC-finite}
We use Notation \ref{n-pic3-CC}. 
Then   $f : X \to Z =\mathbb P^1 \times \mathbb P^1 \times \mathbb P^1$ is a finite surjective morphism. 
\end{lem}

\begin{proof}
It suffices to show that $f$ is a finite morphism. 
  Suppose that $f$ is not a finite morphism. 
Then there exists a curve $C$ on $X$ such that $f(C)$ is a point. 
By $f = \pi^X_1 \times \pi^X_2 \times \pi^X_3$, 
all of $\pi^X_1(C), \pi^X_2(C)$, and $\pi^X_3(C)$ are points. 
Then $[C] \in F_1 \cap F_2 \cap F_3 = \{0\}$, which is a contradiction. 
\end{proof}



\begin{thm}\label{t-pic3-CC}
Let $X$ be a primitive Fano threefold with $\rho(X)=3$. 
Assume that all the extremal rays of $\NE(X)$ are of type $C$. 
Let $R_1$ and $R_2$ be distinct extremal rays of $\NE(X)$. 
Then one and only one of the following holds.  
\begin{enumerate}
\item 
There exists a split double cover $f: X \to \P^1 \times \P^1 \times \P^1$ such that 
$f_*\MO_X/\MO_{\P^1 \times \P^1 \times \P^1} \simeq \MO_{\P^1 \times \P^1 \times \P^1}(1, 1, 1)^{-1}$. 
Furthermore, $\MO_X(-K_X) \simeq f^*\MO_{\P^1 \times \P^1 \times \P^1}(1, 1, 1)$, $(-K_X)^3=12$,  there are exactly three extremal rays of $\NE(X)$, all the extremal rays are of type $C_1$, and 
$\Delta_{\varphi}$ is of bidegree $(4, 4)$ for the contraction $\varphi : X \to \P^1 \times \P^1$ of an arbitrary extremal ray. 
\item 
$X \simeq \P^1 \times \P^1 \times \P^1$. 
Furthermore, $(-K_X)^3=48$, there are exactly three extremal rays $R_1, R_2, R_3$ of $\NE(X)$, 
and all the extremal rays are of type $C_2$. 
\end{enumerate}
\end{thm}

\begin{proof}
We use Notation \ref{n-pic3-CC}. 
Set $d := \deg f$. 
For each $i \in \{1, 2, 3\}$, we have morphisms: 
\[
\pi^X_i : X \xrightarrow{f} Z = \mathbb P^1 \times \mathbb P^1 \times \mathbb P^1 \xrightarrow{\pi_i^Z} \mathbb P^1.  
\]
By  
\[
H_1^Z \cdot H_2^Z \cdot H_3^Z = (\pi_1^Z)^*\MO_{\mathbb P^1}(1) \cdot (\pi_2^Z)^*\MO_{\mathbb P^1}(1)  \cdot (\pi_3^Z)^*\MO_{\mathbb P^1}(1) =1, 
\]
we have 
\[
d = \deg f = 
(f^*H_1^Z) \cdot (f^*H_2^Z) \cdot (f^*H_3^Z) = H_1 \cdot H_2 \cdot H_3. 
\]
Since $H_1 \cdot H_2$ is a fibre of $f_{12} : X \to Y_{12} = 
B_1 \times B_2 = \mathbb P^1 \times \P^1$, 
we obtain $H_1 \cdot H_2 \equiv \frac{2}{\mu_{12}}\ell_{12}$ 
{\cred (Corollary \ref{fib_of_C})}, 
where $\mu_{12}$ denotes the length of $R_{12}$. 
Therefore, we get 
\[
d = H_1 \cdot H_2 \cdot H_3 
=\frac{2}{\mu_{12}} H_3 \cdot \ell_{12}, 
\]
which implies 
\[
-K_X \cdot \ell_{12} = \mu_{12} = \frac{2}{d}H_3\cdot \ell_{12}. 
\]
By symmetry, we have that 
$ -d K_X \cdot \ell_{12} = 2 H_3 \cdot \ell_{12}, 
-d K_X \cdot \ell_{23} = 2 H_1 \cdot \ell_{23},  -d K_X \cdot \ell_{13} = 2 H_2 \cdot \ell_{13}.$ 
By the exact sequence 
\[
0 \to \Pic\,Y_{12} \to \Pic\,X \xrightarrow{\cdot \ell_{12}} \Z, 
\]
$H_1, H_2, H_3$ form a $\Q$-linear basis of $(\Pic\,X) \otimes_{\Z} \Q$. 
Hence we obtain  
\[
-d  K_X \sim 2H_1 + 2H_2 + 2H_3. 
\]
For the positive integer $g \in \Z_{>0}$ satisying $(-K_X)^3 = 2g-2$, we get 
\[
d^3(2g-2) =d^3 (-K_X)^3 = 8(H_1+H_2+H_3)^3 = 48 \deg f = 48d. 
\]
Hence $d^2 (g-1)= 24 = 2^3 \cdot 3$, which implies $d \in \{1, 2\}$.  

Assume $d=1$. 
Then $f : X \to Z = \P^1 \times \P^1 \times \P^1$ is a finite birational morphism 
(Lemma \ref{l-rho3-CC-finite}). Since $Z$ is normal, $f$ is an isomorphism. We have $(-K_X)^3 = \MO_{\P^1 \times \P^1 \times \P^1}(2, 2, 2)^3 = 48$. 
Furthermore, each projection $X \to \P^1 \times \P^1$ is clearly of type $C_2$. 
Hence (2) holds.

Assume $d=2$. 
Then $-K_X \sim H_1 + H_2 + H_3$ and $f : X \to Z = \P^1 \times \P^1 \times \P^1$ is a double cover (Lemma \ref{l-rho3-CC-finite}). 
This double cover is split (Lemma \ref{l-split-criterion}), 
because $\P^1 \times \P^1 \times \P^1$ is $F$-split and 
$H^1(\P^1 \times \P^1 \times \P^1, \MO_{\P^1 \times \P^1 \times \P^1}(E))=0$ 
for {\cred every}  effective divisor $E$ on $\P^1 \times \P^1 \times \P^1$. 
 We have $(-K_X)^3 = \deg f \cdot \MO_{\P^1 \times \P^1 \times \P^1}(1, 1, 1)^3 = 12$.

Let us show that $f_{12} : X \to B_1 \times B_2$ is of type $C_1$ 
{\cred and $\Delta_{f_{12}}$ is of bidegree $(4, 4)$}. 
We can write 
\[
\MO_{\P^1 \times \P^1}(\Delta_{f_{12}}) = \MO_{\P^1 \times \P^1}(a, b)
\]
for some $a, b \in \Z$. 
We obtain 
 \[
4 = 2d = 2 H_1 \cdot H_2 \cdot H_3 =   (-K_X)^2\cdot H_1 
\overset{(\star)}{=} -4 K_{\mathbb{P}^1\times \mathbb{P}^1}\cdot \MO_{\P^1 \times \P^1}(1, 0)- \Delta_{f_{12}} \cdot  \MO_{\P^1 \times \P^1}(1, 0)
\]
\[
= \MO_{\P^1 \times \P^1}(8-a, 8-b)\cdot \MO_{\P^1 \times \P^1}(1, 0)
=8-b.
\]
where ($\star$) follows from Proposition \ref{conic_disc-num}. 
Therefore, $b = 4$. 
Similarly, $a=4$. Hence $f_{12}$ is of type $C_1$ {\cred and $\Delta_{f_{12}}$ is of bidegree $(4, 4)$}. 
By symmetry, $f_{ij} : X \to B_i \times B_j$ is of type $C_1$ 
{\cred  and $\Delta_{f_{ij}}$ is of bidegree $(4, 4)$} for all $1 \leq i < j \leq 3$. 
Thus (1) holds. 
\end{proof}

\begin{rem}
By the above proof, 
the split double cover $f : X \to \P^1 \times \P^1 \times \P^1$ in  Theorem \ref{t-pic3-CC}(1) 
can be chosen to be $f$ as in Notation \ref{n-pic3-CC}. 
\end{rem}

\subsection{Case $C-E$}\label{ss-rho3-CE}


\begin{nota}\label{n-pic3-CE}
Let $X$ be a primitive Fano threefold with $\rho(X)=3$. 
Let $R_1$ and $R_2$ be two distinct extremal rays. 
Assume that 
\begin{itemize}
\item $R_1$ is of type $C$, and 
\item $R_2$ is of type $E_1$. 
\end{itemize} 
For each $i \in \{1, 2\}$, 
let 
\[
f_i : X \to Y_i 
\] 
be the contraction of $R_i$, where $Y_1 =\mathbb P^1 \times \mathbb P^1$ {\cred (Theorem \ref{pic-over2})}. 
Set $D := \Ex(f_2)$, which is a prime divisor on $X$ such that  $D \simeq \mathbb{P}^1\times \mathbb{P}^1$ (Theorem \ref{dimy=3}). 
{\cred For each $i \in \{1, 2\}$,}  
let $\ell_i$ be an extremal rational curve with $R_i = \R_{\geq 0}[\ell_i]$. 
Set $\mu_i := -K_X \cdot \ell_i$, which is the length of $R_i$.  
\end{nota}

\begin{lem}\label{rho3-f_D}
We use Notation \ref{n-pic3-CE}. 
For the induced morphism   
\[
f_1|_D\colon D\simeq \mathbb{P}^1\times \mathbb{P}^1 \to \mathbb{P}^1 \times \mathbb{P}^1=Y_1,
\]
the following hold. 
  \begin{enumerate}
    \item If $R_1$ is of type $C_1$, then $f_1|_D$ is a double cover.
    \item If $R_2$ is of type $C_2$, then $f_1|_D$ is an isomorphism. 
  \end{enumerate}
\end{lem}

\begin{proof}
By Lemma \ref{l-finite-morph}(2), $f_1|_D : D \to Y_1$ is a finite surjective morphism. 
Fix a fibre $\zeta$ of $f_1 : X \to Y_1$ over a closed point of $Y_1$. 
By $\zeta \equiv \frac{2}{\mu_1}\ell_1$ (Corollary~\ref{fib_of_C}), 
it holds that 
  \[
\deg(f_1|_D) = D \cdot \zeta = \frac{2}{\mu_1} D\cdot \ell_1.
\]
Recall that we have the following exact sequence (Theorem \ref{t-ex-surje}):
  \[
  0 \longrightarrow \Pic\,(\mathbb{P}^1\times \mathbb{P}^1) 
  \stackrel{f_1^*}{\to} \Pic\,X \stackrel{ \cdot \ell_1}{\to} \mathbb{Z} \longrightarrow 0.
  \]
Fix a Cartier divisor $E$ on $X$ with $E \cdot \ell_1 =1$, 
whose existence is guaranteed by the surjectivity of the last map $E \mapsto  E\cdot \ell_1$. 
Then there exist $a \in \Z$ and a Cartier divisor $L$ on $Y_1 = \mathbb P^1 \times \mathbb P^1$ such that 
\[
D \sim f_1^*L + aE. 
\]
By $a = D \cdot \ell_1$, we have $a \in \Z_{>0}$. 
 Since $L$ is a divisor on $\mathbb{P}^1\times \mathbb{P}^1$, we have $L^3=0$ and 
$\MO_{\mathbb{P}^1\times \mathbb{P}^1}(L)\simeq \MO_{\mathbb{P}^1\times \mathbb{P}^1}(b,c)$ for some $b,c\in \mathbb{Z}$. 
We have $L^2 \equiv 2bcQ$ for a closed point $Q$ on $\mathbb{P}^1\times \mathbb{P}^1$, which implies $f_1^*L^2 \equiv 2bc \zeta   \equiv \frac{4bc}{\mu_1}\ell_1$ 
and $f_1^*L^2 \cdot E  = \frac{4bc}{\mu_1}$. 
Therefore, we obtain 
  \begin{align*}
    D^3 &= (f^*_1L+aE)^3 \\
    &= f_1^*L^3+3 f_1^*L^2\cdot (aE) + 3f^*_1L\cdot (a^2E^2) + a^3E^3 \\
    &= 3a\cdot \frac{4bc}{\mu_1} + 3a^2(f^*_1L\cdot E^2) + a^3E^3.
  \end{align*}
  On the other hand, for $C:=f_2(D) (\simeq \P^1)$, 
it holds that 
\[
D^3 = \deg_C(\mathcal N^*_{C/Y_2}) = \deg_{\mathbb{P}^1}(\MO_{\mathbb{P}^1}(1)\oplus \MO_{\mathbb{P}^1}(1)) =2, 
\]
where the first and second equalities follow from Proposition \ref{standard_3}(3) and 
Lemma \ref{l-E1-prim}(2), respectively. 
We then obtain 
  \begin{equation}\label{D^3-equation}
    2=D^3 = a\left(\frac{12bc}{\mu_1}+3a(f_1^*L\cdot E^2)+a^2E^3\right).
  \end{equation}
By $a \in \Z_{>0}, b, c \in \Z$, and $\mu_1 \in \{1, 2\}$, it holds that $a=1$ or $a=2$. 
If $a=2$, then the right hand side of (\ref{D^3-equation}) would be contained in $4\Z$, 
which is absurd. 
Then we obtain $a=1$, which implies  $D \cdot \ell_1 =1$ and  $\deg(f_1|_D)=\frac{2}{\mu_1}$.
\medskip

(1) 
Assume that $R_1$ is of type $C_1$. 
Then we have $\deg(f_1|_D)=2$, i.e., $f_1|_D$ is a double cover. 
Hence (1) holds. 
\medskip

(2) 
Assume that $R_1$ is of type $C_2$. 
Then we have $\deg(f_1|_D)=1$, i.e.,  $f_1|_D$ is a finite birational morphism. 
Since $f_1|_D : D=\P^1 \times \P^1 \to Y_1 = \mathbb P^1 \times \mathbb P^1$ is a finite birational morphism of normal varieties, 
$f_1|_D$ is an isomorphism. 
Thus (2) holds. 
\end{proof}

\subsubsection{Case $C_2-E_1$}


\begin{prop}\label{p-pic3-C2E1}
Let $X$ be a primitive Fano threefold with $\rho(X)=3$. 
Assume that there exist two extremal rays $R_1$ and $R_2$ such that 
$R_1$ is of type $C_2$ and $R_2$ is of type $E_1$. 
Then 
\[
X \simeq \mathbb P_{\mathbb P^1 \times \mathbb P^1}(\MO_{\P^1 \times \P^1} \oplus \MO_{\P^1 \times \P^1}(1, 1))
\]
and $(-K_X)^3 = 52$. 
\end{prop}

\begin{proof}
We use Notation \ref{n-pic3-CE}. 
Recall that we have $X \simeq \mathbb P_{\P^1 \times \P^1}(\mathcal E)$ and $\mathcal E = (f_1)_*\MO_X(1)$. 
Set $\MO_D(D) := \MO_X(D)|_D$.

\setcounter{step}{0}

\begin{step}\label{s1-pic3-C2E1}
There is the following exact sequence: 
\[
0 \to (f_1)_*\MO_X\to (f_1)_*\MO_X(D) \to (f_1)_*\MO_D(D) \to 0.
\]
\end{step}

\begin{proof}[Proof of Step \ref{s1-pic3-C2E1}]
We have an exact sequence 
\[0 \to \MO_X \to \MO_X(D) \to \MO_D(D) \to 0,\]
which induces another exact sequence 
\[ 0 \to (f_{1})_*\MO_X\to (f_1)_*\MO_X(D) \to (f_1)_*\MO_D(D) \to R^1(f_1)_*\MO_X. \]
By Proposition \ref{p-cont-ex2}, $R^1(f_1)_*\MO_X = 0$. 
This completes the proof of Step \ref{s1-pic3-C2E1}. 
\end{proof}

\begin{step}\label{s2-pic3-C2E1}
The following hold. 
\begin{enumerate}
\item $(f_1)_*\MO_X(D) \simeq \mathcal E \otimes \mathcal L$ for some $\mathcal L \in \Pic (\P^1 \times \P^1)$. 
\item $(f_1)_*\MO_D(D) \simeq \MO_{\P^1 \times \P^1}(-1, -1)$. 
\end{enumerate}
\end{step}

\begin{proof}[Proof of Step \ref{s2-pic3-C2E1}]
Let us show (1). 
Since $f_1 : X \to Y_1 = \P^1 \times \P^1$ is a $\mathbb P^1$-bundle, it holds that 
\[
\MO_X(D) \simeq \MO_X(n) \otimes f_1^*\mathcal L
\]
for some $n \in \Z$ and $\mathcal L \in \Pic (\P^1 \times \P^1)$. 
Fix a fibre $\zeta$ of $f_1 : X \to Y_1 = \P^1 \times \P^1$. 
Since $f_1|_D : D \to Y_1$ is an isomorphism (Lemma \ref{rho3-f_D}), we have $D \cdot \zeta =1$. 
By $\MO_X(1) \cdot \zeta =1$, we obtain $n=1$. 
It holds that 
\[
(f_1)_*\MO_X(D) \simeq (f_1)_*(\MO_X(1) \otimes f_1^*\mathcal L) \simeq 
 ((f_1)_*\MO_X(1)) \otimes \mathcal L \simeq \mathcal E \otimes \mathcal L. 
\]
Thus (1) holds. 

Let us show (2). 
Since $f_1|_D : D \to Y_1 = \P^1 \times \P^1$ is an isomorphism, 
Lemma \ref{l-E1-prim}(4) implies $(f_1)_*\MO_D(D) \simeq \MO_{\P^1 \times \P^1}(-1, -1)$. 
This completes the proof of Step \ref{s2-pic3-C2E1}. 
\end{proof}

\begin{step}\label{s3-pic3-C2E1}
$X \simeq \mathbb P_{\mathbb P^1 \times \mathbb P^1}(\MO_{\P^1 \times \P^1} \oplus \MO_{\P^1 \times \P^1}(1, 1))$. 
\end{step}

\begin{proof}[Proof of Step \ref{s3-pic3-C2E1}]
By Step \ref{s1-pic3-C2E1}and Step \ref{s2-pic3-C2E1},  we have the following exact sequence for some $\mathcal L \in \Pic\,(\P^1 \times \P^1)$:
\[
0\to \MO_{\mathbb{P}^1\times \mathbb{P}^1}\to \mathcal{E} \otimes \mathcal{L}\to \MO_{\mathbb{P}^1\times \mathbb{P}^1}(-1,-1) \to 0.
\]
It holds that 
\[
\Ext^1_{\mathbb{P}^1\times \mathbb{P}^1}( \MO_{\mathbb{P}^1 \times \P^1}(-1,-1), \MO_{\mathbb{P}^1\times \mathbb{P}^1}) 
\simeq H^1(\P^1 \times \P^1, \MO_{\mathbb{P}^1 \times \P^1}(1,1))=0, 
\]
which implies 
\[
\mathcal{E} \otimes \mathcal{L} \simeq  \MO_{\mathbb{P}^1\times \mathbb{P}^1} \oplus \MO_{\mathbb{P}^1\times \mathbb{P}^1}(-1,-1).
\]
Therefore, 
\[
X \simeq \P_{\P^1 \times \P^1}(\mathcal E) \simeq 
 \P_{\P^1 \times \P^1}(\mathcal E  \otimes \mathcal{L} \otimes \MO_{\mathbb{P}^1\times \mathbb{P}^1}(1, 1))
\simeq \mathbb P_{\mathbb P^1 \times \mathbb P^1}(\MO_{\P^1 \times \P^1} \oplus \MO_{\P^1 \times \P^1}(1, 1)).
\]
This completes the proof of Step \ref{s3-pic3-C2E1}. 
\end{proof}

\begin{step}\label{s4-pic3-C2E1}
$(-K_X)^3 = 52$. 
\end{step}

\begin{proof}[Proof of Step \ref{s4-pic3-C2E1}]
By Proposition~\ref{standard_2}{\cred (3)}, we have 
\[(-K_X)^3=2c_1(\MO_{\P^1 \times \P^1}\oplus \MO_{\P^1 \times \P^1}(1,1))^2-8c_2(\MO_{\P^1 \times \P^1}\oplus \MO_{\P^1 \times \P^1}(1,1))+6K_{\mathbb{P}^1\times \mathbb{P}^1}^2.\]
By \cite[Appendix A, \S 3, C3 and C5]{Har77}, we obtain 
\begin{equation*}
  \begin{aligned}
    c_1(\MO_{\mathbb{P}^1\times \mathbb{P}^1}\oplus \MO_{\mathbb{P}^1\times \mathbb{P}^1}(1,1))&=c_1(\MO_{\mathbb{P}^1\times \mathbb{P}^1})+c_1(\MO_{\mathbb{P}^1\times \mathbb{P}^1}(1,1))=c_1(\MO_{\mathbb{P}^1\times \mathbb{P}^1}(1,1)),\\
    c_2(\MO_{\mathbb{P}^1\times \mathbb{P}^1}\oplus \MO_{\mathbb{P}^1\times \mathbb{P}^1}(1,1))&=c_1(\MO_{\mathbb{P}^1\times \mathbb{P}^1})\cdot c_1(\MO_{\mathbb{P}^1\times \mathbb{P}^1}(1,1))=0.
  \end{aligned}
\end{equation*}
Hence it holds that $(-K_X)^3 = 2\cdot 2-8\cdot 0 +6\cdot 8 = 52.$ 
This completes the proof of Step \ref{s4-pic3-C2E1}. 
\end{proof}
Step \ref{s3-pic3-C2E1} and Step \ref{s4-pic3-C2E1} complete the proof of Proposition \ref{p-pic3-C2E1}
\end{proof}

\subsubsection{Case $C_1-E_1$}

\begin{nota}\label{n-pic3-C2E}
We use Notation \ref{n-pic3-CE}. 
Assume that $R_1$ is of type $C_1$. 
In particular, 
$X$ is a primitive Fano threefold with $\rho(X)=3$, 
$R_1$ is of type $C_1$, and 
$R_2$ is of type $E_1$.  
By Lemma~\ref{conic-embedding}, we have 
a commutative diagram 
\begin{center}
  \begin{tikzcd}
    X \arrow[r, hook, "\iota^{\prime}"] \arrow[rd, "f_1"']&P^{\prime}:=\mathbb{P}(\mathcal{E}^{\prime}) \arrow[d,"\pi^{\prime}"] \\
    & Y_1 =\mathbb{P}^1\times \mathbb{P}^1
  \end{tikzcd}
\end{center}
where $\mathcal E' := (f_{1})_*\MO_X(-K_X)$, $\iota'$ is a closed immersion, 
{\cred and $\pi'$ is the natural projection}. 
We identify $X$ with the smooth prime divisor $\iota'(X)$ on $P'$. 
\end{nota}

\begin{lem}\label{l-D-to-Y_1-12}
We use Notation \ref{n-pic3-C2E}. 
For the induced morphism 
\[
f_1|_D : D =\P^1 \times \P^1 \to Y_1 = \P^1 \times \P^1,
\]
we set $\mathcal L := ((f_1|_D)_*\MO_D/\MO_{Y_1})^{-1}$, which is an invertible sheaf (cf. Remark \ref{r-L-inv}, Lemma \ref{rho3-f_D}(1)).  
Then, after possibly permuting the direct product factors of $Y_1 = \P^1 \times \P^1$, the following hold. 
\begin{enumerate}
\item $\mathcal L \simeq \MO_{\P^1 \times \P^1}(0, 1)$. 
\item $f_1|_D$ is a split double cover. 
\item There exists a double cover $h : \mathbb P^1 \to \mathbb P^1$ such that $f_1|_D = {\rm id} \times h : D =\P^1 \times \P^1 \to Y_1 = \P^1 \times \P^1$. 
\item 
  $(f_1|_D)_*\MO_D(0,-1) \simeq \MO_{\mathbb{P}^1\times \mathbb{P}^1}(0,-1) \oplus \MO_{\mathbb{P}^1\times \mathbb{P}^1}(0,-1)$.
\end{enumerate}
\end{lem}

\begin{proof}
Since $f_1|_D\colon D\simeq \mathbb{P}^1\times \mathbb{P}^1\to \mathbb{P}^1\times \mathbb{P}^1$ is a double cover (Lemma \ref{rho3-f_D}(1)), 
Lemma \ref{l-dc-omega} implies 
\[
\omega_D\simeq (f_1|_D)^*(\omega_{\mathbb{P}^1\times \mathbb{P}^1}\otimes \mathcal{L}).
\]
We can write 
$\mathcal{L} \simeq \MO_{\mathbb{P}^1\times \mathbb{P}^1}(a,b)$ 
for some $a,b\in \mathbb{Z}$. 
We then get 
\begin{equation}\label{D_P1P1_K}
  \MO_D(-2,-2) \simeq (f_1|_D)^*(\MO_{\mathbb{P}^1\times \mathbb{P}^1}(a-2,b-2)). 
\end{equation}
Hence $a <2$ and $b<2$. Furthermore, we get 
\begin{align*}
  8 &= (c_1(\MO_D(-2,-2)))^2 \\
  &= \deg(f_1|_D)\cdot c_1(\MO_{\mathbb{P}^1\times \mathbb{P}^1}(a-2,b-2))^2 \\
  &= 4(a-2)(b-2), 
\end{align*}
i.e., $(a-2)(b-2)=2$. 
Therefore, after possibly permuting the direct product factors of 
$Y_1 = \P^1 \times \P^1$, we obtain $(a, b)=(0, 1)$ and 
$\mathcal{L} \simeq \MO_{\mathbb{P}^1\times \mathbb{P}^1}(0,1)$. 
Thus (1) holds. 
Lemma \ref{l-split-criterion}{\cred (2)} and Lemma \ref{l-trivial-double}  imply (2) and (3), respectively. 

Let us show (4). 
 It is easy to see that $h_*\MO_{\mathbb{P}^1}(1)$ is a locally free sheaf of rank 2 on $\mathbb{P}^1$.
  Hence we can write 
$h_*\MO_{\mathbb{P}^1}(1)\simeq \MO_{\mathbb{P}^1}(c)\oplus \MO_{\mathbb{P}^1}(d)$   for some $c,d\in \mathbb{Z}$.  
We obtain 
\begin{eqnarray*}
h^0(\MO_{\mathbb{P}^1}(c)) + h^0(\MO_{\mathbb{P}^1}(d))
&=& h^0(\P^1, h_*\MO_{\mathbb{P}^1}(1)) = h^0(\P^1, \MO_{\mathbb{P}^1}(1)) =2\\
h^0(\MO_{\mathbb{P}^1}(c-1)) + h^0(\MO_{\mathbb{P}^1}(d-1))
&=& h^0(\P^1, \MO_{\P^1}(-1) \otimes h_*\MO_{\mathbb{P}^1}(1))
= h^0(\P^1, \MO_{\mathbb{P}^1}(-1)) =0. 
\end{eqnarray*}
Therefore, we get $c=d=0$ and 
\[
h_*\MO_{\mathbb{P}^1}(-1) \simeq \MO_{\P^1}(-1) \otimes h_*\MO_{\mathbb{P}^1}(1) \simeq \MO_{\P^1}(-1) \oplus \MO_{\P^1}(-1) . 
\]
This, together with the flat base change theorem, implies the following: 
\begin{align*}
(f_1|_D)_*\MO_D(0,-1) & \simeq ({\rm id} \times h)_* \MO_D(0,-1)\\
&\simeq ({\rm id} \times h)_*  ({\rm pr}_2)^*\MO_{\P^1}(-1)\\
&\simeq({\rm pr}_2)^* h_*  \MO_{\P^1}(-1)\\
&\simeq({\rm pr}_2)^*(\MO_{\P^1}(-1)\oplus \MO_{\P^1}(-1))\\
&\simeq \MO_{\P^1\times \P^1}(0,-1) \oplus  \MO_{\P^1\times \P^1}(0,-1). 
  \end{align*}
  Thus (4) holds. 
  \qedhere

\end{proof}

\begin{lem}\label{l-pic3-C2E-ext}
We use Notation \ref{n-pic3-C2E}. 
Then, after possibly permuting the direct product factors of $Y_1 = \P^1 \times \P^1$, 
the following holds: 
\[
(f_1)_*\MO_X(-K_X)\simeq \MO_{\mathbb{P}^1\times \mathbb{P}^1}(2,1)\oplus \MO_{\mathbb{P}^1\times \mathbb{P}^1}(1,0)^{\oplus 2}. 
\]
\end{lem}

\begin{proof}

\setcounter{step}{0}
\begin{step}\label{s1-pic3-C2E-ext}
The following hold. 
\begin{enumerate}
\item $(-K_X-D) \cdot \zeta =0$ for a fibre $\zeta$ of $f_1: X \to Y_1$. 
\item 
There exists the following exact sequence:  
\[ 
0 \longrightarrow (f_1)_*\MO_X(-K_X-D) \longrightarrow (f_1)_*\MO_X(-K_X) \longrightarrow (f_1)_*\MO_D(-K_X) \longrightarrow 0. 
\]
\end{enumerate}
\end{step}

\begin{proof}[Proof of Step \ref{s1-pic3-C2E-ext}]
The assertion (1) follows from $-K_X \cdot \zeta =2$ and $D \cdot \zeta =2$ (Lemma \ref{rho3-f_D}).  
Let us show (2).  
We have an exact sequence
\[0 \longrightarrow \MO_X(-K_X-D) \longrightarrow \MO_X(-K_X) \longrightarrow \MO_D(-K_X) \longrightarrow 0,\]
which induces another exact sequence
\begin{align*}
  0 \to (f_1)_*\MO_X(-K_X-D) \to (f_1)_*\MO_X(-K_X) \to (f_1)_*\MO_D(-K_X) \to 
  R^1(f_1)_*\MO_X(-K_X-D).
\end{align*}
By (1) and \cite[Theorem 0.5]{Tan15}, 
we obtain $R^1(f_1)_{*}\MO_X(-K_X-D)=0$. 
Thus (2) holds. 
This completes the proof of Step  \ref{s1-pic3-C2E-ext}. 
\end{proof}

\begin{step}\label{s2-pic3-C2E-ext}
$(f_1)_*\MO_X(-K_X-D)\simeq \MO_{\mathbb{P}^1\times \mathbb{P}^1}(2,1)$.
\end{step}

\begin{proof}[Proof of Step \ref{s2-pic3-C2E-ext}]
By Step \ref{s1-pic3-C2E-ext}(1), 
we have $(-K_X-D)\cdot \zeta =0$. 
Therefore, we can write 
\[
\MO_X(-K_X-D) \simeq f_1^*\MO_{\mathbb{P}^1\times \mathbb{P}^1}(s,t)
\]
for some $s,t\in \mathbb{Z}$ (Theorem \ref{t-ex-surje}). 
We obtain 
\begin{align*}
  \MO_D(2,2) &\simeq \MO_D(-K_D) \\
  &\simeq \MO_X(-K_X-D)|_D \\
  &\simeq (f_1|_D)^* \MO_{\mathbb{P}^1\times \mathbb{P}^1}(s,t) \\
  &\simeq ({\rm id}\times h)^*\MO_{\mathbb{P}^1\times \mathbb{P}^1}(s,t) \\
  &\simeq \MO_D(s,2t), 
\end{align*}
where the fourth and fifth isomorphisms follow from Lemma \ref{l-D-to-Y_1-12}. 
Hence we get $s=2$ and $t=1$. 
Then it holds that  
\[
(f_1)_*\MO_X(-K_X-D) \simeq (f_1)_*f_1^*\MO_{\mathbb{P}^1\times \mathbb{P}^1}(s,t) \simeq \MO_{\mathbb{P}^1\times \mathbb{P}^1}(s,t) 
\simeq \MO_{\mathbb{P}^1\times \mathbb{P}^1}(2,1). 
\]
This completes the proof of Step \ref{s2-pic3-C2E-ext}. 
\end{proof}

\begin{step}\label{s3-pic3-C2E-ext}
$(f_1)_*\MO_D(-K_X) \simeq \MO_{\P^1 \times \P^1}(1, 0)^{\oplus 2}$. 
\end{step}

\begin{proof}[Proof of Step \ref{s3-pic3-C2E-ext}]
Recall that we have 
\[
(f_1)_*\MO_D(-K_X) =(f_1|_D)_*(\MO_X(-K_X)|_D).
\]  
By $\MO_X(-K_X-D)|_D \simeq \MO_D(2,2)$ and $\MO_X(D)|_D \simeq\MO_D(-1,-1)$ (Theorem \ref{dimy=3}), it holds that 
\[
(f_1|_D)_*(\MO_X(-K_X)|_D) \simeq (f_1|_D)_*\MO_D(1,1).
\] 
By 
Lemma \ref{l-D-to-Y_1-12}, we obtain 
\begin{align*}
(f_1)_*\MO_D(-K_X) &\simeq (f_1|_D)_*\MO_D(1,1)  \\
  &\simeq ({\rm id}\times h)_*(\MO_D(1,2)\otimes \MO_D(0,-1)) \\
  &\simeq ({\rm id}\times h)_*((({\rm id}\times h)^*\MO_{\mathbb{P}^1\times \mathbb{P}^1}(1,1)) \otimes \MO_D(0,-1)) \\
  &\simeq \MO_{\mathbb{P}^1\times \mathbb{P}^1}(1,1) \otimes ({\rm id}\times h)_*\MO_D(0,-1)\\
&\simeq \MO_{\mathbb{P}^1\times \mathbb{P}^1}(1,0)^{\oplus 2}. 
\end{align*}
This completes the proof of Step  \ref{s3-pic3-C2E-ext}. 
\end{proof}

\begin{step}\label{s4-pic3-C2E-ext}
It holds that 
\[
(f_1)_*\MO_X(-K_X) \simeq \MO_{\mathbb{P}^1\times \mathbb{P}^1}(2,1)\oplus \MO_{\mathbb{P}^1\times \mathbb{P}^1}(1,0)^{\oplus 2}. 
\]
\end{step}

\begin{proof}[Proof of Step \ref{s4-pic3-C2E-ext}]
By Step \ref{s1-pic3-C2E-ext}, we have the following exact sequence:  
\[
0 \longrightarrow (f_1)_*\MO_X(-K_X-D) \longrightarrow (f_1)_*\MO_X(-K_X) \longrightarrow (f_1)_*\MO_D(-K_X) \longrightarrow 0.
\]
By Step \ref{s2-pic3-C2E-ext} and Step \ref{s3-pic3-C2E-ext}, 
this exact sequence can be written as follows: 
\[ 
0 \longrightarrow \MO_{\mathbb{P}^1\times \mathbb{P}^1}(2,1) \longrightarrow (f_1)_*\MO_X(-K_X) \longrightarrow \MO_{\P^1 \times \P^1}(1, 0)^{\oplus 2} \longrightarrow 0. 
\]
Hence it suffices to show that this exact sequence splits, 
which follows from the following computation: 
\[
\mathrm{Ext}^1(\MO_{\P^1\times \P^1}(1,0)^{\oplus 2}, 
\MO_{\P^1 \times \P^1}(2,1))
\simeq 
H^1(\P^1 \times \P^1, \MO_{ \P^1 \times \P^1}(1,1))^{\oplus 2}=0. 
\]
This completes the proof of Step  \ref{s4-pic3-C2E-ext}. 
\end{proof}
Step \ref{s4-pic3-C2E-ext} completes the proof of Lemma \ref{l-pic3-C2E-ext}. 
\end{proof}

\begin{lem}\label{l-pic3-C2E-X}
We use Notation \ref{n-pic3-C2E}. Then it holds that 
\[
\MO_{P^\prime}(X)\simeq \MO_{P^\prime}(2)\otimes \pi^{\prime *}\MO_{\mathbb{P}^1\times \mathbb{P}^1}(-2,1).
\]
\end{lem}

\begin{proof}
We have the following morphisms: 
\[
f_1 : X \overset{\iota'}{\hookrightarrow}
P' = \mathbb P(\mathcal E') \xrightarrow{\pi'} 
Y_1 = \P^1 \times \P^1, 
\]
where $\mathcal{E}^{\prime}=(f_1)_*\MO_X(-K_X) = \MO_{\mathbb{P}^1\times \mathbb{P}^1}(2,1)\oplus \MO_{\mathbb{P}^1\times \mathbb{P}^1}(1,0)^{\oplus 2}$ (Lemma \ref{l-pic3-C2E-ext}). 
We identify $X$ with $\iota'(X)$, 
which is a smooth prime divisor on $P^\prime$. 
Since $\iota^\prime \colon X \hookrightarrow P^\prime$ is the closed immersion 
induced by the natural surjection $f_1^*(f_1)_*\MO_X(-K_X)\to \MO_X(-K_X)$, 
it follows from \cite[Ch. II, the proof of Propsition 7.12]{Har77} that 
\begin{equation}\label{e1-pic3-C2E-X}
\MO_X(-K_X)\simeq \iota^{\prime *}\MO_{P^\prime}(1) = \MO_{P^\prime}(1)|_X. 
\end{equation}
By Proposition~\ref{standard_2}(2), 
\begin{align}\label{e2-pic3-C2E-X}
  \MO_{P^\prime}(-K_{P^\prime}) &\simeq \MO_{P^\prime}(3)\otimes \pi^{\prime *}(\omega_{\mathbb{P}^1\times \mathbb{P}^1}\otimes \det \mathcal{E}^\prime)^{-1} \notag\\
  &\simeq \MO_{P^\prime}(3) \otimes \pi^{\prime *}(\MO_{\mathbb{P}^1\times \mathbb{P}^1}(-2,-2)\otimes \MO_{\mathbb{P}^1\times \mathbb{P}^1}(4,1))^{-1} \notag\\
  &\simeq \MO_{P^\prime}(3)\otimes \pi^{\prime *}\MO_{\mathbb{P}^1\times \mathbb{P}^1}(-2,1).
\end{align}
By (\ref{e1-pic3-C2E-X}) and (\ref{e2-pic3-C2E-X}), we obtain 
\[
\MO_{P'}(X)|_X \simeq \MO_{P'}(-K_{P'})|_X \otimes \MO_X(K_X) 
\simeq \MO_{P'}(2)|_X \otimes \pi^{\prime *}\MO_{\mathbb{P}^1\times \mathbb{P}^1}(-2,1)|_X. 
\]
Set 
\[
\mathcal N := \MO_{P'}(X)^{-1} \otimes \MO_{P'}(2)  
\otimes \pi^{\prime *}\MO_{\mathbb{P}^1\times \mathbb{P}^1}(-2,1). 
\]
We have $\mathcal N|_X \simeq \MO_X$. 
For a fibre $\zeta$ of $f_1 : X \to \P^1 \times \P^1$, we get 
$\mathcal N \cdot \zeta = (\mathcal N|_X) \cdot \zeta =0$, 
which implies $\mathcal N \simeq \pi'^*\mathcal N'$ for some $\mathcal N' \in \Pic(\P^1 \times \P^1)$. 
Since the pullback of $\mathcal N'$ to $X$ is trivial: $f_1^*\mathcal N' \simeq \mathcal N|_X \simeq \MO_X$, 
it holds that $\mathcal N' \simeq \MO_{\P^1 \times \P^1}$. 
Therefore, $\mathcal N \simeq \pi'^*\mathcal N' \simeq \MO_{P'}$, as required. 
\end{proof}

\begin{prop}\label{p-rho3-C1E}
We use Notation \ref{n-pic3-C2E}. 
Then there exists a closed immersion 
\[
\iota : X \hookrightarrow P:=
\P_{\P^1 \times \P^1}(\MO_{\mathbb{P}^1\times \mathbb{P}^1}\oplus \MO_{\mathbb{P}^1\times \mathbb{P}^1}(-1,-1)^{\oplus 2}) 
\]
such that 
\[
\MO_P(\iota(X)) \simeq \MO_P(2) \otimes \pi^*\MO_{\P^1 \times \P^1}(2, 3),
\]
where $\pi : \P_{\P^1 \times \P^1}(\MO_{\mathbb{P}^1\times \mathbb{P}^1}\oplus \MO_{\mathbb{P}^1\times \mathbb{P}^1}(-1,-1)^{\oplus 2})  \to \P^1 \times \P^1$ 
denotes the induced projection. 
Furthermore, it holds that $(-K_X)^3 = 14$. 
\end{prop}

\begin{proof}
By Lemma \ref{l-pic3-C2E-ext} and Lemma \ref{l-pic3-C2E-X}, we may assume that 
\begin{itemize}
\item $\mathcal E' = \MO_{\mathbb{P}^1\times \mathbb{P}^1}(2,1)\oplus \MO_{\mathbb{P}^1\times \mathbb{P}^1}(1,0)^{\oplus 2}$, 
\item $P' = \P(\mathcal E')$, 
\item $X$ is a smooth prime divisor on $P'$, and
\item $\MO_{P'}(X) \simeq \MO_{P'}(2) \otimes \pi'^*\MO_{\P^1 \times \P^1}(-2, 1)$.  
\end{itemize}
Set $\mathcal{M}:=\MO_{\mathbb{P}^1\times \mathbb{P}^1}(-2,-1)$ and 
\[
\mathcal E := \mathcal E' \otimes \mathcal M = \MO_{\P^1 \times \P^1} 
\oplus \MO_{\P^1 \times \P^1}(-1, -1)^{\oplus 2}. 
\]
We have the induced morphisms: 
\[
\pi : P \xrightarrow{\psi, \simeq} P' \xrightarrow{\pi'} \P^1 \times \P^1
\]
\begin{center}
  \begin{tikzcd}[column sep=50]
    X \arrow[r,"\iota:=\psi^{-1}\circ \iota^\prime"] & P=\mathbb{P}(\MO_{\mathbb{P}^1\times \mathbb{P}^1}\oplus \MO_{\mathbb{P}^1\times \mathbb{P}^1}(-1,-1)^{\oplus 2}) \arrow[d,"\pi :=\pi' \circ \psi"] \\
    & \mathbb{P}^1\times \mathbb{P}^1. 
  \end{tikzcd}
\end{center}
By $\MO_P(1) \simeq \psi^*\MO_{P'}(1) \otimes \pi^*\mathcal M$ \cite[Ch. II, Lemma 7.9]{Har77}, 
it holds that 
\begin{align*}
  \MO_P(\iota(X)) &\simeq \psi^*\MO_{P^\prime}(X) \\
  &\simeq \psi^*\MO_{P^\prime}(2)\otimes \psi^*\pi^{\prime *}\MO_{\mathbb{P}^1\times \mathbb{P}^1}(-2,1) \\
  &\simeq (\psi^*\MO_{P^\prime}(1))^{\otimes 2}\otimes \psi^*\pi^{\prime *}\MO_{\mathbb{P}^1\times \mathbb{P}^1}(-2,1) \\
  &\simeq (\MO_P(1)\otimes \pi^*\mathcal{M}^{-1})^{\otimes 2} \otimes \psi^*\pi^{\prime *}\MO_{\mathbb{P}^1\times \mathbb{P}^1}(-2,1) \\
  &\simeq \MO_P(2)\otimes \pi^*\mathcal{M}^{-2} \otimes \pi^*\MO_{\mathbb{P}^1\times \mathbb{P}^1}(-2,1) \\
  &\simeq \MO_P(2) \otimes \pi^*\MO_{\mathbb{P}^1\times \mathbb{P}^1}(2,3).
\end{align*}

What is remaining is to compute $(-K_X)^3$. 
By Lemma~\ref{standard_2}(5), we have 
\begin{equation*}
  \begin{aligned}
    (-K_X)^3 = 2c_1(\mathcal{E})^2-&2c_2(\mathcal{E})+4(c_1(\mathcal{E})\cdot c_1(\mathcal{F}))+6(c_1(\mathcal{E})\cdot K_{\mathbb{P}^1\times \mathbb{P}^1})\\
    &+9(c_1(\mathcal{F})\cdot K_{\mathbb{P}^1\times \mathbb{P}^1})+6K_{\mathbb{P}^1\times \mathbb{P}^1}^2+3c_1(\mathcal{F})^2, 
  \end{aligned}
\end{equation*}
where $\mathcal{F}:=\MO_{\mathbb{P}^1\times \mathbb{P}^1}(2,3)$. 
By \cite[Appendix A, \S 3, C.3 and C.5]{Har77}, it holds that 
\begin{equation*}
  \begin{aligned}
    c_1(\mathcal{E})&=
2c_1(\MO_{\mathbb{P}^1\times \mathbb{P}^1}(-1,-1))\quad {\rm and}\\
    c_2(\mathcal{E})&=c_1(\MO_{\mathbb{P}^1\times \mathbb{P}^1}(-1,-1))^2 =2.
  \end{aligned}
\end{equation*}
Hence we obtain
\begin{align*}
  (-K_X)^3 = 2\cdot 8-2\cdot 2+  4\cdot (-10)+6\cdot 8 + 9\cdot (-10)+6\cdot8 + 3\cdot 12 =14.
\end{align*}
\end{proof}

\begin{thm}\label{t-pic3-CE}
Let $X$ be a primitive Fano threefold with $\rho(X)=3$. 
Let $R_1$ and $R_2$ be extremal rays of $\NE(X)$. 
Assume that $R_1$ is of type $C$ and $R_2$ is of type $E$. 
Then one and only one of the following holds.  
\begin{enumerate}
\item 
$R_1$ is of type $C_1$, $R_2$ is of type $E_1$, $(-K_X)^3=14$, 
and $X$ is isomorphic to a prime divisor $X'$ on 
$P := \P_{\P^1 \times \P^1}(\MO_{\P^1 \times \P^1} \oplus \MO_{\P^1 \times \P^1}(-1, -1)^{\oplus 2})$ 
such that 
$\MO_P(X') \simeq \MO_P(2) \otimes \pi^*\MO_{\P^1 \times \P^1}(2, 3)$, 
where $\pi : P= \P_{\P^1 \times \P^1}(\MO_{\P^1 \times \P^1} \oplus \MO_{\P^1 \times \P^1}(-1, -1)^{\oplus 2}) \to \P^1 \times \P^1$ denotes the natural projection. 
\item 
$R_1$ is of type $C_2$, $R_2$ is of type $E_1$, $(-K_X)^3=52$, 
and 
$X \simeq \P(\MO_{\P^1 \times \P^1} \oplus \MO_{\P^1 \times \P^1}(1, 1))$. 
\end{enumerate}
\end{thm}

\begin{proof}
The assertion follows from Proposition \ref{p-pic3-C2E1} and Proposition \ref{p-rho3-C1E}. 
\end{proof}

\section{Appendix: Computation of Chern classes}

In this section, we compute some intersection numbers on projective space bundles (Proposition \ref{standard_1}) and blowups (Proposition \ref{standard_3}). 
We include the proofs for the sake of completeness, although both results are obtained just by applying standard results on Chern classes.

\begin{prop}\label{standard_1}\label{standard_2}
Let $Y$ be an $n$-dimensional smooth {\cred projective} variety and let $\mathcal{E}$ be a locally free sheaf of rank $r$ on $Y$. 
Let $\pi \colon X=\mathbb{P}(\mathcal{E})\to Y$ be the projection induced by the $\P^{r-1}$-bundle structure. 
Then the following hold. 
  \begin{enumerate}
\item 
If $r=2$, then it holds that $\MO_X(1) \cdot F =1$ for 
{\cred every}  fibre $F$ of  $\pi : X \to Y$ over a closed point. 
    \item 
    The following holds: 
    \[\omega_X \simeq \MO_X(-r)\otimes \pi^*(\omega_Y\otimes \det \mathcal{E}).  \]
    \item If $n=2$ and $r=2$, then 
    \[(-K_X)^3 = 2c_1(\mathcal{E})^2-8c_2(\mathcal{E})+6K_Y^2.\]
    \item If $n=1$ and $r=2$, then 
    \[c_1(\MO_X(1))^2=\deg_Y(\mathcal{E}).\]
    \item If $n=2$, $r=3,$ and $D$ is a smooth prime divisor on $X$ such that $\MO_X(D)\simeq \MO_X(2)\otimes \pi^*\mathcal{F}$ for some invertible sheaf $\mathcal{F}$ on $Y$, then 
    \begin{align*}
      (-K_D)^3 = &2c_1(\mathcal{E})^2-2c_2(\mathcal{E})+4c_1(\mathcal{E})\cdot c_1(\mathcal{F}) +6 c_1(\mathcal{E})\cdot K_Y \\
      & + 9c_1(\mathcal{F})\cdot K_Y + 6K_Y^2+3c_1(\mathcal{F})^2.
    \end{align*}
  \end{enumerate}
\end{prop}


\begin{proof}
Fix a Cartier divisor $\xi$ on $X$ such that $\MO_X(1) \simeq \MO_X(\xi)$. 
The assertion (1) holds by 
\[
\MO_X(1) \cdot F = \deg_F( \MO_X(1)|_F) = \deg_{\mathbb P^1}(\MO_{\mathbb{P}^1}(1))=1.
\]

Let us show (2).  
Since $\pi$ is locally trivial, it follows from \cite[Ch. II, Proposition 8.11 and Exercise 8.3(a)]{Har77} that we have an exact sequence
  \[ 0 \to \pi^*\Omega_Y\longrightarrow \Omega_X \longrightarrow \Omega_{X/Y} \longrightarrow 0.\]  
By \cite[Ch. II, Exercise 5.16(d),(e)]{Har77}, it holds that 
  \begin{align*}
    \omega_X &\simeq \bigwedge^{n+r-1}\Omega_X \\
    &\simeq \bigwedge^n \pi^*\Omega_Y \otimes \bigwedge^{r-1}\Omega_{X/Y} \\
    &\simeq \pi^*\left(\bigwedge^n \Omega _Y\right)\otimes \bigwedge^{r-1}\Omega_{X/Y} \\
    &\simeq \pi^*\omega_Y\otimes \omega_{X/Y}.
  \end{align*}
  Moreover, we have $\omega_{X/Y}\simeq \MO_X(-r)\otimes \pi^*\det \mathcal{E}$ by \cite[Ch. III, Exercise 8.4(b)]{Har77}. 
  Thus (2) holds.

Let us show   (3). 
Set $\mathcal{L} :=\omega_Y\otimes \det \mathcal{E}$. 
It follows from (2) that $\omega_X \simeq \MO_X(-2) \otimes \pi^*\mathcal L$. 
We then obtain 
  \begin{align}\label{standard2-1}
    (-K_X)^3 = (2\xi-c_1(\pi^*\mathcal{L}))^3 = 8\xi^3-12 \xi^2\cdot \pi^*c_1(\mathcal{L})+6\xi\cdot \pi^*c_1(\mathcal{L})^2.
  \end{align}
  By \cite[Appendix A, \S 3, page 429]{Har77}, we have an equation 
  \[\pi^*c_0(\mathcal{E})\cdot \xi^2-\pi^*c_1(\mathcal{E})\cdot\xi+\pi^*c_2(\mathcal{E})=0.\]
  By $c_0(\mathcal{E})=1$, this can be written as 
  \begin{equation}\label{standard2-2}
    \xi^2=\pi^*c_1(\mathcal{E})\cdot \xi-\pi^*c_2(\mathcal{E}).
  \end{equation} 
We then  get 
  \begin{align}\label{standard2-3}
    \xi^3 &= \pi^*c_1(\mathcal{E})\cdot \xi^2-\pi^*c_2(\mathcal{E})\cdot \xi \notag\\
    &= \pi^*c_1(\mathcal{E})\cdot(\pi^*c_1(\mathcal{E})\cdot \xi-\pi^*c_2(\mathcal{E}))-\pi^*c_2(\mathcal{E})\cdot \xi \notag\\
    &= (\pi^*c_1(\mathcal{E}))^2\cdot \xi-\pi^*c_1(\mathcal{E})\cdot \pi^*c_2(\mathcal{E})-\pi^*c_2(\mathcal{E})\cdot \xi\\
    &= \pi^*(c_1(\mathcal{E})^2-c_2(\mathcal{E}))\cdot \xi-\pi^*(c_1(\mathcal{E})\cdot c_2(\mathcal{E})) \notag\\
    &= \pi^*(c_1(\mathcal{E})^2-c_2(\mathcal{E}))\cdot \xi, \notag
  \end{align}
where the first and second equalities follows from (\ref{standard2-2}) and 
the last one holds, because $\dim Y=2$ implies $c_1(\mathcal{E})\cdot c_2(\mathcal{E})=0$. 
  By (\ref{standard2-1}), (\ref{standard2-2}), and (\ref{standard2-3}), we obtain 
 \begin{align*}
    (-K_X)^3 
&=8\pi^*(c_1(\mathcal{E})^2-c_2(\mathcal{E}))\cdot \xi- 12(\pi^*c_1(\mathcal{E})\cdot \xi-\pi^*c_2(\mathcal{E}))\cdot \pi^*c_1(\mathcal{L})+6\xi\cdot \pi^*c_1(\mathcal{L})^2\\
&= 8\pi^*(c_1(\mathcal{E})^2-c_2(\mathcal{E}))\cdot \xi-12\pi^*(c_1(\mathcal{E})\cdot c_1(\mathcal{L}))\cdot \xi  + 6\pi^*c_1(\mathcal{L})^2\cdot \xi,\\
& = 8(c_1(\mathcal{E})^2-c_2(\mathcal{E}))-12c_1(\mathcal{E})\cdot c_1(\mathcal{L})+6c_1(\mathcal{L})^2.
  \end{align*}
where the second equality holds by $c_2(\mathcal{E})\cdot c_1(\mathcal{L}) =0$ and 
the third one follows from  $\pi^*Z\cdot \xi=\deg Z$ for {\cred every}  $0$-cycle $Z\in \mathrm{A}^2(Y)$ 
(note that $\deg$ is dropped by abuse of notation). 
By $c_1(\mathcal L) = c_1(\omega_Y) + c_1(\det \mathcal{E}) =K_Y +c_1(\mathcal{E})$, we obtain 
  \begin{align*}
    (-K_X)^3 &= 8(c_1(\mathcal{E})^2-c_2(\mathcal{E}))-12c_1(\mathcal{E})\cdot (K_Y+c_1(\mathcal{E})) +6(K_Y+c_1(\mathcal{E}))^2 \\
    &= 2c_1(\mathcal{E})^2-8c_2(\mathcal{E})+6K_Y^2.
  \end{align*}
Thus (3) holds.

Let us show (4). 
It follows from \cite[Appendix A, \S 3, page 429]{Har77} that 
\[
\xi^2=\pi^*c_1(\mathcal{E})\cdot \xi-\pi^*c_2(\mathcal{E}).
\]
By $\dim Y=1$, we have $c_2(\mathcal{E})=0$, which implies 
\[
c_1(\MO_X(1))^2 = \xi^2=\pi^*c_1(\mathcal{E})\cdot \xi = \deg_Y(\mathcal{E}).\]
Thus (4) holds. 


Let us show  (5). 
By the adjunction formula,
  \begin{align*}
    (-K_D)^3 = ((-K_X-D)|_D)^3 = (-K_X-D)^3\cdot D.
  \end{align*}
  By (2) and our assumption, 
  we get 
  $K_X\sim -3\xi+\pi^*c_1(\omega_Y\otimes \det \mathcal{E})$ and 
  $D\sim 2\xi+\pi^*c_1(\mathcal{F})$, respectively. 
  Thus 
  \[
  K_X+D \sim -\xi +\pi^*c_1(\omega_Y \otimes \det \mathcal E \otimes \mathcal F) 
  = -\xi +\pi^*c_1(\mathcal G)\]
  \[
  \text{for}\qquad  \mathcal{G}:=\omega_Y\otimes \det \mathcal{E}\otimes \mathcal{F}.
  \]
  Then 
  \begin{align}\label{standard2-4}
    (-K_D)^3 &= (-K_X-D)^3\cdot D\notag\\
&= (\xi -\pi^*c_1(\mathcal G) )^3 \cdot (2\xi+\pi^*c_1(\mathcal{F})) \notag\\
&= (\xi^3  -3 \pi^*c_1(\mathcal G) \cdot \xi^2 +3 \pi^*c_1(\mathcal G)^2 \cdot \xi -(\pi^*c_1(\mathcal G))^3) \cdot (2\xi+\pi^*c_1(\mathcal{F}))\\
&= (\xi^3  -3 \pi^*c_1(\mathcal G) \cdot \xi^2 +3 \pi^*c_1(\mathcal G)^2 \cdot \xi) \cdot (2\xi+\pi^*c_1(\mathcal{F}))\notag\\
&= 2\xi^4 + (-6 \pi^*c_1(\mathcal G)+ \pi^*c_1(\mathcal{F})) \cdot \xi^3  + (6 \pi^*c_1(\mathcal G)^2-3 \pi^*(c_1(\mathcal G) \cdot c_1(\mathcal{F})) ) \cdot \xi^2,  
\notag
  \end{align}
where the fourth equality holds by 
$c_1(\mathcal{G})^3=0$ and 
the last one follows from 
$c_1(\mathcal{G})^2\cdot c_1(\mathcal{F})=0$. 
  By \cite[Appendix A, \S 3, page 429]{Har77}, we have an  equation 
  \[\pi^*c_0(\mathcal{E})\cdot \xi^3-\pi^*c_1(\mathcal{E})\cdot\xi^2+\pi^*c_2(\mathcal{E})\cdot \xi -\pi^*c_3(\mathcal{E})=0.\] 
By $\dim Y=2$, we have $c_3(\mathcal{E})=0$. 
It holds by $c_0(\mathcal{E})=1$ that 
  \begin{equation}\label{standard2-5}
    \xi^3=\pi^*c_1(\mathcal{E})\cdot \xi^2-\pi^*c_2(\mathcal{E})\cdot \xi.    
  \end{equation}
  Moreover, we obtain 
  \begin{align}\label{standard2-6}
    \xi^4 &= \pi^*c_1(\mathcal{E})\cdot \xi^3-\pi^*c_2(\mathcal{E})\cdot \xi^2 \\
    &= \pi^*c_1(\mathcal{E})\cdot(\pi^*c_1(\mathcal{E})\cdot \xi^2-\pi^*c_2(\mathcal{E})\cdot \xi)-\pi^*c_2(\mathcal{E})\cdot \xi^2 \notag\\
    &= \pi^*(c_1(\mathcal{E})^2-c_2(\mathcal{E}))\cdot \xi^2,\notag
  \end{align}
where we used  $c_1(\mathcal{E})\cdot c_2(\mathcal{E})=0$, 
which is guaranteed by $\dim Y=2$.
  By (\ref{standard2-4}), (\ref{standard2-5}) and (\ref{standard2-6}), 
  {\cred we get} 
  \begin{align}\label{standard2-7}
(-K_D)^3 
&= 2\xi^4 + (-6 \pi^*c_1(\mathcal G)+\pi^*c_1(\mathcal{F})) \cdot \xi^3  + (6 \pi^*c_1(\mathcal G)^2-3 \pi^*(c_1(\mathcal G) \cdot c_1(\mathcal{F}) ) \cdot \xi^2\notag\\
&= 2\pi^*(c_1(\mathcal{E})^2-c_2(\mathcal{E}))\cdot \xi^2 \notag\\
&+ 
(-6 \pi^*c_1(\mathcal G)+\pi^*c_1(\mathcal{F})) \cdot 
(\pi^*c_1(\mathcal{E})\cdot \xi^2-\pi^*c_2(\mathcal{E})\cdot \xi)\notag\\
&+ (6 \pi^*c_1(\mathcal G)^2-3 \pi^*(c_1(\mathcal G) \cdot c_1(\mathcal{F}) ) \cdot \xi^2\\
&= 2\pi^*(c_1(\mathcal{E})^2-c_2(\mathcal{E}))\cdot \xi^2 \notag\\
&+ 
(-6 \pi^*c_1(\mathcal G)+\pi^*c_1(\mathcal{F})) \cdot 
\pi^*c_1(\mathcal{E})\cdot \xi^2 \notag\\
&+ (6 \pi^*c_1(\mathcal G)^2-3 \pi^*(c_1(\mathcal G) \cdot c_1(\mathcal{F}) ) \cdot \xi^2 \notag
  \end{align}
where the last equality holds by 
$c_1(\mathcal{G})\cdot c_2(\mathcal{E})=c_1(\mathcal{F})\cdot c_2(\mathcal{E})=0$.

  Since we have $\MO_X(1)|_F =\MO_{\P^2}(1)$ for every fibre $F\simeq \mathbb{P}^2$ of $\pi$, 
  it holds that 
  $\xi^2\cdot F = (\xi|_F)^2 = c_1(\MO_{\mathbb{P}^2}(1))^2 =1$, which implies  $\pi^*Z\cdot \xi^2=\deg Z$ for {\cred every}  $0$-cycle $Z\in \mathrm{A}^2(Y)$. 
It follows from 
(\ref{standard2-7}) 
that  
  \begin{align*}
    (-K_D)^3 &= 
2c_1(\mathcal{E})^2-2c_2(\mathcal{E}) 
+ (-6 c_1(\mathcal G)+c_1(\mathcal{F})) \cdot c_1(\mathcal{E}) 
+ (6 c_1(\mathcal G)^2-3 c_1(\mathcal G) \cdot c_1(\mathcal{F}) ) 
  \end{align*}
By $c_1(\mathcal{G})= K_Y+c_1(\mathcal{E}) + c_1(\mathcal{F})$, we obtain 
  \begin{align*}
&-6 c_1(\mathcal G) \cdot c_1(\mathcal{E}) 
+ 6 c_1(\mathcal G)^2-3 c_1(\mathcal G) \cdot c_1(\mathcal{F}) \\
=& c_1(\mathcal G) \cdot ( -6 c_1(\mathcal{E})  +6 c_1(\mathcal G)
-3 c_1(\mathcal{F}))\\
=& ( K_Y+c_1(\mathcal{E}) + c_1(\mathcal{F})) \cdot (6K_Y+3 c_1(\mathcal{F})). 
\end{align*}
Therefore, 
  \begin{align*}
    (-K_D)^3 &= 
2c_1(\mathcal{E})^2-2c_2(\mathcal{E}) 
+ c_1(\mathcal{F}) \cdot c_1(\mathcal{E}) \\
&+ ( K_Y+c_1(\mathcal{E}) + c_1(\mathcal{F})) \cdot (6K_Y+3 c_1(\mathcal{F}))\\
&= 2c_1(\mathcal{E})^2-2c_2(\mathcal{E}) +4c_1(\mathcal{F}) \cdot c_1(\mathcal{E})\\
&+6K_Y^2 + 9 K_Y \cdot c_1(\mathcal F)+6 K_Y \cdot c_1(\mathcal E)+3c_1(\mathcal F)^2. 
  \end{align*}
  Thus (5) holds. 
\end{proof}

\begin{prop}\label{standard_3}
Let $Y$ be a smooth projective threefold and let $C$ be a smooth curve on $Y$. 
 Let $f\colon X\to Y$ be the blowup along $C$. 
 Set $D := \Ex(f)$. 
  Then the following hold:
  \begin{enumerate}
    \item $K_X\sim f^*K_Y + D$.
    \item The induced morphism $f|_D\colon D\to C$ coincides with the $\mathbb{P}^1$-bundle $\mathbb{P}(\mathcal{N}^*_{C/Y})$ over $C$, and $\MO_X(-D)|_D\simeq \MO_{\mathbb{P}(\mathcal{N}^*_{C/Y})}(1)$.
    \item $D^3 =\deg_C(\mathcal{N}^*_{C/Y})$.
    \item $(-K_X)^2\cdot D  = -K_Y\cdot C +2-2g(C)$.   
  \end{enumerate}
\end{prop}

\begin{proof}
The assertions (1) and (2) follow from  \cite[Ch. II, Exercise 8.5(b)]{Har77} and 
and \cite[Ch. II, Theorem 8.24(b)(c)]{Har77}, respectively. 
Let us show (3). 
 By (2), we have $D\simeq \mathbb{P}(\mathcal{N}^*_{C/Y})$. 
 We then get 
  \begin{align*}
   D^3 &= (\MO_X(-D)|_D) \cdot (\MO_X(-D)|_D) \\
    &= c_1(\MO_D(1))^2 \\
    &= \deg_C(\mathcal{N}^*_{C/Y}), 
  \end{align*}
where the second equality holds by (2) and 
the last one follows from Proposition~\ref{standard_2}(4). 
Thus (3) holds. 

Let us show (4).  
Set $\pi: =f|_D \colon D=\mathbb{P}(\mathcal{N}^*_{C/Y})\to C$ and 
$\mathcal{L} := (\omega_C\otimes \det \mathcal{N}^*_{C/Y})^{-1}$. 
Fix a divisor $\xi$ on $D$ such that $\MO_D(\xi) =\MO_D(1)$.
It holds that 
  \begin{align*}
    \MO_X(-K_X)|_D &\simeq \MO_D(-K_D)\otimes \MO_X(D)|_D \\
    &\simeq \MO_D(2)\otimes \pi^*(\omega_C\otimes \det \mathcal{N}^*_{C/Y})^{-1}\otimes \MO_D(-1) \\
    &\simeq \MO_D(1)\otimes \pi^* \mathcal{L},
  \end{align*}
  where the first isomorphism holds by the adjunction formula and 
the second one follows from (2) and Proposition~\ref{standard_2}(2).
  Hence 
  \begin{align*}
    (-K_X)^2\cdot D   &= (\xi+c_1(\pi^*\mathcal{L}))^2 \\
    &= \xi^2+2\xi\cdot \pi^*c_1(\mathcal{L})+\pi^*c_1(\mathcal{L})^2 \\
    &= \xi^2+2\pi^*c_1(\mathcal{L})\cdot \xi, 
  \end{align*}
  where $c_1(\mathcal{L})^2=0$ holds by $\dim C=1$. 
By \cite[Appendix A, \S 3, page 429]{Har77}, we have 
$\xi^2=\pi^*c_1(\mathcal{N}^*_{C/Y})\cdot \xi$. 
Therefore, we get 
  \begin{align*}
   (-K_X)^2\cdot D &= \pi^*c_1(\mathcal{N}^*_{C/Y})\cdot \xi + 2\pi^*c_1(\mathcal{L})\cdot \xi \\
    &= \pi^*(c_1(\mathcal{N}^*_{C/Y})+2c_1(\mathcal{L}))\cdot \xi\\
&= \deg_C(c_1(\mathcal{N}^*_{C/Y})+2c_1(\mathcal{L}))\\
&= \deg_C(-c_1(\mathcal{N}^*_{C/Y})-2K_C)\\
&= \deg_C(-c_1(\omega_Y|_C))-K_C)\\
&= -K_Y \cdot C +2-2g(C), 
  \end{align*}
where the fourth equality holds by $\mathcal{L} = (\omega_C\otimes \det \mathcal{N}^*_{C/Y})^{-1}$ 
and 
the fifth equality follows from $\omega_C \simeq \omega_Y|_C \otimes (\det \mathcal N_{C/Y})$ 
\cite[Ch. II, Proposition 8.20]{Har77}.  
\end{proof}

\section{Appendix: Description of split double covers}\label{s-dc}


\begin{dfn}\label{d-simple-mu2}
Let $Y$ be a normal variety and let $\mathcal L$ be an invertible sheaf. 
Set 
\[
\mathbb V
 := \Spec_Y \bigoplus_{i \geq 0} \mathcal L^{\otimes -i}, 
\] 
where $\bigoplus_{i \geq 0}\mathcal L^{\otimes -i}$ denotes the natural graded $\MO_Y$-algebra. 
Let $\pi: \mathbb V \to Y$ be the projection. 
Take the tautological section $z \in H^0(\mathbb V, \pi^*\mathcal L)$, which corresponds to 
the element $(0, 1, 0, 0, ...) \in \bigoplus_{i \geq 0} H^0(Y, \mathcal L^{1-i})$ via 
\[
H^0(\mathbb V, \pi^*\mathcal L) \simeq H^0(Y, \pi_*\pi^*\mathcal L) \simeq \bigoplus_{i \geq 0} H^0(Y, \mathcal L^{1-i}). 
\]
\begin{enumerate}
\item 
For $s \in H^0(Y, \mathcal L^{\otimes 2})$, the closed subscheme 
\[
X := \{ z^2 - \pi^*(s) =0 \} \subset \mathbb V
\]
and the induced morphism $X \to Y$ are called the {\em simple $\mu_2$-cover} associated to  $(\mathcal L, s)$. 
\item 
Assume $p=2$. 
For $s \in H^0(X, \mathcal L)$ and $t \in H^0(X, \mathcal L^{\otimes 2})$, 
the closed subscheme 
\[
X := \{ z^2 - \pi^*(s) \cdot z + \pi^*(t) =0 \} \subset \mathbb V
\]
and the induced morphism $X \to Y$ are called the {\em split $\alpha_{\mathcal L, s}$-torsor}  associated to $(\mathcal L, s, t)$. 
Note that 
the split $\alpha_{\mathcal L, 0}$-torsor associated to $(\mathcal L, 0, t)$ 
coincides with the simple $\mu_2$-cover associated to $(\mathcal L, -t)$. 
\end{enumerate}
\end{dfn}




\begin{lem}\label{l-dc-odd}
Assume that $p \neq 2$. 
Let $f: X \to Y$ be a double cover of smooth projective varieties. 
Set $\mathcal L := (f_*\MO_X/\MO_Y)^{-1}$, which is an invertible sheaf (Remark \ref{r-L-inv}). 
Then $f$ is a simple $\mu_2$-cover associated to $(\mathcal L, s)$ for some $s \in H^0(Y, \mathcal L^{\otimes 2})$
\end{lem}

\begin{proof}
See \cite[Page 10]{CD89}. 
\end{proof}


\begin{lem}\label{l-dc-char2}
Assume $p = 2$. 
Let $f: X \to Y$ be a double cover of smooth projective varieties. 
Set $\mathcal L := (f_*\MO_X/\MO_Y)^{-1}$, which is an invertible sheaf (Remark \ref{r-L-inv}). 
Then $f$ is a {\cred split} double cover if and only if 
$f$ is a split $\alpha_{\mathcal L, s}$-torsor associated to some $(\mathcal L, s, t)$. 
\end{lem}

\begin{proof}
If $f$ is the split $\alpha_{\mathcal L, s}$-torsor associated to some $(\mathcal L, s, t)$, 
then it is clear that $f$ is a split double cover. 
Conversely, assume that $f$ is a split double cover. 
By \cite[Page 10-11]{CD89}, $f$ is given by a splittable admissible triple in the sense of \cite[Page 11]{CD89}, which is nothing but 
a split $\alpha_{\mathcal L, s}$-torsor associated to some $(\mathcal L, s, t)$. 
\end{proof}

\begin{lem}\label{l-trivial-double}
Let $f: X \to Y$ be a split double cover of smooth projective varieties. 
Set $\mathcal L := (f_*\MO_X/\MO_Y)^{-1}$, which is an invertible sheaf (Remark \ref{r-L-inv}).  
Let $\pi : Y \to Y'$ be a morphism of smooth projective varieties 
such that ${\cred \pi}_*\MO_Y = \MO_{Y'}$ and $\mathcal L \simeq {\cred \pi}^*\mathcal L'$ for some invertible sheaf $\mathcal L'$ on $Y'$. 
Then there exist a split double cover $f':X' \to Y'$ and a morphism $\pi_X : X \to X'$ 
such that the following diagram is cartesian: 
\[
\begin{tikzcd}
X \arrow[r, "\pi_X"] \arrow[d, "f"] & X' \arrow[d, "f'"]\\
Y \arrow[r, "\pi"] & Y'. 
\end{tikzcd}
\]
In particular, if $Y = Y' \times_k Y''$ and $\pi: Y \to Y'$ is the first projection in addition to the above assumptions, 
then $X \simeq X' \times_k Y''$ and $\pi_X : X \to X'$ coincides with the first projection. 
\end{lem}

\begin{proof}
Assume $p=2$. 
By Lemma \ref{l-dc-char2}, there exist $s \in H^0(Y, \mathcal L)$ and $t \in H^0(Y, \mathcal L^{\otimes 2})$ 
such that $f$ is the split $\alpha_{\mathcal L, s}$-torsor associated to $(\mathcal L, s, t)$. 
By $\pi_*\MO_Y = \MO_{Y'}$ and $\mathcal L \simeq {\cred \pi}^*\mathcal L'$, 
we have the following corresponding elements: 
\begin{itemize}
    \item $H^0(Y, \mathcal L) \simeq H^0(Y', \mathcal L'), \qquad s \mapsto s'$. 
    \item $H^0(Y, \mathcal L^{\otimes 2}) \simeq H^0(Y', \mathcal L'^{\otimes 2}), \qquad t \mapsto t'$. 
\end{itemize}
Set $f' : X' \to Y'$ to be the split $\alpha_{\mathcal L', s'}$ torsor associated to $(\mathcal L', s', t')$. 

Fix an affine open cover $Y' = \bigcup_{i \in I} Y'_i$ such that $\mathcal L'|_{Y'_i} \simeq \MO_{Y'_i}$ for every $i \in I$. 
Then $X'$ is given as follows: 
\[
X' = \bigcup_{i \in I} \Spec\,(\MO_{Y'_i}(Y'_i)[z'_i]/(z'^2_i +s'_i z'_i+t'_i)), 
\] 
where $s'_i := s'|_{Y'_i}$ and $t'_i := t'|_{Y'_i}$. 
By taking an affine open cover of $Y$ which refines $Y = \bigcup_{i \in I} \pi^{-1}(Y'_i)$, 
we can directly check that the induced morphism $X \to X' \times_{Y'} Y$ is an isomorphism (more explicitly, 
if $\widetilde Y$ is an affine open subset contained in $\pi^{-1}(Y'_i)$ for some $i$, 
then both $X$ and $X' \times_{Y'} Y$ are given by $\MO_{\widetilde Y}(\widetilde Y)[z]/(z^2 + (s|_{\widetilde Y}) z + (t|_{\widetilde Y})$). 
If $p\neq 2$, the same argument works by using Lemma \ref{l-dc-odd} instead of Lemma \ref{l-dc-char2}. 
\end{proof}

\section{Classification table of Fano threefolds with $\rho=2$}\label{s table rho2}

{\cred 
\begin{dfn}\label{d No pic2}
We say that $X$ is a Fano threefold of {\em No.}\  2-xx if 
$X$ is a Fano threefold with $\rho(X)=2$ such that 
$(-K_X)^3$ and 
the types of the extremal rays are as in Table \ref{table-pic2} below. 
For example, No.\ 2-1 and No.\ 2-15 are defined as follows. 
\begin{itemize}
\item 
$X$ is a Fano threefold of {\em No.\ 2-1} if 
$X$ a Fano threefold such that $\rho(X)=2$, 
$(-K_X)^3 = 4$, 
one of the extremal rays is of type $D_1$, and the other is of type $E_1$. 
\item 
$X$ is a Fano threefold of {\em  No.\ 2-15} if 
$X$ a Fano threefold such that $\rho(X)=2$, 
$(-K_X)^3 = 22$, one of the extremal rays is of type $E_1$, and the other is of type $E_3$ or $E_4$. 
\end{itemize}
\end{dfn}
}

Let $X$ be a Fano threefold with $\rho(X)=2$. 
Then   one of No.\ 2-1, ..., No.\ 2-36 
 holds. 
In Table \ref{table-pic2}, we use the {\cred following} notation and terminologies.  
\begin{enumerate}
{\cred 
\item We say that $f : X \to Y$ is a split double cover if 
$f$ is a finite surjective morphism of projective normal varieites such that $\MO_Y \to f_*\MO_X$ splits as an $\MO_Y$-module homomorphism and the induced field extension $K(X) \supset K(Y)$ is of degree two. 
For a split double cover $f : X \to Y$, we set $\mathcal L := (f_*\MO_X/\MO_Y)^{-1}$, which is an invertible sheaf on $Y$ (Remark \ref{r-L-inv}). 
\item $Q$ is the smooth quadric hypersurace in $\P^4$. 
\item $W$ is a smooth prime divisor on $\P^2 \times \P^2$ of bidegree $(1, 1)$. 
Note that such a threefold is unique up to isomorphisms (Lemma \ref{l-bideg11}). 
\item For $1 \leq d \leq 7$ with $d \neq 6$, 
let $V_d$ be a Fano threefold of index $2$ satisfying $(-K_{V_d}/2)^3 =d$. 
Note that $V_d$ is not determined uniquely by $d$. 
\item The centres of all the blowups are smooth curves. 
The degree of a curve $B$ on a Fano threefold $Y$ is defined as $-\frac{1}{r_Y}K_Y \cdot B$, 
where $r_Y$ denotes the index of $Y$. 
A line (resp. conic) is a smooth rational curve of degree one (resp. two). 
Note that if $|-\frac{1}{r_Y}K_Y|$ is very ample, then a curve of degree one (resp. two) is automatically a line (resp. conic). 
}
\item 
If an extremal ray is of type $C$, 
then $\deg \Delta$ denotes the degree of the discriminant bundle $\Delta$, which is an invertible sheaf on $\P^2$. 
When the conic bundle is generically smooth, then $\Delta$ coincides with the invertible sheaf associated with the discriminant divisor. 
\item 
If an extremal ray is of type $D$, 
then $X_t$ denotes a fibre over a closed point $t \in \P^1$. 
In particular, $(-K_X)^2 \cdot X_t$ coincides with $(-K_{X_K})^2$ for the generic fibre $X_K$. 
\end{enumerate}

  \begin{center}
\begin{longtable}{ccp{10cm}cc}
      \caption{Fano threefolds with $\rho(X)=2$}\label{table-pic2}\\
No. & $(-K_X)^3$ & descriptions and extremal rays  \\ \hline
\hyperref[p-E1D1]{2-1}\label{table-2-1} & $4$ & $D_1: (-K_X)^2 \cdot X_t=1$ \\ 
&  & 
$E_1:$ blowup of $V_1$ along an elliptic curve of degree $1$  which is a complete intersection of two members of $|-\frac{1}{2}K_{V_1}|$ & \\ \hline
\hyperref[t-pic2-D]{2-2}\label{table-2-2}
& $6$ & $X$ is a split double cover of $\mathbb{P}^2\times \mathbb{P}^1$ with $\mathcal L \simeq \MO(2, 1)$ & \\ 
 &  & $C_1: \deg \Delta=8$, $X \xrightarrow{{\rm 2:1}} \P^2 \times \P^1 \xrightarrow{{\rm pr}_1} \P^2$  & \\ 
&  & $D_1: (-K_X)^2 \cdot X_t = 2$& \\ \hline
\hyperref[p-E1D1]{2-3}\label{table-2-3} & $8$ & $D_1: (-K_X)^2 \cdot X_t=2$ & \\ 
 &  & $E_1:$ blowup of $V_2$ along an elliptic curve of degree $2$  which is a complete intersection of two members of $|-\frac{1}{2}K_{V_2}|$ & \\ \hline
\hyperref[p-E1D1]{2-4}\label{table-2-4} & $10$ & $D_1: (-K_X)^2 \cdot X_t=3$ & \\ 
 &  & $E_1:$ blowup of $\P^3$ along a curve of genus $10$ and degree $9$  which is a complete intersection of two cubic surfaces & \\ \hline
\hyperref[p-E1D1]{2-5}\label{table-2-5} & $12$ & $D_1: (-K_X)^2 \cdot X_t =3$ & \\ 
&  & $E_1:$ blowup of $V_3$ along an elliptic curve of degree $3$  which is a complete intersection of two members of $|-\frac{1}{2}K_{V_3}|$ & \\ \hline
\hyperref[t-pic2-C]{2-6}\label{table-2-6}
& $12$ & 
a smooth divisor on $\mathbb{P}^2\times \mathbb{P}^2$ of bidegree $(2,2)$, or 
a split double cover of $W$ with $\mathcal L^{\otimes 2} \simeq \omega_W^{-1}$ & \\
 &  & $C_1: \deg \Delta=6$, $X \to \P^2 \times \P^2 \xrightarrow{{\rm pr_1}} \P^2$  & \\ 
&  &  $C_1: \deg \Delta=6$, $X \to \P^2 \times \P^2 \xrightarrow{{\rm pr_2}} \P^2$ & \\ \hline
\hyperref[p-E1D1]{2-7}\label{table-2-7} & $14$ & $D_1: (-K_X)^2 \cdot X_t=4$ & \\ 
 &  & $E_1:$ blowup of $Q$ along a curve of genus $5$ and degree $8$ which is a complete intersection of two members of $|\MO_Q(2)|$ & \\ \hline
\hyperref[t-pic2-E]{2-8}\label{table-2-8} & $14$ &  a split double cover of $V_7$ 
 with $\mathcal L^{\otimes 2} \simeq \omega_{V_7}^{-1}$ & \\ 
 &  &  $C_1: \deg \Delta = 6$, $X \xrightarrow{{\rm 2:1}} 
 V_7 = \mathbb{P}_{\P^2}(\MO_{\mathbb{P}^2}\oplus \MO_{\mathbb{P}^2}(1)) \xrightarrow{{\rm pr}} \P^2$  & \\ 
&  & $E_3\,{\rm or}\,E_4$& \\ \hline
\hyperref[p-E1C1]{2-9}\label{table-2-9} & $16$ & $C_1: \deg \Delta =5$ & \\ 
 &  & $E_1:$ blowup of $\P^3$ along a curve of genus $5$ and degree $7$ & \\ \hline
\hyperref[p-E1D1]{2-10}\label{table-2-10} & $16$ & $D_1:$ $(-K_X)^2 \cdot X_t=4$ & \\ 
 &  & $E_1:$ blowup of $V_4$ along an elliptic curve of degree $4$  which is a complete intersection of two members of $|-\frac{1}{2}K_{V_4}|$ & \\ \hline
\hyperref[p-E1C1]{2-11}\label{table-2-11} & $18$ & $C_1:$ $\deg \Delta =5$ & \\ 
 &  & $E_1:$ blowup of $V_3$ along a line & \\ \hline
\hyperref[p-E1E1]{2-12}\label{table-2-12} & $20$ & $E_1:$ blowup of $\P^3$ along a curve of genus $3$ and degree $6$ & \\ 
 &  & $E_1:$ blowup of $\P^3$ along a curve of genus $3$ and degree $6$ & \\ \hline
\hyperref[p-E1C1]{2-13}\label{table-2-13} & $20$ & $C_1:$ $\deg \Delta =4$ & \\ 
 &  & $E_1:$ blowup of $Q$ along a curve of genus $2$ and degree $6$ & \\ \hline
\hyperref[p-E1D1]{2-14}\label{table-2-14} & $20$ & $D_1:$ $(-K_X)^2 \cdot X_t = 5$ & \\ 
 &  & $E_1:$ blowup of $V_5$ along an elliptic curve of degree $5$  which is a complete intersection of two members of $|-\frac{1}{2}K_{V_5}|$ & \\ \hline
\hyperref[p-E1E3E4]{2-15}\label{table-2-15} & $22$ & $E_1:$ blowup of $\P^3$ along a curve of genus ${\cred 4}$ and degree $6$ & \\ 
&  &$E_3\,{\rm or}\,E_4$ & \\ \hline
\hyperref[p-E1C1]{2-16}\label{table-2-16} & $22$ & $C_1: \deg \Delta =4$ & \\ 
 &  & $E_1:$ blowup of $V_4$ along a conic & \\ \hline
\hyperref[p-E1E1]{2-17}\label{table-2-17} & $24$ & $E_1:$ blowup of $\P^3$ along an elliptic curve of degree $5$ & \\
 &  & $E_1:$ blowup of $Q$ along an elliptic curve of degree $5$ & \\ \hline
\hyperref[t-pic2-D]{2-18}\label{table-2-18} & $24$ & a split double cover of $\mathbb{P}^2\times \mathbb{P}^1$ 
      with $\mathcal L \simeq \MO_{\P^2 \times \P^1}(1, 1)$ & \\
 &  & $C_1: \deg \Delta =4$, $X \xrightarrow{{\rm 2:1}} \P^2 \times \P^1 \xrightarrow{\pr_1} \P^2$  & \\ 
&  & $D_2: (-K_X)^2 \cdot X_t = 8$ & \\ \hline
\hyperref[p-E1E1]{2-19}\label{table-2-19} & $26$ & $E_1:$ blowup of $\P^3$ along a curve of genus $2$ and degree $5$ &\\ 
 &  & $E_1:$ blowup of $V_4$ along a line & \\ \hline
\hyperref[p-E1C1]{2-20}\label{table-2-20} & $26$ & $C_1: \deg \Delta=3$ & \\ 
&  & $E_1:$  blowup of $V_5$ along a cubic rational curve & \\ \hline
\hyperref[p-E1E1]{2-21}\label{table-2-21}& $28$ & 
$E_1:$ blowup of $Q$ along a rational curve of degree $4$ & \\
  &  & 
  $E_1:$ blowup of $Q$ along a rational curve of degree $4$ & \\ \hline
\hyperref[p-E1E1]{2-22}\label{table-2-22} & $30$ & $E_1:$ blowup of $\P^3$ along a rational curve of degree $4$ & \\ 
&  & $E_1:$ blowup of $V_5$ along a conic & \\ \hline
\hyperref[p-E1E3E4]{2-23}\label{table-2-23} & $30$ & 
$E_1:$ blowup of $Q$ along an elliptic curve of degree $4$ & \\ 
&  & $E_3\,{\rm or}\,E_4$ &\\ \hline
\hyperref[t-pic2-C]{2-24}\label{table-2-24}
& $30$ & a  divisor on $\mathbb{P}^2\times \mathbb{P}^2$ of bidegree $(1,2)$ & \\
 &  & $C_1: \deg \Delta = 3$ & \\ 
&  & $C_2$ & \\ \hline
\hyperref[p-E1D2]{2-25}\label{table-2-25} & $32$ & $D_2: (-K_X)^2 \cdot X_t =8$ & \\ 
 &  & $E_1:$ blowup of $\P^3$ along an elliptic curve of degree $4$  which is a complete intersection of two quadric  surfaces & \\ \hline
\hyperref[p-E1E1]{2-26}\label{table-2-26} & $34$ & $E_1:$ blowup of $Q$ along a cubic rational curve  & \\ 
& & $E_1:$ blowup of $V_5$ along a line  & \\ \hline
\hyperref[p-E1C2]{2-27}\label{table-2-27} & $38$ &$C_2$ & \\ 
 &  & $E_1:$ blowup of $\P^3$ along a cubic rational curve & \\ \hline
\hyperref[p-E1E5]{2-28}\label{table-2-28} & $40$ & 
$E_1:$ blowup of $\P^3$ along an elliptic curve of degree $3$ & \\ 
&  &${\cred E_5}$ & \\ \hline
\hyperref[p-E1D2]{2-29}\label{table-2-29} & $40$ & $D_2: (-K_X)^2 \cdot X_t =8$ & \\ 
 &  & $E_1:$ blowup of $Q$ along a conic  which is a complete intersection of two members of $|\MO_Q(1)|$ & \\ \hline
\hyperref[p-E1E2]{2-30}\label{table-2-30} & $46$ & $E_1:$ blowup of $\P^3$ along a conic & \\ 
&  & $E_2:$ blowup of $Q$ at a point & \\ \hline
\hyperref[p-E1C2]{2-31}\label{table-2-31} & $46$ & $C_2$ & \\ 
&  & $E_1:$ blowup of $Q$ along a line & \\ \hline
\hyperref[t-pic2-C]{2-32}\label{table-2-32} & $48$ & $W$, a divisor on $\P^2 \times \P^2$ of bidegree $(1, 1)$ &\\
 &  & $C_2: W \hookrightarrow \P^2 \times \P^2 \xrightarrow{\pr_1} \P^2$ & \\ 
 &  & $C_2: W \hookrightarrow \P^2 \times \P^2 \xrightarrow{\pr_1} \P^2$ & \\ \hline
\hyperref[p-E1D3]{2-33}\label{table-2-33} & $54$ & $D_3: (-K_X)^2 \cdot X_t =9$ & \\ 
 &  & $E_1:$ blowup of $\P^3$ along a line & \\ \hline
\hyperref[t-pic2-D]{2-34}\label{table-2-34} & $54$ & $\P^2 \times \P^1$ \\
 & &  $C_2:$ the projection $\P^2 \times \P^1 \to \P^2$ & \\ 
 &  & $D_3:$ the projection $\P^2 \times \P^1 \to \P^1$ & \\ \hline
\hyperref[t-pic2-E]{2-35}\label{table-2-35} & $56$ & $V_7$, i.e., $\P_{\P^2}(\MO_{\P^2} \oplus \MO_{\P^2}(1))$ \\
 & & $C_2:$ the projection $\mathbb{P}_{\P^2} (\MO_{\mathbb{P}^2}\oplus \MO_{\mathbb{P}^2}(1)) \to \P^2$ & \\ 
 &  & $E_2:$ blowup of $\P^3$ at a point &\\ \hline
\hyperref[t-pic2-E]{2-36}\label{table-2-36} & $62$  & 
$\mathbb{P}_{\P^2}(\MO_{\mathbb{P}^2}\oplus \MO_{\mathbb{P}^2}(2))$\\
& & 
$C_2:$ the projection $\mathbb{P}_{\P^2}(\MO_{\mathbb{P}^2}\oplus \MO_{\mathbb{P}^2}(2)) \to \P^2$  & \\ 
&  &  $E_5:$ blowup at the singular point of the cone over the Veronese surface & \\ \hline
      \end{longtable}
  \end{center}

\begin{proof}(of Table \ref{table-pic2}) 
We use Notation \ref{n-pic-2}. 
Almost all the parts follow from Section \ref{s-pic2}. 
It is enough to 
compute 
\begin{enumerate}
\item[(D)]$(-K_X)^2 \cdot X_t$ (type $D$), and 
\item[(C)] $\deg \Delta$ (type $C$). 
\end{enumerate}

{\cred 
(D) Assume that $f_1 : X \to \P^1$ is of type $D$ and let us see how to compute 
$(-K_X)^2 \cdot X_t$. 
Recall that $f_2$ is of type $E_1$ or $C$ (Theorem \ref{pic-over2}(2)). 
If $f_2$ is of type $E_1$, then 
$(-K_X)^2 \cdot X_t$ has been already comuted in Proposition \ref{p-E1D1}. 
Assume that $f_2$ is of type $C$. In this case,  
we have $H_1^3 =H_1^2 \cdot H_2 =0$ and $H_2^2 \equiv (2/\mu_2) \ell_2$ (Corollary \ref{fib_of_C}), 
which imply the following:  
\[
(-K_X)^2 \cdot X_t 
= (\mu_2 H_1 +\mu_1 H_2)^2 \cdot H_1 
= \mu_1^2 H_1 \cdot H_2^2 = \frac{2\mu_1^2}{\mu_2} H_1 \cdot \ell_2 = \frac{2\mu_1^2}{\mu_2}, 
\]
where the last equality holds by Proposition \ref{basis}(2). 
}

{\cred (C)} 
Assume that $f_1 : X \to \P^2$ is of type $C$ and let us see how to compute 
$\deg \Delta$ for $\Delta := \Delta_{f_1}$. 
We have $H_1^3 =0$ and $H_1^2 \cdot H_2 = (\frac{2}{\mu_1} \ell_1) \cdot H_2=\frac{2}{\mu_1}$ (Lemma \ref{pic2-numequiv}, Proposition \ref{basis}), 
which imply the following:  
\[
\deg \Delta \overset{{\rm (a)}}{=} 12 - (-K_X)^2 \cdot   H_1 \overset{{\rm (b)}}{=} 12 - (\mu_2 H_1 + \mu_1 H_2)^2 \cdot H_1 
=12 - 4\mu_2   -  \mu_1^2 H_1 \cdot H_2^2, 
\]
where (a) and (b) follow from 
Proposition \ref{conic_disc-num} and Proposition \ref{basis}, respectively. 
{\cred If $f_2$ is of type $C$ (resp. $D$), then 
we get $H_1 \cdot H_2^2= H_1 \cdot (\frac{2}{\mu_2} \ell_2)
=\frac{2}{\mu_2}$ (resp. $H_1 \cdot H_2^2 =0$).}


{\cred What is remaining is the case when 
$f_2$ is of type $E$. 
If $f_2$ is of type $E_1$, then $\deg \Delta$ has already been computed in 
Proposition \ref{p-E1C1}. 
Hence we may assume that $f_2$ is of type $E_2, E_3, E_4$, or $E_5$. 
By Theorem \ref{t-pic2-E}, 
the problem is reduced to the case when $X$ is of No.\ 2-8.} 
It holds that $-K_X \sim H_1 + H_2$ (Proposition \ref{basis}). For $D:=\Ex(f_2)$, we have $K_X \sim f_2^*K_Y +D$ and $( (-K_X)^2 \cdot D, (-K_X) \cdot D^2, D^3) = (2, -2, 2)$ \cite[Proposition 3.22]{TanII}.  
By Lemma \ref{r1-c*}(1), we have $D \cdot \ell_1 =1$ (note that $\zeta \equiv 2\ell_1$ for a fibre $\zeta$ of $f_1 : X \to \P^2$). 
Hence we obtain $-f_2^*K_Y \cdot \ell_1 = -K_X \cdot \ell_1 + D \cdot \ell_1 =2$. 
By $H_2 \cdot \ell_1 =1$ (Proposition \ref{basis}), we get $-f_2^*K_Y \equiv 2H_2$. Then 
\[
(-K_X)^2 \cdot H_1 = (-K_X)^2 \cdot (-K_X -H_2) = (-K_X)^2 \cdot \left(-K_X - \frac{-K_X+D}{2}\right) 
\]
\[
= \frac{(-K_X)^3 - (-K_X)^2 \cdot D}{2} = \frac{14-2}{2} =6.  
\]
Thus $\deg \Delta = 12 -
(-K_X)^2 \cdot H_1 = 12 - 6 =6$. 
\end{proof}

\bibliographystyle{skalpha}
\bibliography{reference.bib}

\end{document}